\RequirePackage{fix-cm}
\documentclass[twocolumn]{svjour3}          
\smartqed  
\usepackage{graphicx}
\usepackage{mathptmx}      
\usepackage{amsmath}
\usepackage{cite}
\usepackage{algorithm}
\usepackage{algpseudocode}
\usepackage{hyperref}
\hypersetup{%
    colorlinks=True, 
    citecolor=blue,
    linkcolor=black} 
\usepackage{verbatimbox}
    
\usepackage{xargs} 
\usepackage[colorinlistoftodos,prependcaption,textsize=tiny]{todonotes}
\usepackage{listings}
\definecolor{mygreen}{rgb}{0,0.6,0}
\definecolor{mygray}{rgb}{0.5,0.5,0.5}
\definecolor{mymauve}{rgb}{0.58,0,0.82}
\definecolor{mygreen}{rgb}{0,0.6,0}
\definecolor{mygray}{rgb}{0.5,0.5,0.5}
\definecolor{mymauve}{rgb}{0.58,0,0.82}
\newcommand{\eps}{{\varepsilon}}
\newcommand{\overbar}[1]{\mkern 1.5mu\overline{\mkern-1.5mu#1\mkern-1.5mu}\mkern 1.5mu}

\newcommand{\ben}{\begin{equation}}
\newcommand{\een}{\end{equation}}

\newcommand{\benn}{\begin{equation*}}
\newcommand{\eenn}{\end{equation*}}

\newcommand{\ngsolve}{\texttt{NGSolve}}
\newcommand{\normal}{n}
\newcommand{\tangDet}{\omega}

\newcommand{\VR}{{\mathbf{R}}}
\newcommand{\Id}{\mbox{Id}}                     
\newcommand{\curl}{\mbox{{\bf curl}}}  
\newcommand{\Dsf}{\mathsf{D}}
\newcommand{\Cj}{{\cal J}}
\newcommand{\VC}{{\mathbf{C}}}
\newcommand{\Cb}{{\cal B}}
\newcommand{\Div}{\mbox{div}}          
\newcommand{\Det}{\mbox{det}}                    

\begin{document}

\title{Fully and Semi-Automated Shape Differentiation in \ngsolve
}

\author{Peter Gangl \and Kevin Sturm \and Michael Neunteufel \and Joachim Sch\"oberl}

\authorrunning{P. Gangl, K. Sturm, M. Neunteufel, J. Sch\"oberl} 

\institute{P. Gangl (corresponding author) \at
              TU Graz, Steyrergasse 30, 8010 Graz, Austria \\
              \email{gangl(at)math.tugraz.at}           
           \and
           K. Sturm\at
              TU Wien, Wiedner Hauptstr. 8-10, 1040 Vienna\\
              \email{kevin.sturm(at)tuwien.ac.at}
           \and
           M. Neunteufel\at
              TU Wien, Wiedner Hauptstr. 8-10, 1040 Vienna\\
              \email{michael.neunteufel(at)tuwien.ac.at}
           \and
           J. Sch\"oberl\at
              TU Wien, Wiedner Hauptstr. 8-10, 1040 Vienna\\
              \email{joachim.schoeberl(at)tuwien.ac.at}
}

\date{Received: date / Accepted: date}

\maketitle

\begin{abstract}
In this paper we present a framework for automated shape differentiation in the finite element software \ngsolve. Our approach combines the mathematical Lagrangian approach for differentiating PDE constrained 
shape functions with the automated differentiation capabilities of \ngsolve.  
The user can decide which degree of automatisation is required, thus allowing for either a more custom-like or black-box-like behaviour of the software. 
We discuss the automatic generation of first and second order shape derivatives for unconstrained model problems as well as for more realistic problems that are constrained by different types of partial differential equations. We consider linear as well as nonlinear problems and also problems which are posed on surfaces. 
In numerical experiments we verify the accuracy of the computed derivatives via a Taylor test. Finally we present first and second order shape optimisation algorithms and illustrate them for several numerical optimisation examples ranging from nonlinear elasticity to Maxwell's equations. 

\end{abstract}

\section{Introduction}

 Numerical simulation and shape optimisation tools to solve the problems have become an integral part in the design process of many products. Starting out from an initial design, non-parametric shape optimisation techniques based on first and second order shape derivatives can assist in finding shapes of a product which are optimal with respect to a given objective function. Examples include the optimal design of aircrafts \cite{a_SCILSCGA_2013, a_SCILSCGA_2011a}, optimal inductor design \cite{a_HOSO_2003a}, optimisation of microlenses \cite{a_PASAHIHA_2015a}, the optimal design of electric motors \cite{GLLMS2015}, applications to mechanical engineering \cite{AllaireJouveToader2004, a_LA_2018a}, multiphysics problems \cite{Feppon2019Sep} or electrical impedance tomography (EIT) in medical sciences to name only a few \cite{a_HILA_2008a}.

 Shape optimisation algorithms are based on the concept of shape derivatives. Let $\mathcal P(\VR^d)$ denote the set of all subsets of $\VR^d$. Further let $\mathcal A \subset \mathcal P(\VR^d)$ be a set of admissible shapes and $\mathcal J: \mathcal A \rightarrow \VR$ be a shape function. Given an admissible shape $\Omega \in \mathcal A$ and a sufficiently smooth vector field $V$, we define the perturbed domain $\Omega_t := (\Id+t V)(\Omega)$ for a small perturbation parameter $t >0$. The shape derivative is defined as 
\begin{align} \label{eq_def_shapeDer}
    D \mathcal J(\Omega)(V) := \left. \left(\frac{d}{dt} \mathcal J(\Omega_t) \right)\right\rvert_{t=0} = \underset{t \searrow 0}{\mbox{lim }} \frac{\mathcal J(\Omega_t) - \mathcal J(\Omega)}{t}.
\end{align}

\begin{remark}
    We remark that a frequently used definition of shape differentiability is to require the mapping $V\mapsto J((\Id+V)(\Omega))$ being Fr{\'e}chet differentiable in $V=0$; see \cite{b_AL_2007a,b_HEPI_2005a,a_MUSI_1976a}. This stronger notion of differentiability implies that the limit defined in \eqref{eq_def_shapeDer} exists.
\end{remark}

In most practically relevant applications, the objective functional depends on the shape of a (sub-)domain via the solution to a partial differential equation (PDE). Thus, one is facing a problem of PDE-constrained shape optimisation of the form
\begin{equation}\label{P:intro}
    \begin{split}
    \underset{(\Omega,u) \in \mathcal A\times Y }{\mbox{min }}& J(\Omega, u) \\
    \mbox{s.t. } (\Omega, u) \in \mathcal A\times  Y : e(\Omega; u,v) &= 0 \quad  \text{ for all } v \in  Y .
\end{split}
\end{equation}
Here, the second line represents the constraining boundary value problem posed on a Hilbert space $Y $, which we assume to be uniquely solvable for all admissible $\Omega \in \mathcal A$. Denoting the unique solution for a given $\Omega \in \mathcal A$ by $u_\Omega$, we introduce the notation for the reduced functional 
\begin{align*}
    \mathcal J(\Omega) := J(\Omega,u_\Omega).
\end{align*}
In order to be able to apply a shape optimisation algorithm to a given problem of this kind, the shape derivative \eqref{eq_def_shapeDer} has to be computed, see the standard literature \cite{DZ2, SZ} or \cite{a_ST_2015a} for an overview of different approaches. In the following we focus on computing the so-called volume form of the shape derivative which in a finite element context is known to give a better approximation compared to the boundary form; see \cite{MR3348199, MR2642680}.

The convergence of shape optimisation algorithms can be speeded up by using second order shape derivatives. Given two sufficiently smooth vector fields $V$, $W$ and an admissible shape $\Omega \in\mathcal A$, let $\Omega_{s,t} := (\Id +s V + t W)(\Omega)$ be the perturbed domain. Then, the second order shape derivative is defined as
\begin{align} \label{eq_def_shapeDer2}
    D^2 \mathcal J(\Omega)(V)(W) := \left. \left(\frac{d^2}{dsdt} \mathcal J(\Omega_{s,t}) \right)\right\rvert_{s,t=0}. 
\end{align} 
Second order information in Newton-type algorithms has been explored in the articles \cite{a_NORO_2002a,a_ALCAVI_2016a,a_PAST_2019a,a_EPHASC_2007a,a_VSC_2014a}. Since the computation of second order shape derivatives is more involved and error prone, several authors have employed automatic differentiation (AD) tools, see e.g. \cite{a_SC_2018a} and \cite{a_HAMIPAWE_2019a} for two approaches based on the Unified Form Language (UFL) \cite{a_ALMALOOLROWE_2014a}. In \cite{a_HAMIPAWE_2019a}, the authors present a fully automated shape differentiation software which uses the transformation properties on the finite element level. In \cite{a_SC_2018a} (see also the earlier work \cite{a_SC_2014a}) the automated derivatives are computed using UFL. The strategies of \cite{a_HAMIPAWE_2019a} and \cite{a_SC_2018a} differ in that, for the latter, the software computes an unsymmetric shape Hessian since it involves the term $D \Cj(\Omega)(\partial VW)$. Optionally the software allows to make the shape Hessian symmetric by requiring $\partial VW=0$. We will discuss the subtle difference and the relation between the two possible ways of defining shape Hessians in Remark \ref{rem_shapeHess} of Section \ref{sec_secondOrder}. Let us also mention \cite{a_DOMISEPU_2020a} where automated shape derivatives for transient PDEs in FEniCS and Firedrake are presented.

In this paper we present an alternative framework for AD of PDE constrained problems of 
type \eqref{P:intro}. There exist several approaches for the rigorous derivation of the shape derivative of PDE-constrained shape functionals, see \cite{c_ST_2015a} for an overview. The main idea, however, is always similar. After transforming the perturbed setting back to the original domain, shape differentiation in the direction of a given vector field reduces to the differentiation with respect to the scalar parameter $t$ which now enters via the corresponding transformation and its gradient. It is shown in \cite{a_ST_2015a} that the shape derivative for a nonlinear PDE-constrained shape optimisation problem can be computed as the derivative of the Lagrangian with respect to the perturbation parameter. We will illustrate this systematic procedure for a number of different applications and utilise symbolic differentiation provided by the finite element software package \ngsolve{} \cite{Schoeberl2014} to obtain the shape derivative for different classes of PDE-constrained optimisation problems. \ngsolve{} allows for the fast and efficient numerical solution of a large number of different boundary value problems. The aim of this paper is to extend \ngsolve{} by the possibility of semi-automatic and fully automatic shape differentiation and optimisation.

Distinctly from previous approaches we cover the following two points:
\begin{itemize}
    \item a fully automated setting requiring as input the weak formulation of the constraint and the cost function,
    \item a semi-automated setting which offers a highly customisable user interface, but requires mathematical background knowledge.
\end{itemize}

\paragraph{Structure of the paper.}
In Section~\ref{sec_intro_ng} we give a brief introduction on how to solve a PDE in \ngsolve{} and present its built-in auto-differentiation capabilities. The introduced syntax will also lay the foundation for the following sections.
In Section~\ref{sec_intf} we present a first unconstrained shape optimisation problem and show how to solve it in \ngsolve{}. For this purpose we show how to compute the first and second order shape derivative in a semi-automated way. Section~\ref{sec_PDE_constr_shape} extends the preceding section by incorporating a PDE constraint. The strategy is illustrated by means of a simple Poisson equation. We also show how to treat the computation of shape derivatives when the PDE is defined on surfaces. 
While the semi-automated shape differentiation presented in Sections~\ref{sec_intf} and \ref{sec_PDE_constr_shape} requires mathematical background knowledge, in Section~\ref{sec_fullyAutomated} we show how the shape derivatives can be computed in a fully automated fashion. In the last section of the paper we verify the computed formulas by a Taylor test, discuss optimisation algorithms and present several numerical optimisation examples including nonlinear elasticity, Maxwell's equations and Helmholtz's equation.

\section{A brief introduction to \ngsolve{}}\label{sec_intro_ng}
In this section, we give a brief overview of the main concepts of the finite element software \ngsolve{} \cite{Schoeberl2014}. We first describe the main principles for numerically solving boundary value problems in \ngsolve{} before focusing on its built-in automatic differentiation capabilities. In the subsequent sections of this paper, these ingredients will be combined to implement the shape derivative of unconstrained and PDE-constrained shape optimisation problems in an automated way.

\subsection{Solving PDEs with finite elements in \ngsolve{}} \label{sec_intro_ngsolve}
In this section, we illustrate the syntax of \ngsolve{} using the \texttt{python} programming language for the Poisson equation with homogeneous Dirichlet conditions as a model problem. We refer the reader to the online documentation 
\begin{center}
    \href{https://ngsolve.org/docu/latest/}{https://ngsolve.org/docu/latest/}
\end{center}
for a more detailed description of the many features of this package.

Given a domain $\Omega \subset \VR^d$ and a right hand side $f$, we consider the model problem to find $u$ satisfying
\begin{align*}
    - \Delta u  = f &\qquad \mbox{in } \Omega, \\
    u=0 &\qquad\mbox{on } \partial \Omega.
\end{align*}
The weak form of the model problem reads
\ben \label{eq_model_poisson}
    \mbox{Find } u \in H_0^1(\Omega): \; \int_\Omega \nabla u \cdot \nabla w \; \mbox dx = \int_\Omega f w \; \mbox dx \quad \forall w \in H_0^1(\Omega).
\een
We consider a ball of radius $\frac{1}{2}$ in two space dimensions centered at the point $(0.5, 0.5)^\top$, i.e. $\Omega = B((0.5, 0.5)^\top,0.5)$, and the right hand side is defined by $f(x_1, x_2) = 2 x_2(1-x_2)+2x_1(1-x_1)$. We will go through the steps for numerically solving this problem by the finite element method.

We begin by importing the necessary functionalities and setting up a finite element mesh.
\begin{lstlisting}[firstnumber = 1]
from ngsolve import *
from netgen.geom2d import SplineGeometry

geo = SplineGeometry()
geo.AddCircle((0.5,0.5),0.5,bc="circle") (*@\label{lst_defGeo}@*)

mesh = Mesh(geo.GenerateMesh(maxh=0.2)) (*@\label{lst_genMesh}@*)
mesh.Curve(3)   (*@\label{lst_genMesh2}@*)
\end{lstlisting}
The first line imports all modules from the package \ngsolve{}. The second line includes the SplineGeometry function which enables  us to define a mesh via a geometric description, in our case a circle centered at $(0.5,0.5)^\top$ of radius $0.5$. Finally the mesh is generated in line \ref{lst_genMesh} and in line \ref{lst_genMesh2} we specify that we want to use a curved finite element mesh for a more accurate approximation of the geometry. For that purpose, a projection-based interpolation procedure is used, see e.g. \cite{demk04}.

\begin{figure}
    \begin{center}
        \includegraphics[width=.45\textwidth]{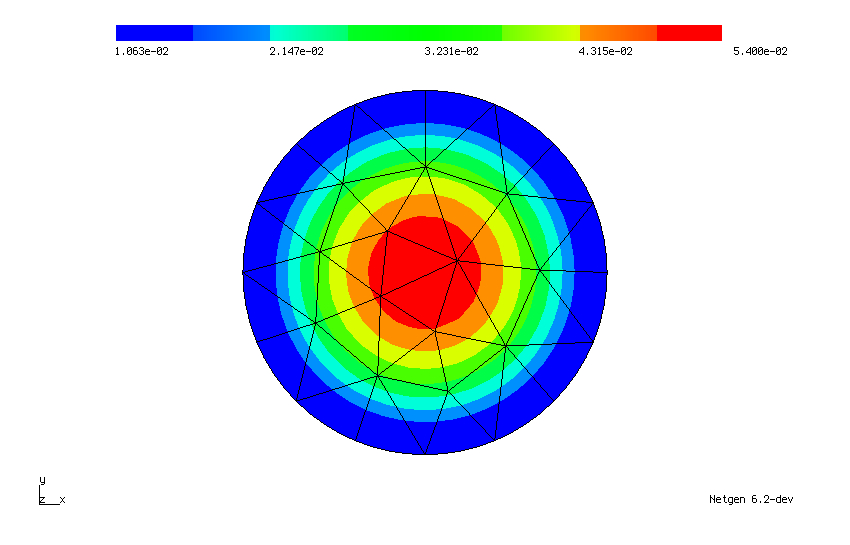}
    \end{center}
    \caption{Solution of problem \eqref{eq_model_poisson} by code fragments of Section \ref{sec_intro_ngsolve} with 29 nodes, 40 (curved) triangular elements and polynomial order $3$.}
    \label{fig_solModelPoisson}
\end{figure}

Next in line \ref{lst_genFes} we define an $H^1$ conforming finite element space of polynomial degree $3$ and include Dirichlet boundary conditions on the boundary of the domain $\partial \Omega$ (referenced by the string \texttt{``circle''} that we assigned in line \ref{lst_defGeo}). On this space we define a trial function \texttt{u}  in line \ref{lst_defTrial} and a test function \texttt{w} in line \ref{lst_defTest}. These are purely symbolic objects which are used to define boundary value problems in weak form.
\begin{lstlisting}[firstnumber=last]
fes = H1(mesh, order=3, dirichlet="circle") (*@\label{lst_genFes}@*)

u = fes.TrialFunction()     (*@\label{lst_defTrial}@*)
w = fes.TestFunction()      (*@\label{lst_defTest}@*)
\end{lstlisting}
For a more compact presentation later on, we define a coefficient function \texttt{X} which combines the three spatial components:
\begin{lstlisting}[firstnumber=last]
X = CoefficientFunction((x,y,z)) (*@\label{lst_defX}@*)
\end{lstlisting}
Now, the left and right hand sides of problem \eqref{eq_model_poisson} can be conveniently defined as a bilinear or linear form, respectively, on the finite element space \texttt{fes} by the following lines.
\begin{lstlisting}[firstnumber=last]
L = LinearForm(fes)
f1 = (2*X[1]*(1-X[1])+2*X[0]*(1-X[0]))
L += f1 * w * dx

a = BilinearForm(fes, symmetric=True)
a += grad(u)*grad(w)*dx
\end{lstlisting}
We assemble the system matrix coming from the bilinear form \texttt{a} and the load vector coming from \texttt{L} and solve the corresponding system of linear equations.
\begin{lstlisting}[firstnumber=last]
a.Assemble()
L.Assemble()

gfu = GridFunction(fes)
gfu.vec.data = a.mat.Inverse(fes.FreeDofs(), inverse="sparsecholesky") * L.vec

Draw (gfu, mesh, "state")
\end{lstlisting}
Here, \texttt{gfu} is defined as a \texttt{GridFunction} over the finite element space \texttt{fes}. A \texttt{GridFunction} object is used to save the results by containing the corresponding finite element coefficient vectors. Further, it can evaluate the stored finite element solution at a given mesh point. The Dirichlet conditions are incorporated into the direct solution of the linear system and the numerical solution is drawn in the graphical user interface. The numerical solution is depicted in Figure~\ref{fig_solModelPoisson}.

\subsection{Automatic Differentiation in \ngsolve{}} \label{sec_autoDiff}
\begin{figure}
    \begin{center}
        \begin{tabular}{cc}
            \includegraphics[width=.2\textwidth]{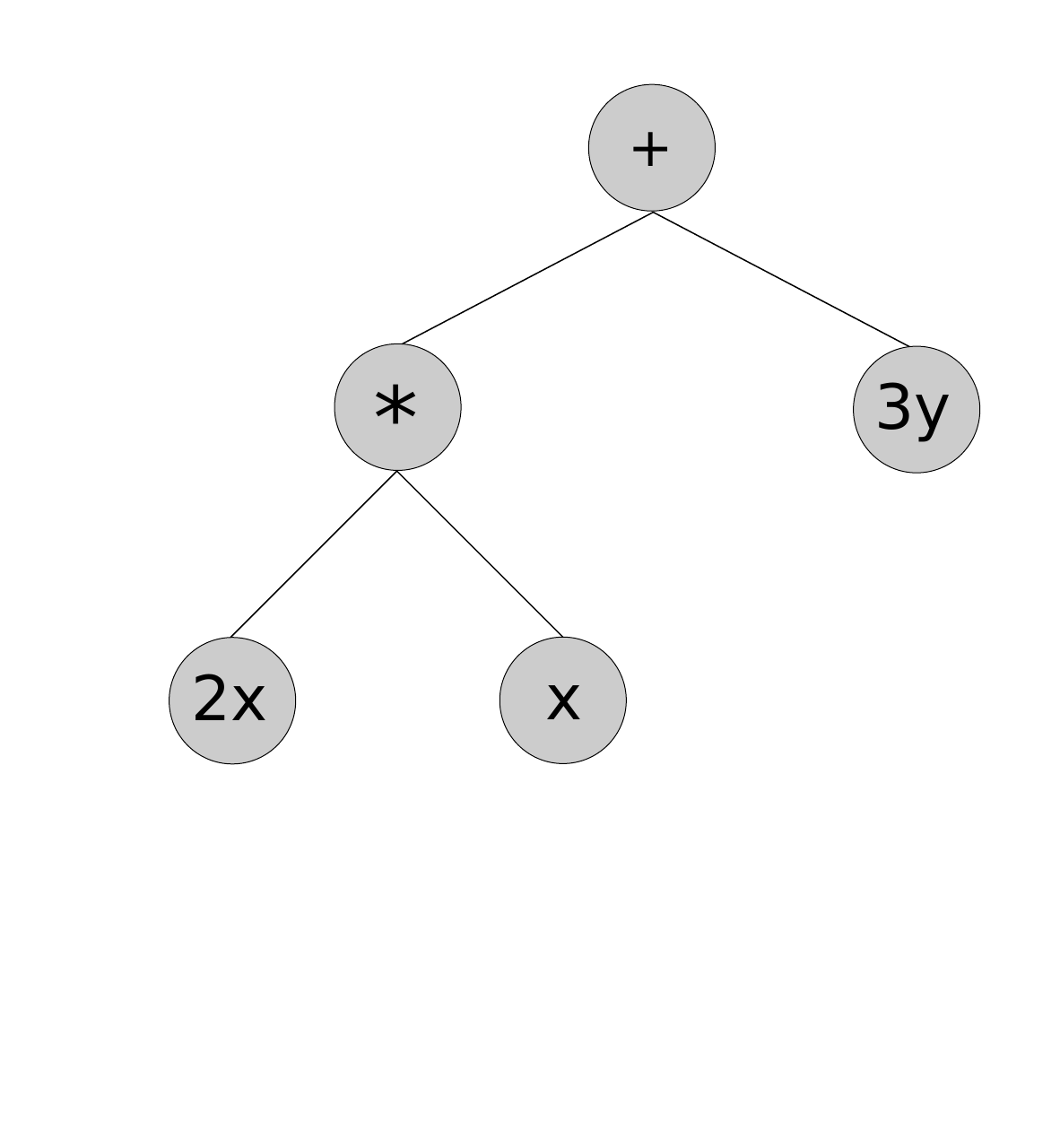} & \includegraphics[width=.2\textwidth]{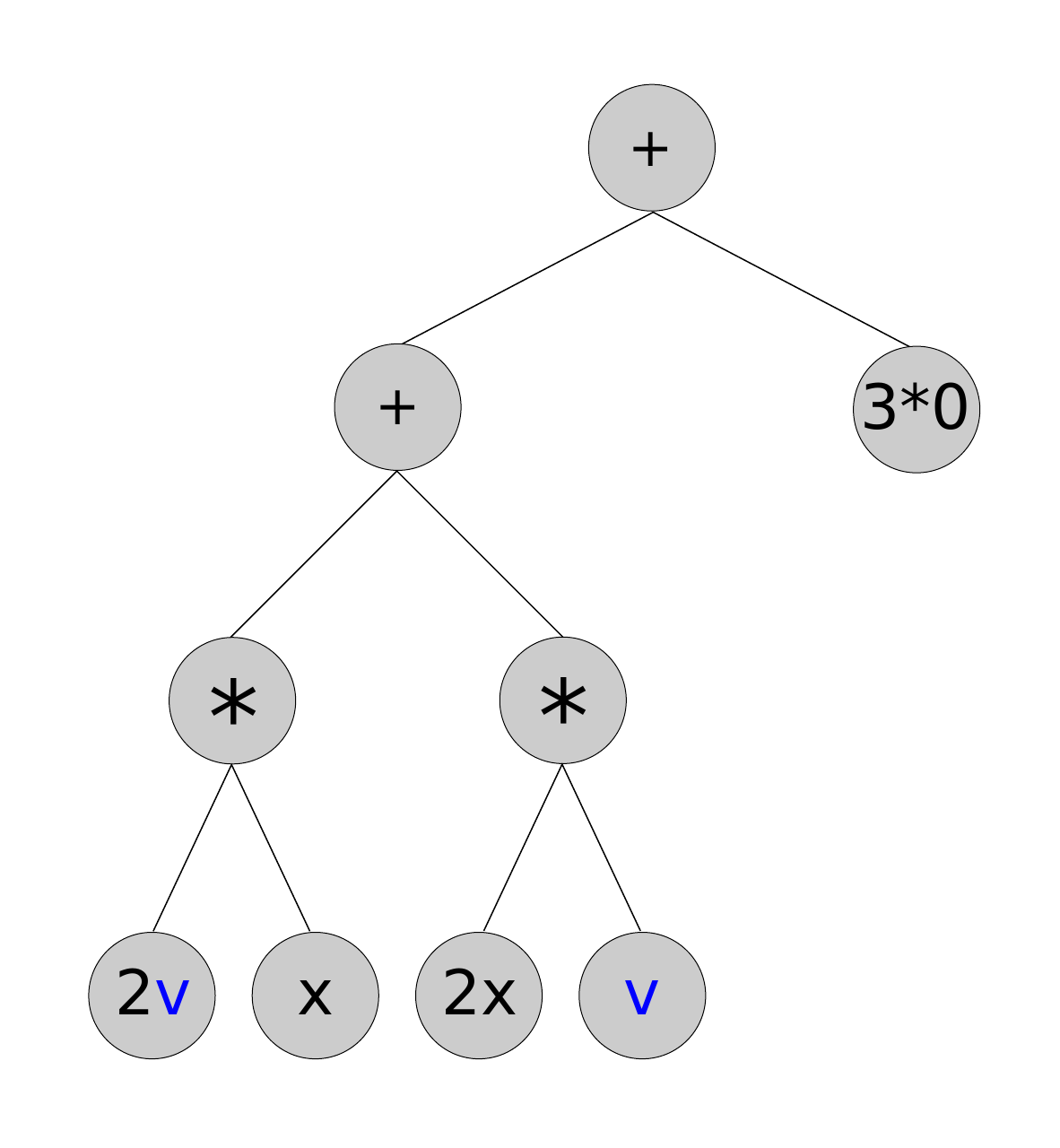} \\
            (a) & (b)
        \end{tabular}
    \end{center}
    \caption{Illustration of \texttt{Diff} command for example \textit{expr= 2x*x+3y}. (a) Expression tree for \textit{expr}. (b) Expression tree for expression obtained by call of {\texttt{expr.Diff(x, v)}}.}
    \label{fig_expressionTreeExample}
\end{figure}

In \ngsolve, symbolic expressions are stored in expression trees, see Figure~\ref{fig_expressionTreeExample} for an example. It is possible to differentiate an expression \textit{expr} with respect to a variable \textit{var} appearing in \textit{expr} into a direction \textit{dir} by the command
\begin{center}
    \texttt{expr.Diff(var, dir)}.
\end{center}
Mathematically this line corresponds to the directional derivative of g:=\emph{expr} at $x:=var$ in direction $v:=dir$, that is,
\ben
Dg(x)(v).
\een
When calling the \texttt{Diff} command for \texttt{expr}, the expression tree of \textit{expr} is gone through node by node, and for each node the corresponding differentiation rules such as product rule or chain rule are applied. When a node represents the variable with respect to which the differentiation is carried out, it is replaced by the direction \textit{dir} of differentiation.

Figure \ref{fig_expressionTreeExample} shows the differentiation of the expression  \textit{expr= 2x*x+3y} with respect to \textit{x} into the direction given by \textit{v}:

\begin{lstlisting}[firstnumber=last]
v = Parameter(1)
expr = 2*x*x+3*y
dexpr = expr.Diff(x,v)
print(expr)
print(dexpr)
\end{lstlisting}
The output of \texttt{print(expr)} reads
\begin{verbnobox}[\footnotesize]
coef binary operation '+', real
  coef binary operation '*', real
    coef scale 2, real
      coef coordinate x, real
    coef coordinate x, real
  coef scale 3, real
    coef coordinate y, real
\end{verbnobox}
which translates to $2x*x+3y$ and corresponds to the expression tree depicted in Figure \ref{fig_expressionTreeExample}(a). The output of \texttt{print(dexpr)} reads
\begin{verbnobox}[\footnotesize]
coef binary operation '+', real
  coef binary operation '+', real
    coef binary operation '*', real
      coef scale 2, real
        coef N5ngfem28ParameterCoefficientFunctionE, real
      coef coordinate x, real
    coef binary operation '*', real
      coef scale 2, real
        coef coordinate x, real
      coef N5ngfem28ParameterCoefficientFunctionE, real
  coef scale 3, real
    coef 0, real
\end{verbnobox}
which translates to $(2v*x + 2x*v) +3*0$ and corresponds to the expression tree depicted in Figure \ref{fig_expressionTreeExample}(b). The coefficient $$ \mathrm{\texttt{N5ngfem28ParameterCoefficientFunctionE}}$$ appearing therein is the C++ internal class name of the Python object \texttt{Parameter}.

\ngsolve{} trial and test functions are purely symbolic objects used for defining bilinear and linear forms. Therefore, they do not depend on the spatial variables $x$, $y$, $z$ as can be seen by differentiating them. \ngsolve{} \texttt{GridFunction}s on the other hand represent functions in the finite element space. However, also for these objects, the space dependency is omitted when performing symbolic differentiation. The code segments

\begin{lstlisting}[firstnumber=last]
u = fes.TrialFunction()  # symbolic object
w = fes.TestFunction()   # symbolic object
gf = GridFunction(fes)
gf.Set(x*x*y)

print("Diff u w.r.t. x", u.Diff(x))
print("Diff w w.r.t. x", w.Diff(x))
print("Diff gf w.r.t. x", gf.Diff(x))
\end{lstlisting}
will give the following output:
\begin{verbnobox}[\footnotesize]
Diff u w.r.t. x:  ConstantCF, val = 0
Diff w w.r.t. x:  ConstantCF, val = 0
Diff gf w.r.t. x:  ConstantCF, val = 0
\end{verbnobox}
Here, the \texttt{GridFunction.Set} method takes a \texttt{CoefficientFunction} object and performs a (local) $L^2$ best-approximation into the underlying finite element space with respect to its natural norm and stores the resulting coefficient vector. 

\section{Semi-automatic shape differentiation without constraints}
\label{sec_intf}
We will illustrate the steps to be taken in order to obtain the shape derivative of a shape function in a semi-automatic way for a simple shape optimisation problem.  For $\Omega \subset \VR^d$ bounded and open and a continuously differentiable function $f \in C^1(\VR^d)$, we consider the shape differentiation of the shape function
\begin{align}\label{eq:shape_func_f}
    \mathcal J(\Omega) = \int_\Omega f(x) \, \mbox dx.
\end{align}
Clearly the minimiser of $\mathcal J$ over all measurable sets in $\VR^d$ is given by $\Omega^*=\{x\in \VR^d:\; f(x)<0\}$. We also refer to \cite{ubt_epub3251} for the computations of first and 
second order variations of functions of type \eqref{eq:shape_func_f} where $\Omega$ is a submanifold of
$\VR^d$.

\subsection{First order shape derivative} \label{sec_firstordershape}

Henceforth we denote by $C^{0,1}(\VR^d)^d$ the space of bounded and Lipschitz continuous vector fields $V:\VR^d\to\VR^d$. 
In view of Rademachers' theorem \cite[Thm.6, p.296]{b_EV_2010a} the space $C^{0,1}(\VR^d)^d$ corresponds to the Sobolev space $W^{1,\infty}(\VR^d)^d$. 

Given a vector field $V \in  C^{0,1}(\VR^d)^d$, we define the transformation
\begin{align*}  
    T_t(x) := (\Id + t \, V)(x) , \quad x \in \VR^d,\quad t\ge 0.
\end{align*}
\begin{definition}
    The first order shape derivative of a shape function $\mathcal J$ at $\Omega$ in direction $V\in C^{0,1}(\VR^d)^d$ is defined by 
    \ben
    D\mathcal J(\Omega)(V) = \lim_{t \searrow  0} \frac{\mathcal J(T_t(\Omega)) - \mathcal J(\Omega)}{t}. 
    \een
\end{definition}

\subsubsection{Shape differentiation of unconstrained volume integrals} 
Using the transformation $y=T_t(x)$ and the notation $F_t := \partial T_t = I + t \partial V$ for the Jacobian of the transformation $T_t$, we get for $\mathcal J$ as in \eqref{eq:shape_func_f},
\begin{align}\label{eq:cost_int_f_trafo}
\mathcal J(\Omega_t) = \int_{\Omega_t} f(x') \; \mbox dx' = \int_\Omega (f\circ T_t)(x) \, \mbox{det}(F_t(x)) \,\mbox{d}x.
\end{align}

Now let us explain how to compute the shape derivative of $\Cj$. Denoting 
\ben \label{eq_defG}
G(T_t,F_t) :=  \int_\Omega (f\circ T_t)(x) \,\mbox{det}(F_t(x)) \,\mbox{d}x,
\een
the chain rule gives (formally)
\begin{align}
   \left. \frac{d}{dt} \mathcal J(\Omega_t) \right|_{t=0}
   &= \left.\frac{d}{dt} G(T_t,F_t) \right|_{t=0} \nonumber \\
   &= \left. \left( \frac{dG}{dT_t} \frac{d T_t}{dt} + \frac{dG}{dF_t} \frac{d F_t}{dt}\right)\right|_{t=0}. \label{eq_dGdTdGdF}
\end{align}
Using that $\frac{d T_t}{dt}(x) = V(x)$ and $\frac{d F_t}{dt}(x) = \partial V(x)$, we get for the shape derivative
\begin{align*}  
   D \Cj(\Omega)(V) =  \left. \frac{d}{dt} \mathcal J(\Omega_t) \right|_{t=0} = \left. \left( \frac{dG}{dT_t} V + \frac{dG}{dF_t} \partial V\right)\right|_{t=0}.
\end{align*}
This is the form we use for defining the first order shape derivative in \ngsolve. Note that a Lipschitz vector field is differentiable almost everywhere and hence $\partial V(x)$ is defined almost everywhere and bounded.

Given the function 
$f(x_1, x_2) = (x_1 - 0.5)^2/a^2 + (x_2 - 0.5)^2/b^2 - R^2$ with $a=1.3$, $b=1/a$ and $R=0.5$, we implement the transformed cost function \eqref{eq:cost_int_f_trafo} as follows:
\begin{lstlisting}[firstnumber=last,language=Python]
f = ((X[0]-0.5)/1.3)**2+(1.3*(X[1]-0.5))**2 - 0.5**2

F = Id(2)               # symbolic identity matrix (*@\label{lst_defF}@*)
G_f = f * Det(F) * dx   # F only acts as a dummy variable (*@\label{lst_defG}@*)
\end{lstlisting}
Here, we introduce the symbol \texttt{F} and assign to it the value of the identity matrix in line \ref{lst_defF}. This allows us to differentiate with respect to \texttt{F}. Then we define the function $G$ of \eqref{eq_defG} in line \ref{lst_defG}. The shape derivative is a bounded linear functional on a space of vector fields. We introduce a vector-valued finite element space \texttt{VEC} and define the object representing the shape derivative \texttt{dJOmega\_f} as a linear functional on \texttt{VEC}. In line \ref{lst_diff1}, we differentiate with respect to the spatial variables in the direction given by \texttt{V}. Note that \texttt{X} is the coefficient function we introduced in line \ref{lst_defX}. In line \ref{lst_diff2}, we deal with the differentiation with respect to \texttt{F}.

\begin{remark}
    Defining $\xi_t := \mbox{det}(F_t)$ and using $\frac{d}{dt}\xi_t|_{t=0} = \Div V$, it holds
    \begin{align*}
       \left. \frac{dG}{d F_t} \frac{d F_t}{dt} \right|_{t=0} &= \left. \frac{d G}{d \xi_t} \frac{d \xi_t}{d F_t} \frac{d F_t}{dt} \right|_{t=0} = \left. \frac{d G}{d \xi_t} \frac{d \xi_t}{dt }  \right|_{t=0} \\
       &= \left.\frac{d G}{d \xi_t} \mbox{div}V \right|_{t=0} = \int_\Omega f \Div V \, \mbox dx.
    \end{align*}
    Therefore, we obtain for the first order shape derivative the well-known formula
    \begin{align} \label{eq_dJ_intf}
        D \mathcal J(\Omega)(V) = \int_\Omega \nabla f \cdot V + f \, \mbox{div} V \, \mbox dx.
    \end{align}
 Finally if $\Omega$ is smooth enough (for instance $C^1$), it follows by integration by parts in \eqref{eq_dJ_intf} that the shape derivative is given by 
\ben
D\mathcal J(\Omega)(V) = \int_{\partial \Omega} f V\cdot \normal \;\mbox ds,
\een
where $\normal$ denotes the outward pointing normal along $\partial \Omega$.
\end{remark}
\begin{lstlisting}[firstnumber=last]
VEC = VectorH1(mesh, order=1, dirichlet="") #vectorial FE space of order 1  (*@\label{lst_def_VEC}@*)
V = VEC.TestFunction()

dJOmega_f = LinearForm(VEC)
dJOmega_f += G_f.Diff(X, V)  (*@\label{lst_diff1}@*)
dJOmega_f += G_f.Diff(F, grad(V)) (*@\label{lst_diff2}@*)
\end{lstlisting}

\subsubsection{Shape differentiation of unconstrained boundary integrals} 
For $\Omega$ and $f$ as in the previous section we consider
\ben \label{eq:shape_func_f_bnd}
    \mathcal J_{bnd}(\Omega) = \int_{\partial \Omega} f(x) \; \mbox dx.
\een
Then we get
\begin{align}
\mathcal J_{bnd}(\Omega_t) &= \int_{\partial \Omega_t} f(x') \; \mbox ds_{x'} \\
& = \int_{\partial \Omega} (f\circ T_t)(x) \, \mbox{det}(F_t(x)) | F_t(x)^{-\top} \normal(x) | \,\mbox ds_x, \label{eq:cost_int_f_trafo_bnd}
\end{align}
see e.g. \cite[Prop. 2.47]{SZ}, with the outer unit normal vector $\normal$ and $|\cdot|$ denoting the Euclidean norm.
It is shown in \cite[Prop. 2.50]{SZ} that the shape derivative of \eqref{eq:shape_func_f_bnd} is given by
\begin{align*}
    D \mathcal J_{bnd}(\Omega)(V) = \int_{\partial \Omega} \nabla f \cdot V + f (\mbox{div}V - n^\top \partial V n) ds_x.
\end{align*}

Again, we can compute the shape derivative in \ngsolve{} as the total derivative of expression \eqref{eq:cost_int_f_trafo_bnd} with respect to the parameter $t$. In \ngsolve, the only difference lies in the necessity to use the trace of the gradient of a test vector field \texttt{V}.
\begin{lstlisting}[firstnumber=last]
G_f_bnd = f * Det(F) * Norm( Inv(F).trans*specialcf.normal(2) ) * ds

dJOmega_f_bnd = LinearForm(VEC)
dJOmega_f_bnd+=G_f_bnd.Diff(X, V)                             # no trace needed
dJOmega_f_bnd+=G_f_bnd.Diff(F,grad(V).Trace())                # trace needed
\end{lstlisting}

Note that the trace operator for gradients on the boundary is obligatory in NGSolve, whereas for direct evaluation of $H^1$ trial and test functions itself it is optional.

\subsection{Second order shape derivatives} \label{sec_secondOrder}
For second order shape derivatives, we consider perturbations of the form
\begin{align*}
    T_{s,t}(x) = ( \Id + s V + t W)(x), \quad x \in \VR^d,
\end{align*}
for $s, t \geq 0$ and define $\Omega_{s,t} := T_{s,t}(\Omega)$.
\begin{definition}
The second order shape derivative of a shape function $\mathcal J$ at $\Omega$ in direction $(V,W)\in C^{0,1}(\VR^d)^d \times C^{0,1}(\VR^d)^d$ is defined by 
\ben\label{eq:second_derivative}
D^2\mathcal J(\Omega)(V)(W) =     \left. \frac{d^2}{dsdt} \mathcal J(\Omega_{s,t}) \right\rvert_{s=t=0}.
\een
\end{definition}
\begin{remark} \label{rem_shapeHess}
We remark that if $\mathcal J$ is smooth enough, the second order derivative as defined in \eqref{eq:second_derivative} is symmetric by definition:
\ben
D^2\Cj(\Omega)(V)(W) = D^2\Cj (\Omega)(W)(V).
\een
We stress that this derivative is not the same as the shape derivative obtained by repeated 
shape differentiation, that is, it does not coincide with (see, e.g., \cite[Chap. 9, Sec. 6]{b_DEZO_2011a})
\ben \label{eq_d2J_nonsymm}
d^2\Cj(\Omega)(V)(W) := \lim_{t  \searrow  0}\frac{D\Cj(T_t^W(\Omega))(V)- D\Cj(\Omega)(V)}{t}
\een
which is in general asymmetric. 

The derivative defined in \eqref{eq_d2J_nonsymm} is only symmetric if $\partial VW=0$ since it holds
\ben
d^2\Cj(\Omega)(V)(W) = D^2\Cj(\Omega)(V)(W) + D\Cj(\Omega)(\partial VW),
\een
see also the early work of Simon \cite{Simon1989} on this topic. 
However, in \ngsolve{}, when repeating the shape differentiation procedure introduced in Section \ref{sec_firstordershape}, we compute directly the second order shape derivative as defined in \eqref{eq:second_derivative}. Here, we exploit the fact that trial functions are independent of the spatial coordinates, see also Section \ref{sec_autoDiff} and the example below. 
\end{remark}

Let us now exemplify the computation of the second order shape derivative for the shape function $\mathcal J$ defined in \eqref{eq:shape_func_f}. Similarly to the computations of the first derivative, we use the notation $F_{s,t} := \partial T_{s,t} = I + s \partial V + t \partial W$. 
Then we get 
\begin{align*}
    \frac{d^2}{dsdt} \mathcal J(\Omega_{s,t})& \left. \vphantom{\frac{d^2}{dsdt}} \right\rvert_{s=t=0} = \left. \frac{d^2}{dsdt} \int_{\Omega_{s,t}} f(x) \,\mbox{d}x \right\rvert_{s=t=0} \\
    &=\left. \frac{d^2}{dsdt} \int_\Omega (f\circ T_{s,t})(x) \, \mbox{det}(F_{s,t}(x)) \,\mbox{d}x \right\rvert_{s=t=0}.
\end{align*}
Again, using the notation $$G(T_{s,t},F_{s,t}) = \int_\Omega (f\circ T_{s,t})(x) \, \mbox{det}(F_{s,t}(x)) \,\mbox{d}x,$$ we get
\begin{align*}
    \left. \frac{d^2}{dsdt} \mathcal J(\Omega_{s,t}) \right\rvert_{s=t=0}
    &=\left.  \frac{d^2}{dsdt}  G(T_{s,t},F_{s,t}) \right\rvert_{s=t=0}\\
    &=\left.  \frac{d}{ds} \left( \frac{dG}{dT_{s,t}} \frac{d T_{s,t}}{dt} + \frac{dG}{dF_{s,t}} \frac{d F_{s,t}}{dt}   \right) \right\rvert_{s=t=0}.
\end{align*}
Using that $\frac{d^2 T_{s,t} }{dsdt} = 0$ and $\frac{d^2 F_{s,t} }{dsdt} = 0$, we get further
\begin{align}
     \frac{d^2}{dsdt}& \mathcal J(\Omega_{s,t}) \left. \vphantom{\frac{d^2}{dsdt}} \right\rvert_{s=t=0} \nonumber \\ 
    =&\left.  \frac{d}{ds} \left(\frac{dG}{dT_{s,t}} \right) \frac{d T_{s,t}}{dt} + \frac{d}{ds}\left(\frac{dG}{dF_{s,t}}  \right) \frac{d F_{s,t}}{dt}  \right\rvert_{s=t=0} \nonumber\\                                                                         =&  \left( \frac{d^2 G}{d T_{s,t}^2} \frac{d T_{s,t}}{ds} + \frac{d^2 G}{d F_{s,t} dT_{s,t}} \frac{d F_{s,t} }{d s}  \right) \frac{d T_{s,t}}{dt} \nonumber \\
    &+  \left(  \frac{d^2 G}{d T_{s,t} dF_{s,t}} \frac{d T_{s,t} }{d s}  + \frac{d^2 G}{d F_{s,t}^2} \frac{d F_{s,t}}{ds} \right) \frac{d F_{s,t}}{dt} \left.  \vphantom{\frac{d^2}{dsdt}} \right\rvert_{s=t=0}.  \label{eq_d2J_general}
\end{align}
Formula \eqref{eq_d2J_general} is used for the automatic derivation of the second order shape derivative in \ngsolve. Using $\frac{d T_{s,t}}{ds}(x) = V(x)$, $\frac{d T_{s,t}}{dt}(x) = W(x)$ and $\frac{d F_{s,t}}{ds}(x) = \partial V(x)$, $\frac{d F_{s,t}}{dt}(x) = \partial W(x)$, we get
\begin{align} \label{eq:formula_second_int_f}
\left. \frac{d^2}{dsdt} \mathcal J(\Omega_{s,t}) \right\rvert_{s=t=0} =&   \left( \frac{d^2 G}{d T_{s,t}^2} V + \frac{d^2 G}{d F_{s,t} dT_{s,t}} \partial V  \right) W \nonumber \\
 &+  \left(  \frac{d^2 G}{d T_{s,t} dF_{s,t}} V  + \frac{d^2 G}{d F_{s,t}^2} \partial V \right) \partial W \left. \vphantom{\frac{d^2}{dsdt}}  \right\rvert_{s=t=0}. 
\end{align}

\begin{remark}
We remark that the formula \eqref{eq:formula_second_int_f} can be evaluated explicitly and reads
 \begin{align*}
        D^2\mathcal J(\Omega)(V,W) =& \int_\Omega \nabla^2f V \cdot W + \nabla f \cdot W \, \mbox{div}V+ \nabla f \cdot V \, \mbox{div}W \\
        &+ f \, \mbox{div}V \, \mbox{div}W - f \partial V^\top : \partial W \, \mbox dx.
    \end{align*}
\end{remark}
Formula \eqref{eq:formula_second_int_f} can be implemented in NGSolve as follows:
\begin{lstlisting}[firstnumber=last]
d2JOmega_f = BilinearForm(VEC)
W = VEC.TrialFunction()

d2JOmega_f+=(G_f.Diff(X,W)+G_f.Diff(F,grad(W))).Diff(X,V)(*@\label{lst_diffdiff1}@*)f
d2JOmega_f+=(G_f.Diff(X,W)+G_f.Diff(F,grad(W))).Diff(F,grad(V))(*@\label{lst_diffdiff2}@*)
\end{lstlisting}
Notice that since \texttt{W} is a trial function it is not affected by 
the differentiation with respect to \texttt{X}, see Section \ref{sec_autoDiff}. Therefore, the terms coming from differentiating \texttt{W} with respect to the spatial coordinates \texttt{X} into the direction of \texttt{V} disappear and thus, although code lines \ref{lst_diffdiff1}--\ref{lst_diffdiff2} look like the ``derivative of the derivative'', we actually compute formula \eqref{eq:second_derivative} and not \eqref{eq_d2J_nonsymm}.

In the same fashion, second order derivatives of boundary integrals of the form \eqref{eq:shape_func_f_bnd} can be computed.
\begin{lstlisting}[firstnumber=last]
d2JOmega_f_bnd = BilinearForm(VEC)

d2JOmega_f_bnd+=(G_f_bnd.Diff(X, W)+G_f_bnd.Diff(F,grad(W).Trace())).Diff(X, V)(*@\label{lst_diffdiff3}@*)
d2JOmega_f_bnd+=(G_f_bnd.Diff(X, W)+G_f_bnd.Diff(F,grad(W).Trace())).Diff(F,grad(V).Trace())(*@\label{lst_diffdiff4}@*)
\end{lstlisting}
Again note that the trace operator is necessary when dealing with gradients on the boundary.

\section{Semi-automatic shape differentiation with PDE constraints} \label{sec_PDE_constr_shape}
In this section, we describe the automatic computation of the shape derivative for the following type of equality constrained shape optimisation problems:
\ben \label{eq:PDE_abstract_J}
\min_{(\Omega, u)} J(\Omega, u)  
\een
subject to $(\Omega,u)\in \mathcal A\times Y $ solves
\ben\label{eq:PDE_abstract}
e(\Omega, u) = 0 ,
\een
where $e:\mathcal A\times Y \to Y ^*$ with $e(\Omega,\cdot):Y (\Omega)\to Y (\Omega)^*$ represents an abstract PDE constraint with $ Y =\cup_{\Omega\in \mathcal A}Y (\Omega)$ being the union of Banach spaces $Y (\Omega)$ and $\mathcal A$ a set of admissible shapes. For any given $\Omega \in \mathcal A$ we assume the PDE constraint \eqref{eq:PDE_abstract} to admit a unique solution which we denote by~$u_\Omega$. Moreover, let $\Cj(\Omega) := J(\Omega, u_\Omega)$ denote the reduced cost functional.
By introducing a Lagrangian function, we can henceforth deal with an unconstrained shape function $\mathcal L$ rather than a shape function $\mathcal J$ and a PDE constraint. We introduce the Lagrangian 
\ben
\mathcal L(\Omega, u, p) := J(\Omega,u) + \langle e(\Omega, u),p\rangle.
\een
Now an initial shape $\Omega$ is perturbed by a family of transformations $T_t$, resulting
in a new shape $\Omega_t:= T_t(\Omega)$. Transforming back to the initial shape $\Omega$ 
leads to the Lagrangian:
\ben \label{eq_Gtup}
G(t,u,p) := \mathcal L(T_t(\Omega),\Phi_t(u), \Phi_t(p)), \quad u,p\in Y (\Omega),
\een
where $\Phi_t:Y (\Omega) \to Y (\Omega_t)$ is a bijective mapping. Here the transformation $\Phi_t$ depends on the differential operator involved. For instance 
\begin{itemize}
    \item if $Y (\Omega)=H^1_0(\Omega)$, then $\Phi_t(u) = u\circ T_t^{-1}$,
    \item if $Y (\Omega)=H(\curl,\Omega)$, then $\Phi_t(u) = \partial T_t^{-\top}(u\circ T_t^{-1})$,
    \item if $Y (\Omega)=H(\Div,\Omega)$, then $\Phi_t(u) = \frac{1}{\Det(\partial T_t)} \partial T_t (u\circ T_t^{-1})$.
\end{itemize}
Intuitively the transformations $\Phi_t$ are chosen in such a way that the transformed function $\Phi_t(u)$ still belongs to the same space, but on a different domain. For the above three examples this essentially requires to check how the differential operators $\nabla $, $\curl$ and $\Div$ transform under the change of variables $T_t$, respectively. In fact one can check that 
\begin{align*}
    (\nabla u)\circ T_t & = \partial T_t^{-\top} \nabla(u\circ T_t), \quad u\in H^1_0(\Omega), \\
    (\curl u)\circ T_t & = \frac{1}{\xi(t)} \partial T_t \curl\left( \partial T_t^\top (u\circ T_t) \right), \quad u \in H(\curl,\Omega),\\
    (\Div u)\circ T_t & = \frac{1}{\xi(t)} \Div\left( \xi(t) \partial T_t^{-1} (u\circ T_t) \right), \quad u\in H(\Div,\Omega),
\end{align*}
where $\xi(t):= \Det(\partial T_t)$, see also \cite[Section 3.9]{Monk2003}. The transformation rules are precisely given by 
the respective $\Phi_t$. We also note that for smooth functions this can be checked by direct computation.

Now the shape differentiability of \eqref{eq:PDE_abstract_J}--\eqref{eq:PDE_abstract} is reduced to proving that (see \cite{c_ST_2015a})
\begin{align} \label{eq_dGdt}
    D\mathcal J(\Omega)(V)= \frac{d}{dt} G(t, u^t, 0)|_{t=0}  = \partial_t G(0, u, p),
\end{align}
where $u^t:=u_t\circ T_t$ and $u_t\in Y (\Omega_t)$ solves $e(\Omega_t,u_t)=0$ and $p$ is the solution to the adjoint equation 
\begin{equation}
  p\in Y(\Omega), \quad   \partial_u G(0,u,p)(\varphi) =0 \quad \text{ for all } \varphi \in Y(\Omega). 
\end{equation}
We stress that the choice of $p$ as the solution of the adjoint equation is important in order for the second equality in \eqref{eq_dGdt} to hold. The verification of this equality depends on the specific PDE under consideration and can be accomplished by different methods. We refer the reader 
to \cite{c_ST_2015a} for an overview and remark that \eqref{eq_dGdt} holds for a large class of nonlinear PDE constrained shape optimisation problems; see \cite{a_ST_2015a}. 

The rest of this section is organised as follows: We introduce a model problem, which is the minimisation of a tracking-type cost functional subject to Poisson's equation in Section \ref{sec_laplace}. We illustrate how the first and second order shape derivative for this PDE-constrained model problem can be obtained in \ngsolve{} in Sections \ref{sec_PDE_constr_shape_1st} and \ref{sec_PDE_constr_shape_2nd}. Finally, we also briefly discuss the extension to partial differentiation equations on surfaces. 

\subsection{PDE-constrained model problem} \label{sec_laplace}
We will illustrate the derivation of the first and second order shape derivative for the minimisation of a tracking-type cost functional subject to Poisson's equation on the unknown domain $\Omega$. Let $d=2$ or $3$, $f, u_d \in H^1(\VR^d)$ and $\mathcal A \subset \mathcal P(\VR^d)$ be a set of admissible shapes. Here, $\mathcal P(\VR^d)$ denotes the power set of all subsets of $\VR^d$. We consider the problem
\begin{subequations} \label{eq_modelPDEopti}
    \begin{align}
        &\underset{(\Omega,u) }{\mbox{min }} J(\Omega, u) = \int_\Omega |u-u_d|^2 \; \mbox dx \label{eq_modelPDEopti_J} 
    \end{align}
subject to $(\Omega,u) \in \mathcal A \times H^1_0(\Omega)$ solves 
\begin{align}
         & \langle e(\Omega,u),\psi\rangle := \int_\Omega \nabla u \cdot \nabla \psi \, \mbox dx - \int_\Omega f \psi \, \mbox dx = 0 \label{eq_modelPDEopti_PDE}
    \end{align}
\end{subequations}
for all $\psi \in H_0^1(\Omega)$. The Lagrangian is given by
\ben
\mathcal L(\Omega,\varphi,\psi) := \int_\Omega|\varphi - u_d|^2 \;\mbox dx + \int_{\Omega} \nabla \varphi \cdot\nabla \psi \, \mbox dx - \int_{\Omega} f \psi \,\mbox dx. 
\een
Given an admissible shape $\Omega$, a vector field $V \in C^{0,1}(\VR^d)^d$ and $t>0$ small, let $\Omega_t := (\Id + t V)(\Omega)$ be the perturbed domain. Therefore the parametrised Lagrangian is given by 
\ben
G(t,\varphi,\psi) := \mathcal L(T_t(\Omega) ,\varphi\circ T_t^{-1} ,\psi\circ T_t^{-1}), \quad \varphi,\psi\in H^1_0(\Omega). 
\een
Changing variables yields 
\begin{align} 
    G&(t,\varphi, \psi) = \int_\Omega|\varphi - u_d^t|^2 \, \mbox{det}(F_t) \mbox{d} x \nonumber \\
    &+ \int_{\Omega} (F_t^{-\top}\nabla \varphi) \cdot (F_t^{-\top} \nabla \psi) \,\mbox{det}(F_t) \, \mbox dx - \int_{\Omega} f^t \psi \, \mbox{det}(F_t)\,\mbox dx \label{eq_Lagrangian} \\
    &=: \tilde G(T_t, F_t, \varphi, \psi), \nonumber
\end{align}
where $u_d^t = u_d \circ T_t$ and $f^t = f \circ T_t$. Here we also transformed the gradient according to $(\nabla w)\circ T_t = F_t^{-\top} \nabla (w\circ T_t)$ for $w \in H^1_0(\Omega)$. Recall that, for a given $\Omega \in \mathcal A$, $u_\Omega$ denotes the corresponding unique solution to \eqref{eq_modelPDEopti_PDE} and $\mathcal J(\Omega)$ the reduced cost functional, $\mathcal J(\Omega) := J(\Omega, u_\Omega)$. Let $u^t \in H_0^1(\Omega)$ be the solution of the perturbed state equation brought back to the original domain $\Omega$, that is, $u^t \in H_0^1(\Omega)$ is the unique solution to
\begin{align} \label{eq_pertState}
    \partial_\psi G(t, u^t, 0)(\psi) = 0 \quad \text{ for all } \psi \in H_0^1(\Omega).
\end{align}
Note that, for $u^t$ defined by \eqref{eq_pertState}, it holds $\mathcal J(\Omega_t) = G(t, u^t, \psi)$ for all $\psi \in H_0^1(\Omega)$ and therefore also $D\mathcal J(\Omega)(V) = \frac{d}{dt} G(t, u^t, \psi) $ for all $\psi \in H_0^1(\Omega)$.

It can easily be shown that \eqref{eq_dGdt} holds and thus the shape derivative in the direction of a vector field $V \in C^{0,1}(\VR)^d$ is given by
$$D \mathcal J(\Omega)(V) = \partial_t G(0, u, p), $$
where $p \in H_0^1(\Omega)$ denotes the adjoint state and is defined as the unique solution $p \in H_0^1(\Omega)$ to 
\begin{align}
     \partial_\varphi G(0, u, p)( \hat \varphi) = 0 \quad \text{ for all }  \hat \varphi \in H_0^1(\Omega),
\end{align}
or explicitly
\begin{align} \label{eq_adjoint}
    \int_\Omega \nabla \hat \varphi \cdot \nabla p \, \mbox dx = - 2 \int_\Omega (u-u_d) \hat \varphi \; \mbox dx  \quad \text{ for all }  \hat \varphi \in H_0^1(\Omega).
\end{align}

\subsection{First order shape derivative} \label{sec_PDE_constr_shape_1st}
By the discussion above, the first order shape derivative is given by $\partial_t G(0, u, p)$ with $G$ defined in \eqref{eq_Lagrangian} and $u$ and $p$ the unique solutions to the boundary value problems \eqref{eq_modelPDEopti_PDE} and \eqref{eq_adjoint}, respectively.

Writing $\tilde G(T_t, F_t) := \tilde G(T_t, F_t, u, p) =  G(t,u,p) $ we obtain in analogy to the unconstrained problem
\begin{align*}  
    D \mathcal J(\Omega)(V) =  \left. \frac{d}{dt} \mathcal J(\Omega_t) \right|_{t=0} = \left. \left( \frac{d \tilde G}{dT_t} V + \frac{d \tilde G}{dF_t}\partial V  \right)\right|_{t=0}.
\end{align*}
We can compute explicitly
\begin{align}
\frac{d \tilde G}{dF_t} |_{t=0} \partial V =& \int_\Omega \Div(V)(u-u_d)^2 - (\partial V + \partial V^\top)\nabla u\cdot \nabla p \nonumber \\
&- \Div(V) \nabla u\cdot \nabla p - fp \Div(V)\;\mbox dx, \\
\frac{d \tilde{G}}{dT_t} |_{t=0}V =& \int_\Omega - 2(u-u_d) \nabla u_d\cdot V - \nabla f\cdot V p\;\mbox dx.
\end{align}

Now we are in a position to compute the first order shape derivative for the PDE-constrained shape optimisation problem \eqref{eq_modelPDEopti} in \ngsolve{}. After solving the state equation as shown in Section \ref{sec_intro_ngsolve}, the adjoint equation can be solved as follows.
\begin{lstlisting}[firstnumber = last]
ud = X[0]*(1-X[0])*X[1]*(1-X[1])
def Cost(u):
    return (u-ud)**2 * Det(F) * dx

#solve adjoint equation
gfp = GridFunction(fes)
dCostdu = LinearForm(fes)
dCostdu += Cost(gfu).Diff(gfu, w)
dCostdu.Assemble()
gfp.vec.data = -a.mat.Inverse(fes.FreeDofs(), inverse="sparsecholesky").T * dCostdu.vec

Draw(gfp, mesh, "adjoint")
\end{lstlisting}
We can now define the Lagrangian \eqref{eq_Lagrangian} such that the shape derivative can be obtained by the same procedure as in the unconstrained setting. Note that lines \ref{lst_diff1constr}--\ref{lst_diff2constr} coincide with lines \ref{lst_diff1}--\ref{lst_diff2}.
\begin{lstlisting}[firstnumber = last]
def Equation(u,w):
    return ((Inv(F).trans * grad(u)) * (Inv(F).trans * grad(w)) - f*w)*Det(F)*dx
    
G_pde = Cost(gfu) + Equation(gfu, gfp) (*@\label{lst_def_G_pde} @*)

dJOmega_pde = LinearForm(VEC)
dJOmega_pde += G_pde.Diff(X, V)  (*@\label{lst_diff1constr}@*)
dJOmega_pde += G_pde.Diff(F, grad(V)) (*@\label{lst_diff2constr}@*)
\end{lstlisting}

\subsection{Second order shape derivative} \label{sec_PDE_constr_shape_2nd}

Let us introduce the notation 
\begin{align}
    \langle E_{V,W}(s,t)\varphi,\psi\rangle :=& \int_{\Omega} (F_{s,t}^{-\top}\nabla \varphi) \cdot (F_{s,t}^{-\top} \nabla \psi) \,\mbox{det}(F_{s,t}) \; \mbox dx \nonumber \\
    &- \int_{\Omega} f\circ T_{s,t} \psi \, \mbox{det}(F_{s,t})\; \mbox dx
    \\
    J_{V,W}(s,t;\varphi):= &\int_\Omega|\varphi - u_d\circ T_{s,t} |^2 \, \mbox{det}(F_{s,t}) \mbox{d} x
\end{align}
and
\begin{align}\label{eq_Lagrangian_hess}
 G_{V,W}(s,t,u, p):=&  \langle E_{V,W}(s,t)u,p\rangle  + J_{V,W}(s,t;u),
\end{align}
where $T_{s,t}(x) = x + s V(x) + tW(x)$ and $F_{s,t}:= \partial T_{s,t}$. We observe that
\ben\label{eq:diff_second_J}
\Cj(T_{s,t}(\Omega)) = G_{V,W}(s,t,u^{s,t}, p^{s,t})
\een
with $(u^{s,t},p^{s,t})\in H^1_0(\Omega)\times H^1_0(\Omega)$ being the solution to 
\begin{align}
    \partial_p G_{V,W}(s,t,u^{s,t},0)(\varphi) & = 0 \quad \text{ for all } \varphi\in H^1_0(\Omega), \label{eq_dpG} \\
   \partial_u G_{V,W}(s,t,u^{s,t},p^{s,t})(\psi) & = 0 \quad \text{ for all } \psi\in H^1_0(\Omega) \label{eq_duG}
\end{align}
for $s, t \geq 0$.
In case $t=0$ we write $u^s:= u^{s,t}|_{t=0}$ and $p^s:= p^{s,t}|_{t=0}$ and similarly for $t=s=0$ we write $u:=u^{s,t}|_{s=t=0}$ and $p:=p^{s,t}|_{s=t=0}$.  Therefore, consecutive differentiation of \eqref{eq:diff_second_J} first with respect to $t$ at zero and then with respect to $s$ at zero yields
\begin{align}
    D^2 &\Cj(\Omega)(V)(W)  =  \frac{d^2}{dsdt} G_{V,W}(s,t,u^{s,t},p^{s,t})|_{s=t=0} \nonumber  \\
    =& \frac{d}{ds}\partial_t G_{V,W}(s,0,u^s,p^s)|_{s=0} \nonumber  \\
                       =&  \partial_s\partial_t G_{V,W}(0,0,u,p)
                        + \partial_u \partial_t G_{V,W}(0,0,u,p)(\partial_s u^0) \nonumber \\
                        &+ \partial_p \partial_t  G_{V,W}(0,0,u,p)(\partial_s p^0), \label{eq:diff_second_J2}
\end{align}
where $\partial_s u^0 \in H^1_0(\Omega)$ solves the material derivative equation
\ben\label{eq:material_derivative_1}
\partial_u\partial_p G_{V,W}(0,0,u,0)(\psi)(\partial_s u^0)  = - \partial_s\partial_p G_{V,W}(0,0,u,0)(\psi) 
\een
for all $\psi \in H^1_0(\Omega)$ or, equivalently
\ben\label{eq:material_}
\langle \partial_u E_{V,W}(0,0)(\partial_s u^0),\psi\rangle = - \langle \partial_s E_{V,W}(0,0)u,\psi\rangle 
\een
for all $\psi \in H^1_0(\Omega)$. Note that \eqref{eq:material_} is obtained by differentiating \eqref{eq_dpG} with respect to $s$ and setting $s=t=0$. Similarly the function $\partial_s p^0 \in H^1_0(\Omega)$ solves the material derivative equation obtained by differentiating \eqref{eq_duG} with respect to $s$ for $s=t=0$,
\begin{align}
\partial_p\partial_u G_{V,W}(0,0,u,p)(\psi)(\partial_s p^0) =& - \partial_u^2 G_{V,W}(0,0,u,p)(\psi)(\partial_s u^0) \nonumber \\
&- \partial_s\partial_u G_{V,W}(0,0,u,p)(\psi) \label{eq:material_derivative_2}
\end{align}
for all $\psi \in H^1_0(\Omega)$. The introduction of the adjoint variable $p$ 
is analogous to the computation of the first order shape derivative. However in contrast to the first order derivative the evaluation of $D^2 \Cj(\Omega)(V)(W)$ requires 
the computation of the material derivatives $\partial_s u^0$ and $\partial_s p^0$. 

Formally \eqref{eq:material_derivative_1} and \eqref{eq:material_derivative_2} can be written as an operator equation with $x=(0,0,u,p)$,
\ben\label{eq:material_derivatrive_system}
\begin{pmatrix}
    \partial_u^2 G_{V,W}(x) & \partial_p\partial_u G_{V,W}(x) \\
    \partial_u\partial_p G_{V,W}(x) & 0 
  \end{pmatrix}  \begin{pmatrix}
  \partial_s u^0 \\ \partial_s p^0
   \end{pmatrix} = -\begin{pmatrix}
 \partial_s\partial_u G_{V,W}(x) \\
 \partial_s\partial_p G_{V,W}(x) 
\end{pmatrix}.
\een
So to evaluate the second derivative \eqref{eq:diff_second_J2} in some direction $(V,W)$ we have to solve the system \eqref{eq:material_derivatrive_system}.

This is realised in \ngsolve{} by setting up a combined finite element space which we denote by \texttt{X2}. We define trial and test functions as well as grid functions representing the deformation vector fields $V$ and $W$, which we initialise with some functions.
\begin{lstlisting}[firstnumber=last]
X2 = FESpace([fes, fes])
dsu, dsp = X2.TrialFunction()
uTest, pTest = X2.TestFunction()
gfV = GridFunction(VEC)
gfW = GridFunction(VEC)
gfV.Set((X[0]*X[0]*X[1]*exp(X[1]),X[1]*X[1]*X[0]*exp(X[0]))) 
gfW.Set((X[1]*X[1]*X[0]*exp(X[0]),X[0]*X[0]*X[1]*exp(X[1]))) 
\end{lstlisting}
We define a 2$\times$2 block bilinear form as well as a 2$\times$1 block linear form which will represent the left and right hand sides of \eqref{eq:material_derivatrive_system}, respectively. The operator equation in \eqref{eq:material_derivatrive_system} can be conveniently defined by differentiating the Lagrangian with respect to the corresponding variables.
\begin{lstlisting}[firstnumber = last]
shapeHessLag2 = BilinearForm(X2)
shapeGradLag2 = LinearForm(X2)

shapeHessLag2 += (G_pde.Diff(gfu, uTest)).Diff(gfu, dsu)   #block (1,1)
shapeHessLag2 += (G_pde.Diff(gfu, uTest)).Diff(gfp, dsp)   #block (1,2)
shapeHessLag2 += (G_pde.Diff(gfp, pTest)).Diff(gfu, dsu)   #block (2,1)
    
#line 1
shapeGradLag2 += (G_pde.Diff(gfu, uTest)).Diff(F, grad(gfV)) 
shapeGradLag2 += (G_pde.Diff(gfu, uTest)).Diff(X, gfV) 

#line 2
shapeGradLag2 += (G_pde.Diff(gfp,pTest)).Diff(F, grad(gfV))
shapeGradLag2 += (G_pde.Diff(gfp,pTest)).Diff(X, gfV) 
\end{lstlisting}
We can solve this combined system for $\partial_s u^0$ and $\partial_s p^0$ and access and visualise the two components in the following way:
\begin{lstlisting}[firstnumber = last]
gfCombined2 = GridFunction(X2)
shapeHessLag2.Assemble()
shapeGradLag2.Assemble()
gfCombined2.vec.data = shapeHessLag2.mat.Inverse(X2.FreeDofs(), inverse = "umfpack") * shapeGradLag2.vec

gfdsu = GridFunction(fes)
gfdsp = GridFunction(fes)
gfdsu.vec.data = gfCombined2.components[0].vec
gfdsp.vec.data = gfCombined2.components[1].vec

Draw(gfdsu, mesh, "dsu")
Draw(gfdsp, mesh, "dsp")
\end{lstlisting}
In order to obtain the second order shape derivative in the direction given by $(V,W)$, it remains to evaluate the term \eqref{eq:diff_second_J2}. We define the three terms of \eqref{eq:diff_second_J2} as bilinear forms, assemble them and perform vector-matrix-vector multiplications:
\begin{lstlisting}[firstnumber = last]
w1 = fes.TrialFunction()
q1 = fes.TrialFunction()

shapeHess11 = BilinearForm(VEC)
shapeHess11 += (G_pde.Diff(F,grad(W))+G_pde.Diff(X, W)).Diff(F,grad(V)) (*@\label{lst_11block1}@*) 
shapeHess11 += (G_pde.Diff(F,grad(W))+G_pde.Diff(X, W)).Diff(X, V) (*@\label{lst_11block2}@*)
shapeHess11.Assemble()

shapeHess12 = BilinearForm(trialspace = fes, testspace = VEC)
shapeHess12 += (G_pde.Diff(F, grad(V)) + G_pde.Diff(X, V)).Diff(gfu,w1)    
shapeHess12.Assemble()

shapeHess13 = BilinearForm(trialspace = fes, testspace = VEC)
shapeHess13 += (G_pde.Diff(F, grad(V)) + G_pde.Diff(X, V)).Diff(gfp,q1)  
shapeHess13.Assemble()

av = gfV.vec.CreateVector()
av.data = shapeHess11.mat * gfV.vec

adsu = gfV.vec.CreateVector()
adsu.data = shapeHess12.mat * gfdsu.vec

adsp = gfV.vec.CreateVector()
adsp.data = shapeHess13.mat * gfdsp.vec

d2J =  InnerProduct(gfW.vec, av) + InnerProduct(gfW.vec, adsu) + InnerProduct(gfW.vec, adsp) (*@\label{lst_d2J}@*)
\end{lstlisting}

\subsection{PDEs on surfaces} \label{sec_surfaces}
The automated shape differentiation is not restricted to partial differential equations on domains $\Omega$, but is readily extended to surface PDEs. We consider a two dimensional closed surface $M\subset \VR^3$ and denote by $\normal$ the normal field along $M$. Let $u_d\in H^1(\VR^d)$ be given and define
\ben \label{eq_surface_J}
J(M,u) = \int_M |u-u_d|^2\;\mbox ds,
\een
where $u\in H^1(M)$ solves the surface equation
\ben \label{eq_surface_pde}
\int_M \nabla^M u\cdot \nabla^M \psi + u \psi \;\mbox ds = \int_M f \psi \;\mbox ds \quad \text{ for all } \psi \in H^1(M),
\een
where $\nabla^M \psi $ denotes the tangential gradient of $\psi $; see \cite[p.493, Def.5.1]{b_DEZO_2011a}.  We assume that the function $f\in H^1(\VR^3)$ is given. The Lagrangian is given by
\begin{align*}
\mathcal L(M,\varphi,\psi) :=& \int_M |\varphi - u_d|^2 \;\mbox ds + \int_M \nabla^M \varphi \cdot\nabla^M \psi  + \varphi\psi\; \mbox  ds \\
&- \int_M f \psi \; \mbox ds. 
\end{align*}
As in the previous section we fix an admissible shape $M$ and let $M_t := (\Id + t V)(M)$ be a small 
perturbation of $M$ by means of a vector field $V \in C^1(\VR^d)^d$ for $t>0$ small. The parametrised Lagrangian is given by 
\ben
G(t,\varphi,\psi) := \mathcal L(T_t(M) ,\varphi\circ T_t^{-1}, \psi\circ T_t^{-1}), \quad \varphi,\psi\in H^1(M). 
\een
Define the density $\tangDet(F_t) := \Det(F_t) |F_t^{-\top}\normal|$. Changing variables and using 
\begin{align}
(\nabla^{M_t}\varphi)\circ T_t =& B(F_t)\nabla^M(\varphi\circ T_t), \nonumber \\
B(F_t) =& \left(I-\frac{F_t^{-\top}\normal}{|F_t^{-\top}\normal|}\otimes \frac{F_t^{-\top}\normal}{|F_t^{-\top}\normal|}\right) F_t^{-\top},  \label{eq_defB}
\end{align}
yields 
\begin{align} \label{eq_Lagrangian_manifold}
    \begin{aligned}
    G(t,\varphi, \psi) =& \int_M|\varphi - u_d^t|^2 \, \tangDet(F_t) \;\mbox ds \\
    &+ \int_{M}  \left( (B(F_t)\nabla \varphi) \cdot (B(F_t) \nabla \psi) + \varphi\psi \right) \tangDet(F_t)\; \mbox ds \\
    &- \int_M   f^t \psi  \,\tangDet(F_t)\;\mbox ds,
    \end{aligned}
\end{align}
where $u_d^t = u_d \circ T_t$ and $f^t = f \circ T_t$. 

Writing $\tilde G(T_t, F_t) := G(t,u,p) $ we obtain in analogy to the domain case 
\begin{align}  \label{eq_shapeDerivativeSurface}
    D \mathcal J(\Omega)(V) =  \left. \left( \frac{d \tilde G}{dT_t} V + \frac{d \tilde G}{dF_t} \partial V\right)\right|_{t=0}.
\end{align}
We can compute explicitly
\begin{align}
    \frac{d \tilde G}{dF_t}|_{t=0}V = \int_M& \Div^M(V)(u-u_d)^2 \nonumber \\
    &- (\partial^M  V + \partial^M V^\top)\nabla^M u\cdot \nabla^M p \nonumber \\
    &+ \Div^M(V) (\nabla^M u\cdot \nabla^M p + up) \nonumber\\
    &- fp \Div^M(V)\; \mbox ds, \\
    \frac{d \tilde G}{dT_t} |_{t=0} \partial V = \int_M& - 2(u-u_d) \nabla u_d\cdot V - \nabla f\cdot V p\; \mbox ds,
\end{align}
where $\partial^M V$ denotes the tangential Jacobian of $V$ defined by $(\partial^M V)_{ij} := (\nabla^M V_i)_j$ for $i,j=1,\dots,d$, and $\Div^M(V):= \partial^M V:I$ the tangential divergence, which is defined as the trace of the tangential Jacobian; see \cite[p.495]{b_DEZO_2011a}.

The implementation is analogous to the previous sections. We will only illustrate first order derivatives here. We first define the geometry of the unit sphere, create a surface mesh and define a finite element space on the surface mesh:
\begin{lstlisting}[firstnumber=last]
from netgen.csg import *
from netgen.meshing import *
from ngsolve.internal import visoptions
from ngsolve import *

geo_surf = CSGeometry()
sphere = Sphere(Pnt(0,0,0),1).bc("outer")
geo_surf.Add(sphere)
mesh_surf = Mesh(geo_surf.GenerateMesh(perfstepsend=MeshingStep.MESHSURFACE,optsteps2d=3,maxh=0.2))
mesh_surf.Curve(3)
fes_surf = H1(mesh_surf,order = 3)
\end{lstlisting}
Next we define the transformed cost function and partial differential equation needed for setting up the Lagrangian \eqref{eq_Lagrangian_manifold}. Here, we again make use of a symbolic object \texttt{F} to which we assign the identity matrix. We define the tangential determinant $\tangDet$ and the matrix $B$ defined in \eqref{eq_defB} as functions of the deformation gradient $F_t$.
\begin{lstlisting}[firstnumber=last]
X = CoefficientFunction((x,y,z))
func = CoefficientFunction(X[0]*X[1]*X[2])
F = Id(3)
tangDet = Det(F) * Norm( Inv(F).trans*specialcf.normal(3) )
Bmat = (Id(3) - 1/Norm(Inv(F).trans*specialcf.normal(3))**2 * OuterProduct(Inv(F).trans*specialcf.normal(3), Inv(F).trans*specialcf.normal(3)) ) * Inv(F).trans

def Equation_surf(u,w): (*@\label{lst_eqnSurf1}@*)
    return ( (Bmat*grad(u).Trace()) * (Bmat*grad(w).Trace()) + u*w - func * w) * tangDet * ds   (*@\label{lst_eqnSurf2}@*)

def Cost_surf(u):
    return u**2 * tangDet * ds
\end{lstlisting}
Now we can define the bilinear form and solve the state equation. Here, the right hand side of the equation is included in the bilinear form and the boundary value problem -- although linear -- is solved by Newton's method (which terminates after only one iteration) for convenience.
\begin{lstlisting}[firstnumber = last]
#set up and solve state equation
u_surf, w_surf = fes_surf.TnT()
a = BilinearForm(fes_surf)
a += Equation_surf(u_surf, w_surf) (*@\label{lst_a_Eqn}@*)
gfu_surf = GridFunction(fes_surf)
solvers.Newton(a, gfu_surf, printing = False)
Draw(gfu_surf, mesh_surf, "gfu_surf")
\end{lstlisting}
Using Newton's method for solving the linear boundary value problem allows us to define both the left and right hand side of the PDE using only one \texttt{BilinearForm} \texttt{a} (which, strictly speaking, is not bilinear any more). This way, we can reuse \texttt{Equation\_surf} as defined in lines \ref{lst_eqnSurf1}--\ref{lst_eqnSurf2} to define the boundary value problem in line \ref{lst_a_Eqn}. 

The adjoint equation is solved as usual:
\begin{lstlisting}[firstnumber=last]
#solve adjoint equation
lfcost_surf = LinearForm(fes_surf)
lfcost_surf += Cost_surf(gfu_surf).Diff(gfu_surf, w_surf)
lfcost_surf.Assemble()
inva = a.mat.Inverse(fes_surf.FreeDofs(), inverse="sparsecholesky")
gfp_surf = GridFunction(fes_surf)
gfp_surf.vec.data = -inva.T * lfcost_surf.vec
Draw(gfp_surf, mesh_surf, "gfp_surf")
\end{lstlisting}
The shape derivative is obtained as in the case of PDEs posed on volumes by the evaluation of \eqref{eq_shapeDerivativeSurface}:
\begin{lstlisting}[firstnumber=last]
G_surf = Cost_surf(gfu_surf) + Equation_surf(gfu_surf, gfp_surf)

VEC3d = VectorH1(mesh_surf, order=1)
V3d = VEC3d.TestFunction()
dJOmega_surf = LinearForm(VEC3d)
dJOmega_surf += G_surf.Diff(X, V3d) + G_surf.Diff(F, Grad(V3d).Trace())
\end{lstlisting}

\section{Fully automated shape differentiation} \label{sec_fullyAutomated}
In the previous sections we used the automatic differentiation capabilities of \ngsolve{} to alleviate the shape differentiation procedure. However, so far we still had to include some knowledge about the problems at hand. So far, it was necessary to define the objective function or Lagrangian $G$ in the correct way, accounting for the correct transformation rules between perturbed and unperturbed domain. In this section, we will show that also this step can be automated since all necessary information is already included in the functional setting. The fully automated shape differentiation is incorporated by the command 
\begin{center}
    \texttt{DiffShape(...)}.
\end{center}
In particular, in the fully automated setting it is enough to set up the cost function or Lagrangian for the unperturbed setting. For a shape function of the type \eqref{eq:shape_func_f} we can define the shape derivative of the cost function in the following way:
\begin{lstlisting}[firstnumber=last]
G_f_0 = f * dx (*@\label{lst_DiffShape1}@*)
dJOmega_f_0 = LinearForm(VEC)
dJOmega_f_0 += G_f_0.DiffShape(V) (*@\label{lst_diffShape}@*)
\end{lstlisting}
Note that there is no term of the form \texttt{Det(F)} showing up in line \ref{lst_DiffShape1}. Here, the transformation of the domain is taken care of automatically. It can be checked that this really gives the same result as \texttt{dJOmega\_f} defined in lines \ref{lst_diff1}--\ref{lst_diff2}.
\begin{lstlisting}[firstnumber=last]
dJOmega_f.Assemble()
dJOmega_f_0.Assemble()
differenceVec = dJOmega_f.vec.CreateVector()
differenceVec.data = dJOmega_f.vec - dJOmega_f_0.vec
print("|dJOmega_f - dJOmega_f_0| = ", Norm(differenceVec) )
\end{lstlisting}
The above code gives the output 
\begin{verbnobox}[\footnotesize]
|dJOmega_f - dJOmega_f_0| = 1.571008573810619e-17
\end{verbnobox}    
which confirms our claim. The same holds true for second order shape derivatives. The lines \ref{lst_diffdiff1}--\ref{lst_diffdiff2} can be replaced by a repeated call of \texttt{DiffShape(...)}:
\begin{lstlisting}[firstnumber=last]
d2JOmega_f_0 = BilinearForm(VEC)
d2JOmega_f_0+=G_f_0.DiffShape(V).DiffShape(W)
\end{lstlisting}

Again, it can be verified that \texttt{d2JOmega\_f\_0} coincides with the previously defined quantity \texttt{d2JOmega\_f}. Note that slightly different results may occur due to different integration rules used. This can be cured by enforcing an integration rule of higher order for \texttt{G\_f}, i.e. by replacing the symbol \texttt{dx} in the definition of \texttt{G\_f} with \texttt{dx(bonus\_intorder=2)}.

In the more general setting of PDE-constrained shape optimisation, the procedure is very similar. Here the idea exploited in the implementation of the command \texttt{DiffShape(...)} is to just differentiate the general expression \eqref{eq_Gtup} with respect to the parameter $t$. The transformations $\Phi_t$ appearing in \eqref{eq_Gtup}, which depend on the functional setting of the PDE, are identified automatically from the finite element space from which the corresponding functions originate. The shape derivative of lines \ref{lst_diff1constr}--\ref{lst_diff2constr} can be obtained by the following code.
\begin{lstlisting}[firstnumber=last]
def Cost_0(u):
    return (u-ud)**2 * dx

def Equation_0(u,w):
    return (grad(u) * grad(w) - f1*w) *dx

G_pde_0 = Cost_0(gfu) + Equation_0(gfu, gfp)

dJOmega_pde_0 = LinearForm(VEC)
dJOmega_pde_0 += G_pde_0.DiffShape(V)
\end{lstlisting}
Here, \texttt{gfu} and \texttt{gfp} represent the solutions to the state and adjoint equation, respectively, and must have been computed previously. The bilinear form \texttt{shapeHess11} used in Section \ref{sec_PDE_constr_shape_2nd} (see lines \ref{lst_11block1}--\ref{lst_11block2}) can be obtained similarly: 
\begin{lstlisting}[firstnumber=last]
shapeHess11_0 = BilinearForm(VEC)
shapeHess11_0 += G_pde_0.DiffShape(W).DiffShape(V)
\end{lstlisting}
The same holds true for boundary integrals 
\begin{lstlisting}[firstnumber=last]
G_f_bnd_0 = f * ds
dJOmega_f_bnd_0 = LinearForm(VEC)
dJOmega_f_bnd_0 += G_f_bnd_0.DiffShape(V)
\end{lstlisting}
and surface PDEs
\begin{lstlisting}[firstnumber=last]
def Cost_surf_0(u):
    return u**2 * ds
def Equation_surf_0(u,w):
    return  (grad(u).Trace()*grad(w).Trace() + u*w - func * w) * ds
G_surf_0 = Cost_surf_0(gfu_surf) + Equation_surf_0(gfu_surf, gfp_surf)
dJOmega_surf_0 = LinearForm(VEC3d)
dJOmega_surf_0 += G_surf_0.DiffShape(V3d) (*@\label{lst_diffShapeSurf}@*)
\end{lstlisting}
as well as their respective second order derivatives.

\begin{remark}
We remark that the fully automated differentiation using \texttt{DiffShape(...)} should be seen to complement the semi-automated shape differentiation techniques introduced in Sections \ref{sec_intf} and \ref{sec_PDE_constr_shape} rather than to replace them. Using the semi-automated differentiation, the user has the possibility to, on the one hand, keep control over the involved terms, and on the other hand also to adjust the shape differentiation to their custom problems which may be non-standard. As an example where the semi-automated differentiation may be beneficial compared to the fully automated differentiation we mention the case of time-dependent PDE constraints considered in a space-time setting when a shape deformation is only desired in the spatial coordinates, see Section \ref{sec_spacetime}. Of course, when one is interested in the shape derivative for a more standard problem, the fully automated way appears to be more convenient and less error prone.
\end{remark}

\begin{remark}
    We have seen that the command \texttt{DiffShape(\dots)} allows to compute the shape derivative of unconstrained shape optimisation problems in a fully automated way without specifying any transformation rules, see line \ref{lst_diffShape}.
    For the practically more relevant case of PDE-constrained shape optimisation problems, the state and adjoint equations have to be solved beforehand also in the fully automated context using \texttt{DiffShape(\dots)}. We remark that this can be easily achieved by defining a custom function \texttt{solvePDE()} as it is done for the case of a linear PDE in lines \ref{lst_solvePDE1}--\ref{lst_solvePDE2}.    
    Since the purpose of this paper is to illustrate a convenient way of computing shape derivatives and performing shape optimisation rather than to provide a tool for black-box optimisation, this step is left to the user and is not automated, leaving more freedom in the choice of, e.g., solvers for the arising linear systems.
\end{remark}

\section{Optimisation algorithms}
In this section we discuss how to use optimisation algorithms in conjunction with the 
automated shape differentiation explained in the previous sections. The starting point of our 
discussion is a fixed initial shape $\Omega$. Then we consider the mapping
\ben
V\mapsto g(V) := \Cj((\Id+V)(\Omega))
\een
defined on a suitable space of vector fields $\Theta\subset C^{0,1}(\Dsf)^d$. Since the mapping $g$ is defined on an open subset $\Theta$ of the Banach space $C^{0,1}(\Dsf)^d$ we can employ standard algorithms to minimise $g$ over $\Theta$. The only constraint we must impose is that  $\Id+V$ remains invertible, which can be difficult in practice. 
In view of $g(V+t W) = \mathcal J((\Id+V+tW)(\Omega)) = \mathcal J((\Id + t W\circ (\Id+V)^{-1} )((\Id+V)(\Omega)))$ for $V,W\in \Theta$ and $t$ small, we find by differentiating with respect to $t$ at $t=0$, that 
\ben
\partial g(V)(W) = D \mathcal J((\Id+V)(\Omega))(W\circ (\Id + V)^{-1})
\een
for $V,W\in \Theta$ and $\Id+V$ invertible.

\subsection{Gradient computation}

The gradient of $\partial g(V)$ in a Hilbert space $H\subset C^{0,1}(\Dsf)^d$ is defined by 
\ben \label{eq_shapeGrad}
\partial g(V)(W) = (\nabla^H g(V),W)_H \quad \text{ for all }W\in H. 
\een
Typical choices for $H$ are
\begin{align}
    H=H^1_0(\Dsf)^d, \; (W,V)_H:=& \int_\Dsf \partial W:\partial V+V\cdot W\;\mbox dx, \label{eq_shapeGrad_H1} \\
    H=H^1_0(\Dsf)^d, \; (W,V)_H:=& \int_\Dsf \eps(W):\eps(V)+V\cdot W\;\mbox dx,  \label{eq_shapeGrad_ela} \\
    H=H^1_0(\Dsf)^d, \; (W,V)_H:=& \int_\Dsf \eps(W):\eps(V)+ V\cdot W \nonumber \\
    & \quad + \gamma_{CR} \Cb V\cdot \Cb W\;\mbox dx, \label{eq_shapeGrad_ela_CR}
\end{align}
where $\eps(V):= \frac12(\partial V+\partial V^\top)$, $\gamma_{CR}>0$ and
\ben
\Cb := \begin{pmatrix} -\partial_x &  \partial_y \\ \partial_y & \partial_x  \end{pmatrix}.
\een
The last choice, which is restricted to the spatial dimension $d=2$, corresponds to a penalised Cauchy-Riemann gradient and results in a gradient which is approximately conformal and hence preserves good mesh quality. We refer to \cite{a_IGSTWE_2018a} for a detailed description. We also refer to \cite{a_GO_2006a,a_BU_2002a} and \cite{a_ALDAJO_2020a} for the use of different inner products.

\subsection{Basic algorithm}
Let $\Omega$ be an initial shape and let $H\subset C^{0,1}(\Dsf)^d$ be a Hilbert space. 
Then a basic shape optimisation algorithm reads as follows.
\begin{algorithm}[H] 
    \begin{algorithmic}[1]
        \State{{\bf Input:} domain $\Omega_0$, $n=0$, $N_{max}>0$, $\epsilon >0$, $\gamma \geq 0$  }
        \State {\bf Output:} optimal shape $\Omega^*$ 
        \While{ $n\le N_{max}$ and $|\nabla \Cj(\Omega_n)|>\epsilon$}
        \If{$\Cj((\Id - \alpha \nabla \Cj(\Omega_n))(\Omega_n)) < \Cj(\Omega_n) - \gamma \alpha |\nabla \Cj(\Omega_n)|^2$  }
			\State $\Omega_{n+1} \gets (\Id - \alpha \nabla \Cj(\Omega_n))(\Omega_n)$
            \State $n\gets n+1$ 
            \State increase $\alpha$
		\Else
		\State reduce $\alpha$
        \EndIf
	\EndWhile
\caption{gradient algorithm}
\label{alg:gradient}
\end{algorithmic}
\end{algorithm}

We present and explain the numerical realisation of Algorithm \ref{alg:gradient} in \ngsolve{} for the case of a PDE-constrained shape optimisation problem in two space dimensions. The simpler case of an unconstrained shape optimisation problem or the case of three space dimensions can be realised by small modifications of the presented code.

First of all, we mention that we realise shape modifications in \ngsolve{} by means of deformation vector fields without actually modifying the coordinates of the underlying finite element grid. Recall the vector-valued finite element space \texttt{VEC} over a given mesh as introduced in code line \ref{lst_def_VEC}. We define a vector-valued \texttt{GridFunction} with the name \texttt{gfset} which will represent the current shape. We initialise it with some vector-valued coefficient function $V(x_1, x_2) = (x_1^2 \, x_2, x_2^2 \, x_1)^\top$ and obtain the deformed shape $(\Id + V)(\Omega)$ by the command \texttt{mesh.SetDeformation(gfset)}:
\begin{lstlisting}[firstnumber=last]
gfset = GridFunction(VEC)
Draw(gfset, mesh, "gfset")
SetVisualization (deformation=True)
gfset.Set((X[0]*X[0]*X[1],X[1]*X[1]*X[0]))
mesh.SetDeformation(gfset)
Redraw()
\end{lstlisting}
Any operation involving the mesh such as integration or assembling of matrices is now carried out for the deformed configuration. To be more precise, a change of variables is performed internally by accounting for the corresponding Jacobi determinant and transforming the derivatives accordingly with the Jacobian of the deformation. Therefore, all resulting coefficient vectors (which are stored in \texttt{GridFunction}s) correspond to the shape functions in reference configuration. The deformation can be unset by the command\\ \texttt{mesh.UnsetDeformation()}. Integrating the constant function over the mesh in the perturbed and unperturbed setting,
\begin{lstlisting}[firstnumber=last]
print(Integrate(1, mesh))
mesh.UnsetDeformation()
print(Integrate(1, mesh))
\end{lstlisting}
gives the output
\begin{verbnobox}[\footnotesize]
1.7924529046862627
0.7854072970684544
\end{verbnobox}
respectively.

In the course of the optimisation algorithm the state equation as well as the adjoint equation have to be solved for every new shape. We define the following function, which computes the state and adjoint state for a linear PDE constraint:
\begin{lstlisting}[firstnumber=last]
def solvePDE(): (*@\label{lst_solvePDE1}@*)
    a.Assemble()
    L.Assemble()
    dCostdu.Assemble()

    inva = a.mat.Inverse(fes.FreeDofs(), inverse="sparsecholesky")    
    gfu.vec.data = inva * L.vec
    gfp.vec.data = -inva.T * dCostdu.vec (*@\label{lst_solvePDE2}@*)
\end{lstlisting}
The shape derivative \texttt{dJOmega} for some problem at hand can be defined as illustrated in Sections \ref{sec_laplace} and Section \ref{sec_fullyAutomated}.
Finally, we need to define the shape gradient, which is the solution to a boundary value problem of the form \eqref{eq_shapeGrad}. We choose the bilinear form defined in \eqref{eq_shapeGrad_ela_CR} with $\gamma_{CR} = 10$:
\begin{lstlisting}[firstnumber=last]
def eps(u):
    return 1/2 * (grad(u)+grad(u).trans)

aX = BilinearForm(VEC)
W, V = VEC.TnT()    # define trial function W and test function V

aX += InnerProduct(eps(W), eps(V))*dx + InnerProduct(W, V)*dx
aX += 10 * (grad(W)[1,1] - grad(W)[0,0])*(grad(V)[1,1] - grad(V)[0,0])*dx
aX += 10 * (grad(W)[1,0] + grad(W)[0,1])*(grad(V)[1,0] + grad(V)[0,1])*dx
\end{lstlisting}
Now we can run Algorithm \ref{alg:gradient} for problem \eqref{eq_modelPDEopti}:
\begin{lstlisting}[firstnumber=last]
alpha = 1 (*@\label{lst_algo1_l1}@*)
alpha_incr_factor = 1.2
gamma = 1e-4
Nmax = 100
epsilon = 1e-7

isConverged = False
gfset.Set((0,0))
gfX = GridFunction(VEC)
gfsetTemp = GridFunction(VEC)

solvePDE()
Jnew = Integrate(Cost(gfu), mesh)
Jold = Jnew

for k in range(Nmax):
    mesh.SetDeformation(gfset)
    aX.Assemble()
    dJOmega_pde.Assemble()
    invaX = aX.mat.Inverse(VEC.FreeDofs(), inverse="sparsecholesky")
    gfX.vec.data = invaX * dJOmega_pde.vec    
    currentNormGFX = Norm(gfX.vec)
    
    while True:
        if currentNormGFX < epsilon:
            isConverged = True
            break
        
        gfsetTemp.vec.data = gfset.vec - alpha * gfX.vec
        mesh.SetDeformation(gfsetTemp)
        solvePDE()
        Jnew = Integrate (Cost(gfu), mesh)
        mesh.UnsetDeformation()
        if Jnew < Jold - gamma * alpha * currentNormGFX**2 :
            Jold = Jnew
            gfset.vec.data = gfsetTemp.vec
            alpha *= alpha_incr_factor
            break
        else:
            alpha = alpha / 2
    Redraw(blocking=True)   (*@\label{lst_algo1_lN}@*)
\end{lstlisting}

\paragraph{Mesh movement and mesh optimisation}
As an alternative to realizing the deformations via \texttt{mesh.SetDeformation(...)}, where the underlying mesh is not modified, one could also just move every mesh node in the direction of the given descent vector field by changing its coordinates. This can be realised by invoking the following method:
\begin{lstlisting}[firstnumber = last]
def moveNGmesh2D(displ, mesh):
    for p in mesh.ngmesh.Points():
        v = displ(mesh(p[0],p[1]))(*@\label{lst_evalgf}@*)
        p[0] += v[0]
        p[1] += v[1]
    mesh.ngmesh.Update()    (*@\label{lst_meshupdate}@*)
\end{lstlisting}
Here, the displacement vector field \texttt{displ}, which is of type \texttt{GridFunction}, is evaluated for each mesh node and, subsequently, the mesh nodes are updated. At the end of the procedure, the mesh structure needs to be updated, see line \ref{lst_meshupdate}. 
Note that \texttt{GridFunction}s can only be evaluated at points inside the mesh (but not necessarily vertices of the mesh).
Therefore, in order to evaluate \texttt{displ} at the point given by the coordinates \texttt{p[0], p[1]}, we need to pass \texttt{mesh(p[0],p[1])} in line \ref{lst_evalgf}.

\begin{figure}
    \begin{tabular}{cc}
        \includegraphics[width=.22\textwidth]{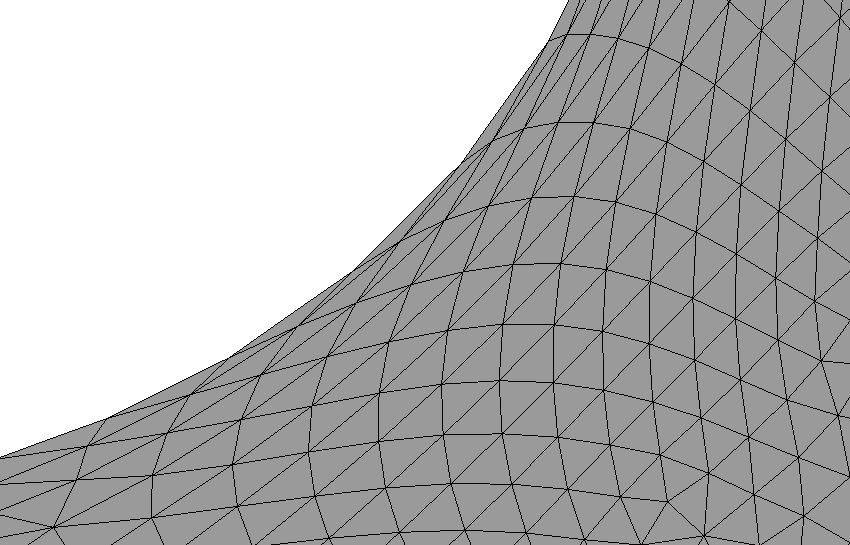} &
        \includegraphics[width=.22\textwidth]{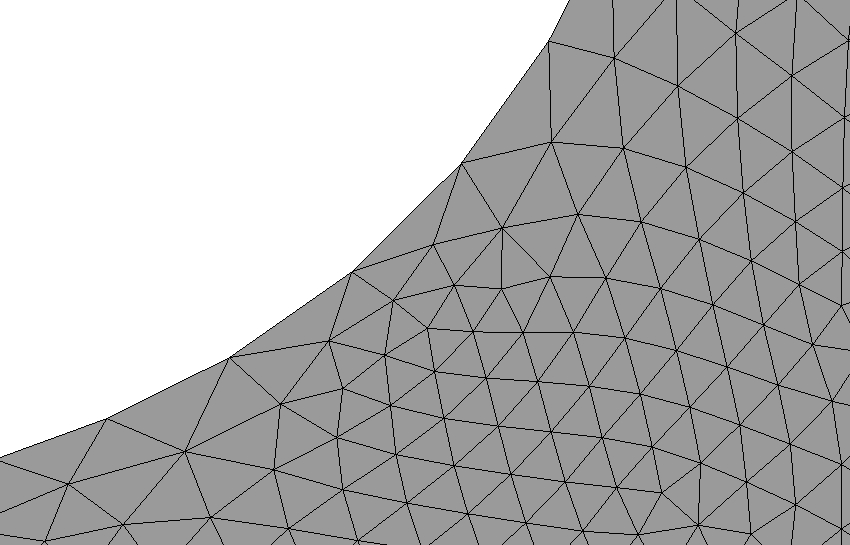} 
    \end{tabular}
    \caption{Before and after mesh optimisation by \texttt{mesh.ngmesh.OptimizeMesh2d()}. }
    \label{fig_clover_results_meshOpti}
 \end{figure}

One advantage of this strategy is that a ill-shaped mesh can easily be repaired by a call of the method \\ \texttt{mesh.ngmesh.OptimizeMesh2d()} followed by \\ \texttt{mesh.ngmesh.Update()}. Figure \ref{fig_clover_results_meshOpti} shows a ill-shaped mesh and the result of a call of \texttt{mesh.ngmesh.OptimizeMesh2d()}.

\subsection{Newton's method for unconstrained problems} \label{sec_newtonUnConstr}

The particular choice $H=H^1_0(\Dsf)^d$ and 
\ben \label{eq_bilinearFormShapeHessian}
(V,W)_H := D^2\Cj(\Omega)(V)(W),
\een
for a given shape function $\Cj$ leads to Newton's method. We refer to \cite{a_NORO_2002a,a_ALCAVI_2016a,a_PAST_2019a,a_EPHASC_2007a} where shape Newton methods were used previously and to \cite[Chapter 2]{b_HIULUL_2009a} and \cite[Chapter 5]{b_ITKU_2008a} for Newton's method in an optimal control setting. This bilinear form is only positive semi-definite on 
$H^1_0(\Dsf)^d$ since $ D^2\Cj(\Omega)(V)(W)=0$ for $V,W$ with $V=W=0$ on $\partial \Omega$. Moreover, from the structure theorem for second shape derivatives proved in \cite{a_NOPI_2002a} we know that at a stationary point $\Omega$, that is, $D\Cj(\Omega)(V)=0$ for all $V\in C^{0,1}(\Dsf)^d$, we have
\ben 
D^2\Cj(\Omega)(V)(W) = \ell_\Omega(V\cdot \normal,W\cdot\normal),
\een
where $\ell_\Omega:C^0(\partial \Omega)\times C^0(\partial \Omega)\to \VR$ is a bilinear function. 
Hence we also have $D^2\Cj(\Omega)(V)(W)=0$ for all $V,W$ such that $V\cdot\normal=W\cdot\normal=0$. As a result the gradient
\ben \label{eq_newtonEq}
(\nabla \Cj(\Omega),V)_H = D\Cj(\Omega)(V) \quad\text{ for all }V\in H^1_0(\Dsf)^d
\een
according to \eqref{eq_bilinearFormShapeHessian} is not uniquely determined. To get around this difficulty, the shape Hessian is often regularised by an $H^1$ term, i.e. \eqref{eq_bilinearFormShapeHessian} is replaced by 
\ben \label{eq_regH1}
    D^2\Cj(\Omega)(V)(W) + \delta \int_\Omega \partial V : \partial W + V \cdot W \; \mbox dx,
\een
see, e.g. \cite{a_SC_2018a}, which, however, impairs the convergence speed of Newton's method. 
\paragraph{Alternative regularisation strategy.}
Here, we propose the following strategy: We regularise the shape Hessian only on the boundary $\partial \Omega$ and only in tangential direction, i.e., we choose
\ben \label{eq_regBnd}
    (V,W)_H := D^2\Cj(\Omega)(V)(W) + \delta \int_{\partial \Omega} (V\cdot \tau)(W \cdot \tau)
\een
with a regularisation parameter $\delta$. To exclude the part of the kernel corresponding to interior deformations, we solve the (regularised) Newton equation \eqref{eq_newtonEq} only on the boundary $\partial \Omega$. This is realised by setting Dirichlet boundary conditions for all degrees of freedom except those on the boundary.
\begin{lstlisting}[firstnumber=last]
VEC2 = VectorH1(mesh, order=1, dirichlet = "circle") #auxiliary space for boundary conditions
aX = BilinearForm(VEC)
aX += G_f_0.DiffShape(W).DiffShape(V)
aX += 100 * InnerProduct(W, specialcf.tangential(2)) * InnerProduct(V, specialcf.tangential(2))*ds
aX.Assemble()
invAX = aX.mat.Inverse(~VEC2.FreeDofs(), inverse="umfpack")

gfX_bnd = GridFunction(VEC)
gfX_bnd.vec.data = invAX * dJOmega_f_0.vec
\end{lstlisting}
As a result, we get a shape gradient $\tilde{\nabla \Cj(\Omega)}$ which is nonzero only on the boundary. We extend this vector field to the interior by solving an additional boundary value problem (of linearised elasticity type), where we use the deformation given by $\tilde{\nabla \Cj(\Omega)}$ as Dirichlet boundary conditions.
\begin{lstlisting}[firstnumber = last]
def getExtension(gfX_bnd, freedofs, gfX_ext):
    u,v = VEC.TnT()
    aX_ext = BilinearForm(VEC)
    aX_ext += InnerProduct(grad(u)+grad(u).trans,grad(v))*dx+InnerProduct(u,v)*dx

    gfX_ext.Set(gfX_bnd)
    aX_ext.Assemble()

    r = gfX_bnd.vec.CreateVector()
    r.data = (-1) * aX_ext.mat * gfX_ext.vec

    gfX_ext.vec.data += aX_ext.mat.Inverse(freedofs=freedofs) * r
    
getExtension(gfX_bnd, VEC2.FreeDofs(), gfX) 
gfset.Set((0,0))
gfset.vec.data = gfset.vec - 1 * gfX.vec
\end{lstlisting}

The Newton algorithm reads as follows.
\begin{algorithm}[H] 
    \begin{algorithmic}[1]
        \State{{\bf Input:} domain $\Omega_0$, $n=0$, $N_{max}>0$, $\epsilon >0$}
        \State {\bf Output:} optimal shape $\Omega^*$ 
        \While{ $n\le N_{max}$ and $|\nabla \Cj(\Omega_n)|>\epsilon$}
            \State solve \eqref{eq_newtonEq} to get $\nabla \Cj(\Omega_n)$
            \State $\Omega_{n+1} \gets (\Id - \nabla \Cj(\Omega_n))(\Omega_n)$
            \State $n\gets n+1$ 
        \EndWhile
\caption{Newton algorithm}
\label{alg:newton}
\end{algorithmic}
\end{algorithm}

\subsection{Newton's method for PDE-constrained problems} \label{sec_newtonPDEConstr}

We consider the PDE-constrained model problem of Section \ref{sec_laplace} which is subject to the Poisson equation. The unregularised Newton system reads 
\ben \label{eq_shapeNewtonPDEconstr}
 D^2\Cj(\Omega)(V)(W) = - D\Cj(\Omega)(V) \quad \text{ for all } V\in H^1_0(\Omega). 
\een
In Subsection~\ref{sec_PDE_constr_shape_2nd} we discussed how the second order shape derivative can be evaluated along a fixed given direction. In this section, we want to assemble the whole shape Hessian and eventually solve a regularised version of \eqref{eq_shapeNewtonPDEconstr}. Recalling that $D\Cj(\Omega)(V) = \partial_s G_{V,0}(x)$ with $x = (0,0,u,p)$ we see that \eqref{eq:material_derivatrive_system} and \eqref{eq:diff_second_J2} lead to 
\ben \label{eq_newtonSystem}
\mathcal H \; \begin{pmatrix}
\tilde V \\ \partial_s u^0 \\ \partial_s p^0
         \end{pmatrix} = -\begin{pmatrix}
         \partial_t G_{V,W}(x) \\ 
    0\\
    0
     \end{pmatrix}.
\een
with
\ben
\mathcal H(x) = \begin{pmatrix}
    \partial_s\partial_t G_{V,W}(x)    &  \partial_u \partial_t G_{V,W}(x)   & \partial_p \partial_t  G_{V,W}(x) \\
    \partial_s \partial_u G_{V,W}(x)         &  \partial_u^2G_{V,W}(x)    & \partial_p\partial_u G_{V,W}(x) \\
    \partial_s \partial_p G_{V,W}(x)     &  \partial_u \partial_p G_{V,W}(x)    & 0 
         \end{pmatrix}.
\een
The component $\tilde V$ then represents the direction which we use for the shape Newton optimisation step. The matrix in \eqref{eq_newtonSystem} can be realised in \ngsolve{} by using a combined finite element space \texttt{X3} consisting of three components as follows:
\begin{lstlisting}[firstnumber = last]
X3 = FESpace([VEC, fes, fes])      
PHI, u1, p1= X3.TrialFunction()
PSI, uTest1, pTest1 = X3.TestFunction() (*@\label{lst_def_testFuncs}@*)

shapeHessLag3 = BilinearForm(X3)
shapeHessLag3 += G_pde_0.DiffShape(PHI).DiffShape(PSI)        #block (1,1)
shapeHessLag3 += G_pde_0.DiffShape(PSI).Diff(gfu,u1)         #block (1,2)
shapeHessLag3 += G_pde_0.DiffShape(PSI).Diff(gfp,p1)         #block (1,3)
shapeHessLag3 += G_pde_0.Diff(gfu, uTest1).DiffShape(PHI)    #block (2,1) 
shapeHessLag3 += (G_pde_0.Diff(gfu, uTest1)).Diff(gfu, u1)   #block (2,2)
shapeHessLag3 += (G_pde_0.Diff(gfu, uTest1)).Diff(gfp, p1)   #block (2,3)
shapeHessLag3 += G_pde_0.Diff(gfp,pTest1).DiffShape(PHI)     #block (3,1)
shapeHessLag3 += (G_pde_0.Diff(gfp,pTest1)).Diff(gfu, u1)    #block (3,2)
\end{lstlisting}
The right hand side of \eqref{eq_newtonSystem} can be defined as follows: 
\begin{lstlisting}[firstnumber = last]
shapeGradLag3 = LinearForm(X3)
shapeGradLag3 += (-1) * G_pde_0.DiffShape(PSI)
\end{lstlisting}
Recall that the system \eqref{eq_newtonEq} has a nontrivial kernel as discussed in Section \ref{sec_newtonUnConstr}. This problem can be circumvented by proceeding like in the unconstrained case. We add a regularisation only on the boundary,
\begin{lstlisting}[firstnumber = last]
delta = 1
shapeHessLag3 += delta * InnerProduct(PHI, specialcf.tangential(2)) * InnerProduct(PSI, specialcf.tangential(2))*ds
\end{lstlisting}
and exclude the interior degrees of freedom in the first row and column of the 3$\times$3 block system. This can be realised by setting Dirichlet boundary conditions for the interior degrees of freedom, i.e. by dealing with the free degrees of freedom,
\begin{lstlisting}[firstnumber = last]
# copy of VEC with Dirichlet boundary conditions on whole boundary:
VEC2 = VectorH1(mesh, order = 1, dirichlet = ".*")    
freeDofsCombined = BitArray(VEC2.ndof + 2*fes.ndof)
for i in range(VEC2.ndof):
    freeDofsCombined[i] = not VEC2.FreeDofs()[i]
for i in range(fes.ndof):
    freeDofsCombined[VEC2.ndof+i] = fes.FreeDofs()[i]
    freeDofsCombined[VEC2.ndof+fes.ndof+i] = fes.FreeDofs()[i]
\end{lstlisting}
and solving the regularised system using these free dofs:
\begin{lstlisting}[firstnumber = last]
gfCombined3 = GridFunction(X3)
shapeHessLag3.Assemble()
shapeGradLag3.Assemble()
gfCombined3.vec.data = shapeHessLag3.mat.Inverse(freedofs=freeDofsCombined, inverse="umfpack") * shapeGradLag3.vec
\end{lstlisting}
The newton direction is then given as the first of the three components of the obtained solution.
\begin{lstlisting}[firstnumber = last]
Vtilde_bnd = GridFunction(VEC)
Vtilde = GridFunction(VEC)
Vtilde_bnd.vec.data = gfCombined3.components[0].vec
getExtension(Vtilde_bnd, VEC2.FreeDofs(), Vtilde) 

gfset.vec.data = gfset.vec + 1 * Vtilde.vec
\end{lstlisting}

\section{Numerical Experiments}
In this section we first verify the copmuted shape derivatives by performing a Taylor test, and then apply the automated shape differentiation and the numerical algorithms introduced in the preceding sections in numerical examples.

\subsection{Code verification}
We verify the expressions that we obtained in a semi-automatic or fully automatic way for the first and second order shape derivatives by looking at the Taylor expansions of the perturbed shape functionals. We illustrate our findings in two examples in $\VR^2$. On the one hand, we consider a shape function as introduced in \eqref{eq:shape_func_f} with an additional boundary integral as in \eqref{eq:shape_func_f_bnd}, henceforth denoted by $\mathcal J_1$; on the other hand, we consider the PDE-constrained shape optimisation problem defined by \eqref{eq_modelPDEopti}, the reduced form of which will be denoted by $\mathcal J_2(\Omega)$. More precisely, we consider
\begin{align}
    \mathcal J_1(\Omega) &= \int_\Omega f(x) \; \mbox dx + \int_{\partial \Omega} f(x) \; \mbox ds, \\
    \mathcal J_2(\Omega) &= \int_\Omega |u_\Omega-u_d|^2 \; \mbox dx \; \mbox{ where } u_\Omega \mbox{ solves \eqref{eq_modelPDEopti_PDE}. }
\end{align}
In the case of $\mathcal J_1$, we used the function $$f(x_1, x_2) = \left(0.5+\sqrt{x_1^2 +x_2^2}\right)^2 \left(0.5-\sqrt{x_1^2+x_2^2} \right)^2$$ and for $\mathcal J_2$, we used $u_d(x_1,x_2) = x_1(1-x_1)x_2(1-x_2)$ and $f(x_1, x_2) = 2 x_2(1-x_2)+2x_1(1-x_1)$ for the function $f$ in the PDE constraint \eqref{eq_modelPDEopti_PDE}.

\begin{figure}[H]
    \begin{center}
        \begin{tabular}{cc}
            \includegraphics[width=.25\textwidth]{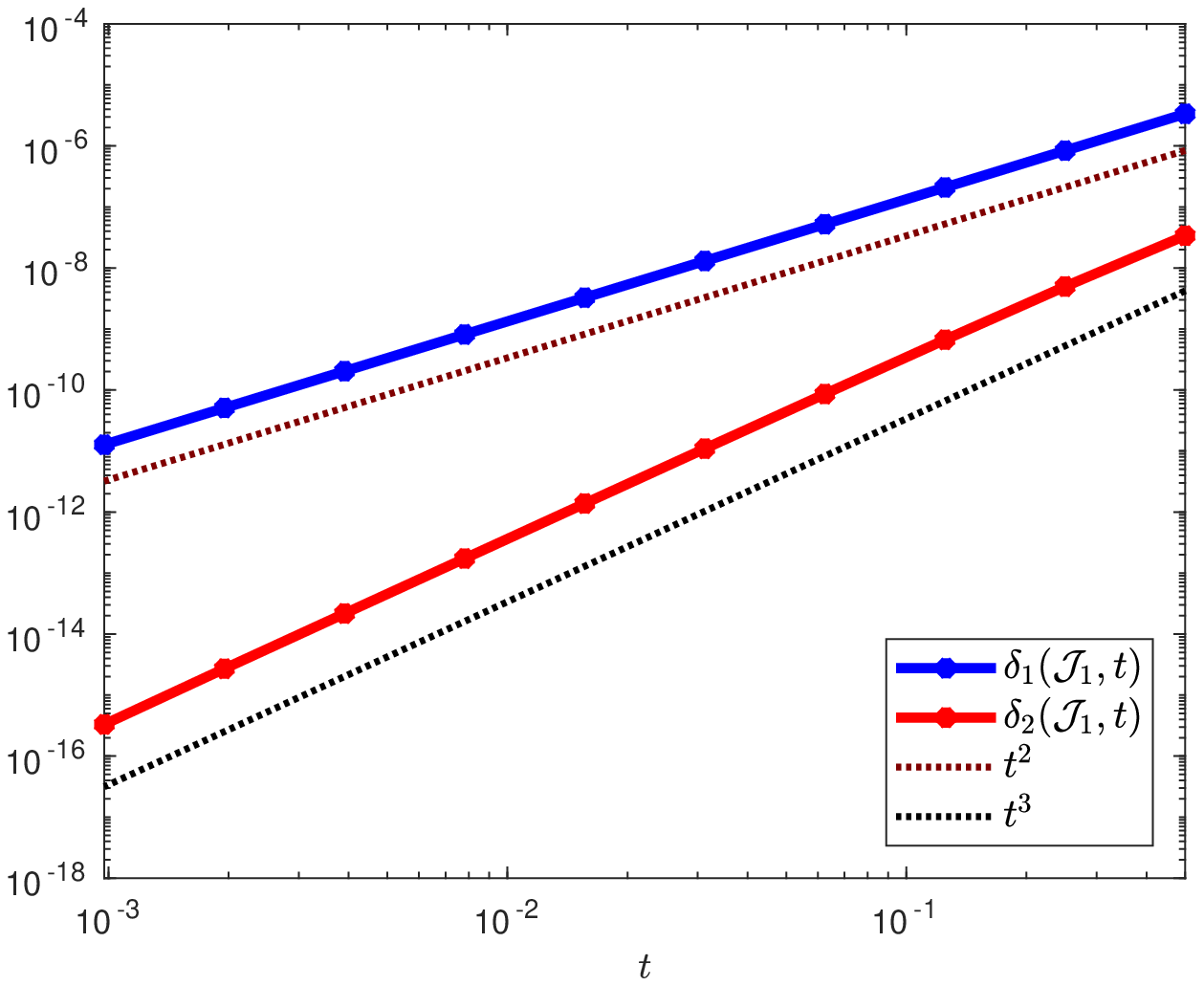} &
            \includegraphics[width=.25\textwidth]{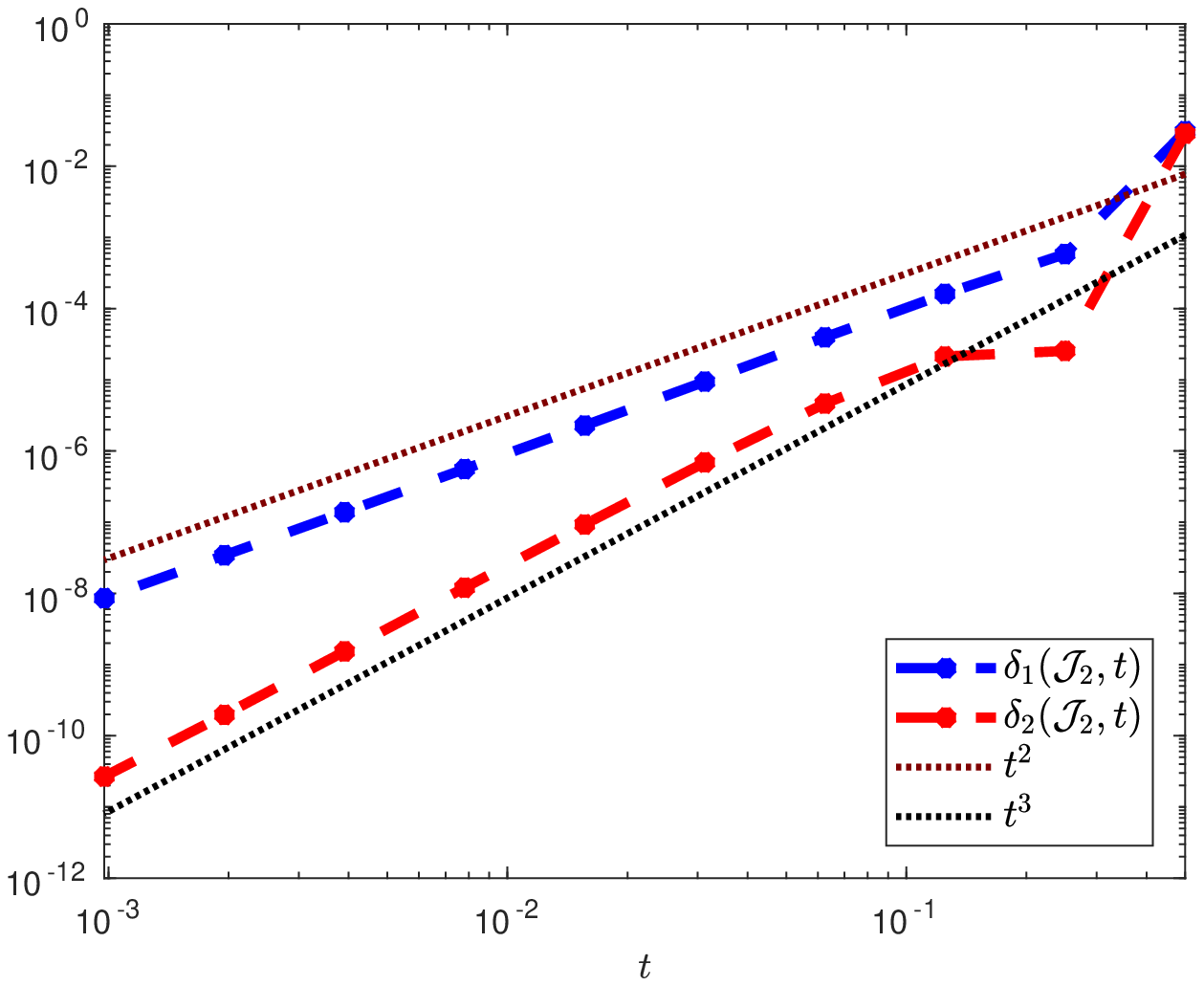} \\
            (a) & (b)
        \end{tabular}
        \caption{Taylor test for functions $\mathcal J_1$ and $\mathcal J_2$.}
        \label{fig_taylortests}
    \end{center}
\end{figure}

For the test of the first order shape derivatives $D \mathcal J_i(\Omega)(V)$ we choose a fixed shape $\Omega$ and a vector field $V \in C^{0,1}(\VR^2)^2$ and observe the quantity
\ben
    \delta_1(\mathcal J_i, t) := \left\lvert \mathcal J_i( (\Id + t V)(\Omega) ) - \mathcal J_i(\Omega) - t  \; D\mathcal J_i(\Omega)(V) \right\rvert,
\een
for $t \searrow 0$. Likewise, for the second order shape derivative, we consider the remainder
\begin{align*}
    \delta_2(\mathcal J_i, t) := \left\lvert \vphantom{\frac{1}{2}}\right. \mathcal J_i( (\Id + t V )(\Omega) ) &- \mathcal J_i(\Omega) - t  \; D\mathcal J_i(\Omega)(V)  \\
    &- \frac{1}{2} t^2 D^2 \mathcal J_i(\Omega)(V)(V) \left. \vphantom{\frac{1}{2}} \right\rvert
\end{align*}
as $t \searrow 0$. By the definition of first and second order shape derivatives, it must hold that
\ben
    \delta_1(\mathcal J_i, t) = \mathcal O(t^2) \quad \mbox{ and } \quad \delta_2(\mathcal J_i, t) = \mathcal O(t^3) \quad \mbox{ as } t \searrow 0.
\een
This behavior can be observed in Figure \ref{fig_taylortests}(a) for $\mathcal J_1$ and in Figure \ref{fig_taylortests}(b) for $\mathcal J_2$, where we used $V(x_1, x_2) = (x_1^2 x_2 e^{x_2}, x_2^2 x_1 e^{x_1})$ in both cases.

The experiments for shape function $\mathcal J_1$ was conducted on a mesh consisting of 13662 vertices, 26946 elements and with polynomial order 2 (resulting in 54269 degrees of freedom), and the experiment for $\mathcal J_2$ with 95556 vertices and 190062 elements and polynomial degree 1 (95556 degrees of freedom).
We conducted these experiments for a number of different problems with different vector fields $V$, in particular with different PDE constraints and boundary conditions, and obtained similar results in all instances provided a sufficiently fine mesh was used.

 \subsection{A first shape optimisation problem}
 In this section, we revisit problem \eqref{eq:shape_func_f} introduced in Section \ref{sec_intf}, i.e. the problem of finding a shape $\Omega$ such that the cost function $\mathcal J(\Omega) =  \int_\Omega f(x) \; \mbox dx$ is minimised.
 
 \subsubsection{First order methods}
\begin{figure}
    \begin{center}
        \includegraphics[width=.4\textwidth]{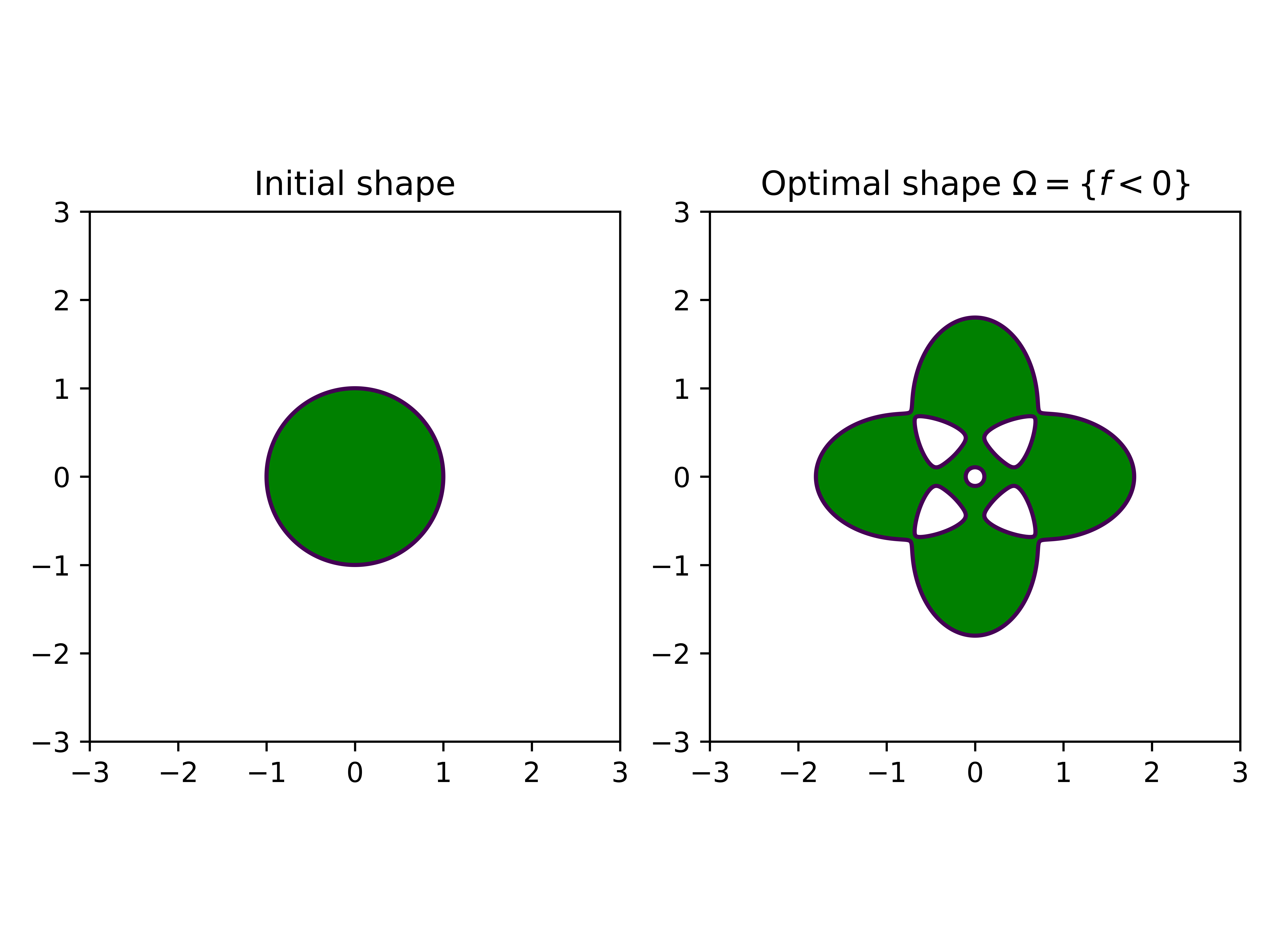} \vspace{-8mm}
    \end{center}
    \caption{Initial domain $\Omega_0$ and optimal domain $\Omega^*$ for problem \eqref{eq:shape_func_f} with $f$ chosen according to \eqref{eq_f_clover}.}
    \label{fig_intf_init_opt}
 \end{figure}

 We illustrate our first order methods in a problem which was also considered in \cite{a_IGSTWE_2018a} and reproduce the results obtained there.  We choose the function
 \begin{align}
    \begin{aligned} \label{eq_f_clover}
    f(x_1,& x_2) = \left(\sqrt{(x_1 - a)^2 + b x_2^2} - 1\right) \left(\sqrt{(x_1 + a)^2 + b x_2^2} - 1\right) \\
    & \cdot  \left(\sqrt{b x_1^2 + (x_2 - a)^2} - 1\right)  \left(\sqrt{b x_1^2 + (x_2 + a)^2} - 1\right) - \epsilon
    \end{aligned}
 \end{align}
 with $a=\frac{4}{5}$, $b=2$ and $\epsilon=0.001$. Recall that the optimal shape is given by $\{(x_1,x_2) \in \VR^2: f(x_1, x_2)<0 \}$ which is depicted in Figure \ref{fig_intf_init_opt} (right). We start our optimisation algorithm with the unit disk, $\Omega^0 = B_1(0)$ as an initial design. Note that the optimal design cannot be reached by means of shape optimisation using boundary perturbations. However, we expect the outer curve of the optimal shape to be reached very closely.
 
\begin{figure}
    \begin{tabular}{cc}
        \includegraphics[width=.2\textwidth, trim = 150 0 150 0, clip]{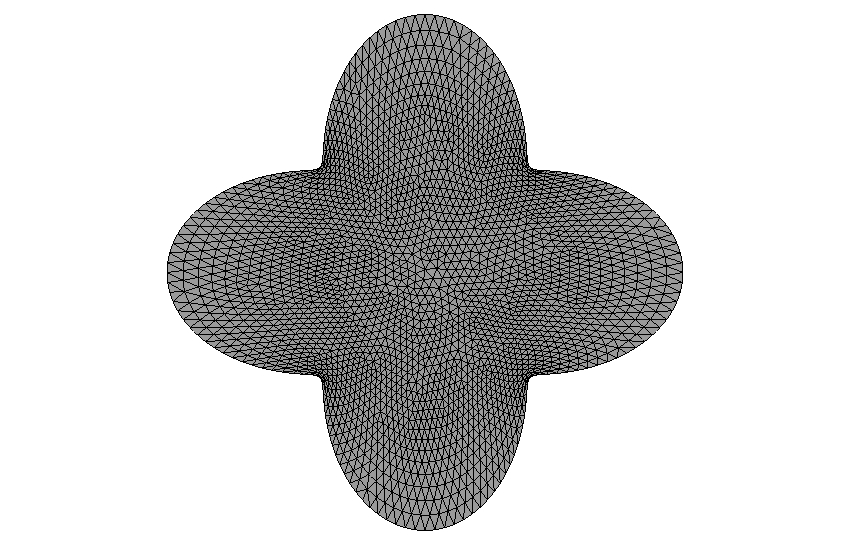} &
        \includegraphics[width=.25\textwidth]{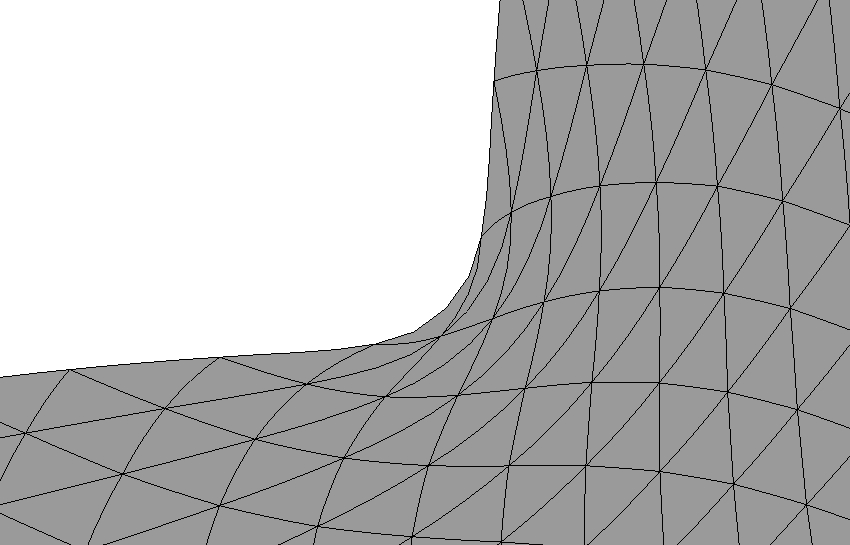} 
    \end{tabular}
    \caption{Results of problem \eqref{eq:shape_func_f} with $f$ as in \eqref{eq_f_clover} and the shape gradient associated to the $H^1$ inner product \eqref{eq_shapeGrad_H1}.}
    \label{fig_clover_results_H1}
 \end{figure}

\begin{figure}
    \begin{tabular}{cc}
        \includegraphics[width=.2\textwidth, trim = 150 0 150 0, clip]{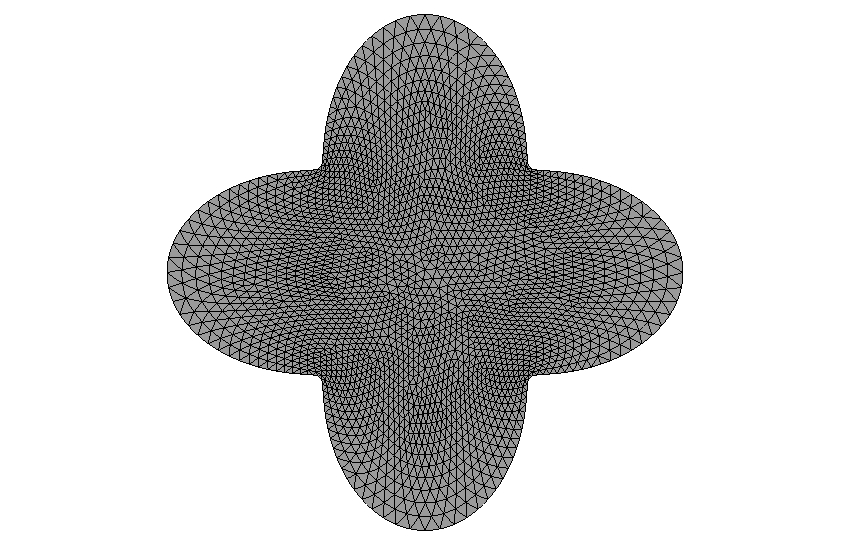} &
        \includegraphics[width=.25\textwidth]{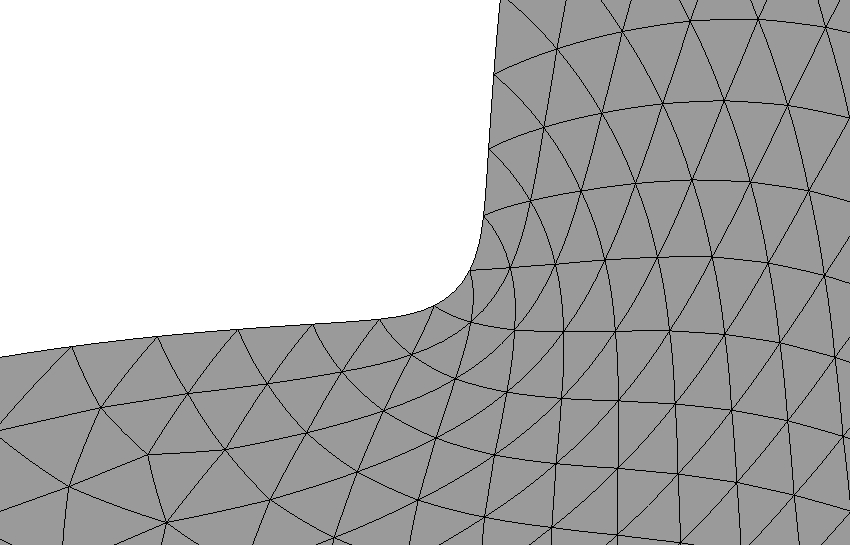} 
    \end{tabular}
    \caption{Results of problem \eqref{eq:shape_func_f} with $f$ as in \eqref{eq_f_clover} and the shape gradient associated to the elasticity bilinear form \eqref{eq_shapeGrad_ela}.}
    \label{fig_clover_results_ela}
 \end{figure}

\begin{figure}
    \begin{tabular}{cc}
        \includegraphics[width=.2\textwidth, trim = 150 0 150 0, clip]{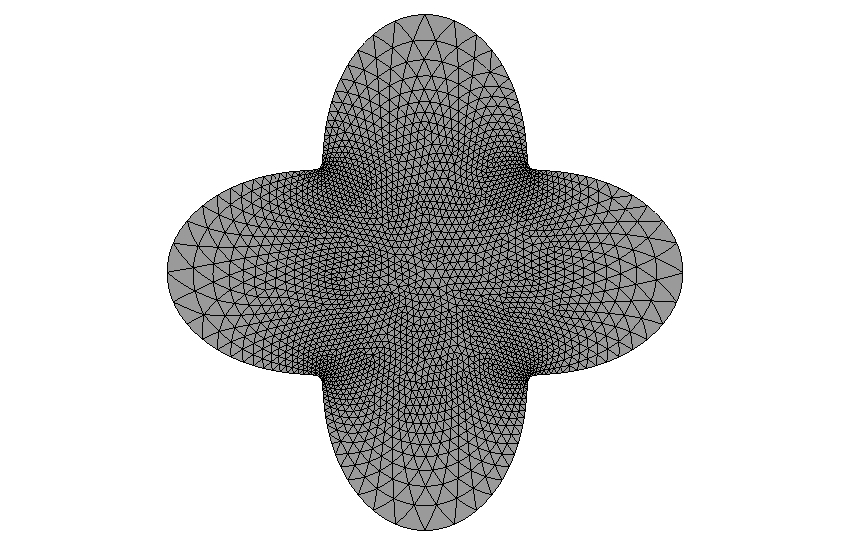} &
        \includegraphics[width=.25\textwidth]{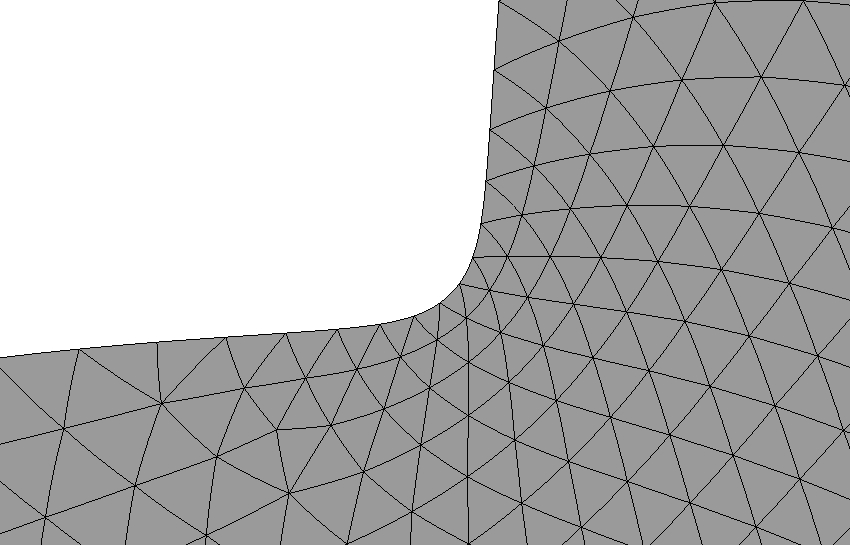} 
    \end{tabular}
    \caption{Results of problem \eqref{eq:shape_func_f} with $f$ as in \eqref{eq_f_clover} and the shape gradient associated to the elasticity bilinear form with Cauchy-Riemann term \eqref{eq_shapeGrad_ela_CR}.}
    \label{fig_clover_results_ela_CR}
 \end{figure}

 We apply Algorithm \ref{alg:gradient} with the shape gradient $\nabla \mathcal J$ associated to the $H^1$ inner product \eqref{eq_shapeGrad_H1}, to the bilinear form of linearised elasticity \eqref{eq_shapeGrad_ela} and including the additional Cauchy-Riemann term \eqref{eq_shapeGrad_ela_CR}. We chose the algorithmic parameters $\gamma = 1e-4$, $\epsilon = 1e-7$, a mesh consisting of 2522 vertices and 4886 elements and a globally continuous vector-valued finite element space \texttt{VEC} of order 3. The results can be seen in Figures \ref{fig_clover_results_H1}, \ref{fig_clover_results_ela} and \ref{fig_clover_results_ela_CR}, respectively.


 \subsubsection{Second order method}
\begin{figure}
    \begin{tabular}{ccc}
        \includegraphics[width=.15\textwidth, trim = 150 0 150 0, clip]{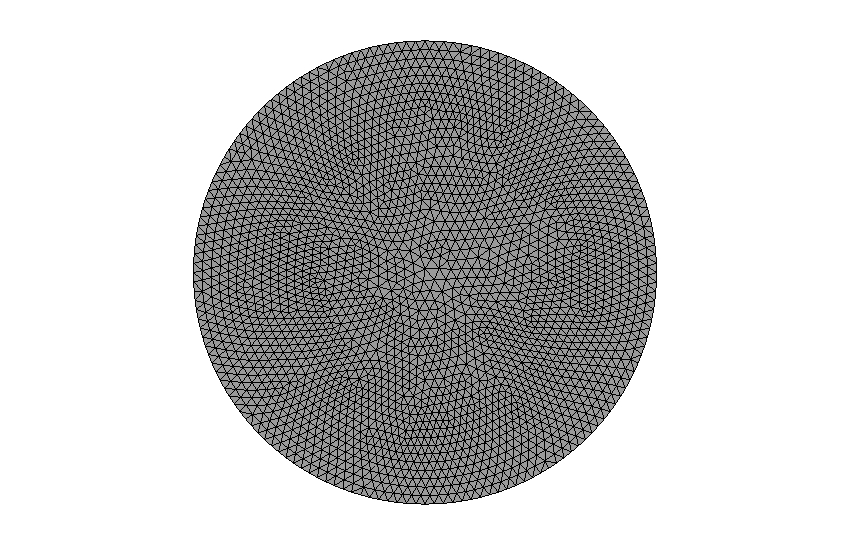}&
        \includegraphics[width=.15\textwidth, trim = 150 0 150 0, clip]{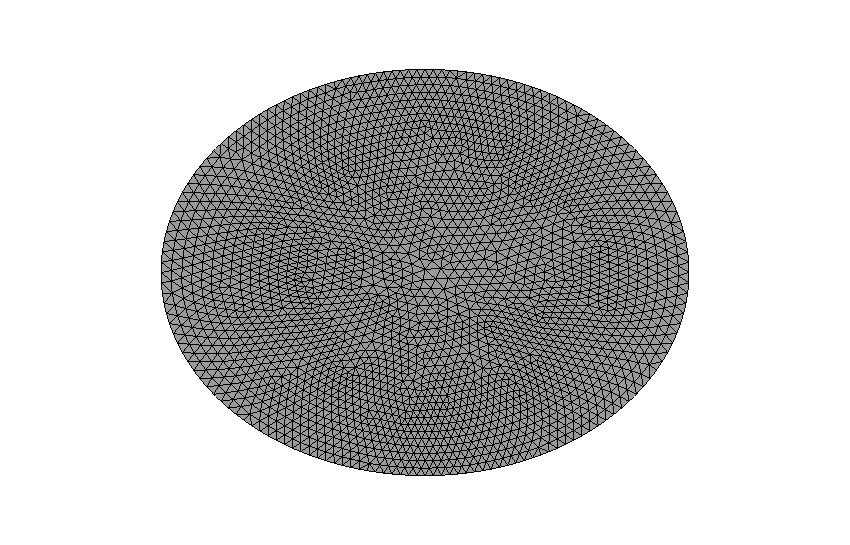}&
        \includegraphics[width=.17\textwidth, trim = 20 0 15 0, clip]{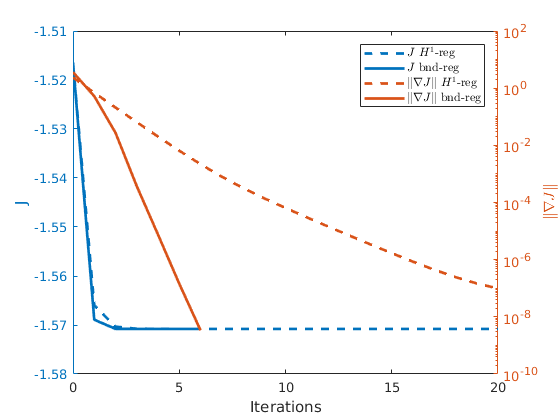}
    \end{tabular}
    \caption{Numerical results for problem \eqref{eq_regBnd} with $f$ as in \eqref{eq_f_ellipse} using second order method. Left: Initial design. Center: Optimised design after six iterations  using \eqref{eq_newtonEq}/\eqref{eq_regBnd}. Right: Objective value $\mathcal J$ and norm of shape gradient $\| \nabla \Cj(\Omega)\|$ in the course of second order optimisation using \eqref{eq_regH1} with $\delta = 0.5$ and \eqref{eq_regBnd} with $\delta = 100$.}
    \label{fig_ellipseSecondOrder}
 \end{figure}

 Since Newton's method converges quadratically only in a neighborhood of the optimal solution, we choose a simpler optimal design here. We choose
 \ben \label{eq_f_ellipse}
    f(x_1, x_2) = \frac{x^2}{a^2} + \frac{y^2}{b^2} - 1
 \een
 which yields an ellipse with the lengths of the two semi-axes $a$ and $b$. We choose $a= 1.3$ and $b=1/a$ and again start the optimisation with the unit disk as initial shape. Figure \ref{fig_ellipseSecondOrder} shows the initial and optimised design after only six iterations of Algorithm \ref{alg:newton} with $(\cdot,\cdot)_H$ chosen as in \eqref{eq_regBnd} with $\delta= 100$. A comparison of the convergence histories between the choice \eqref{eq_regBnd} with $\delta = 100$ and \eqref{eq_regH1} with $\delta = 0.5$ is shown in the right picture of Figure \ref{fig_ellipseSecondOrder}. In both cases, we tested a range of different values for $\delta$ and compared the convergence histories for the values which yielded the fastest convergence. 
 The experiments were conducted on a finite element mesh consisting of 2522 nodes and 4886 triangular elements with a finite element space \texttt{VEC} of order 3, with the algorithmic parameter $\epsilon = 10^{-7}$.
 
 \subsection{Shape optimisation subject to the Poisson equation}
 In this section, we revisit the model problem introduced in Section \ref{sec_laplace} with $f(x_1, x_2) = 2 x_2(1-x_2)+2x_1(1-x_1)$ and $u_d(x_1,x_2) = x_1(1-x_1)x_2(1-x_2)$. Note that the data is chosen in such a way that, for $\Omega^* = (0,1)^2$ it holds $\mathcal J(\Omega^*) = 0$ and thus $\Omega^*$ is a global minimiser of $\mathcal J$. We show results obtained by first and second order shape optimisation methods exploiting automated differentiation. 
 
\begin{figure}
    \begin{tabular}{cc}
        \includegraphics[width=.25\textwidth]{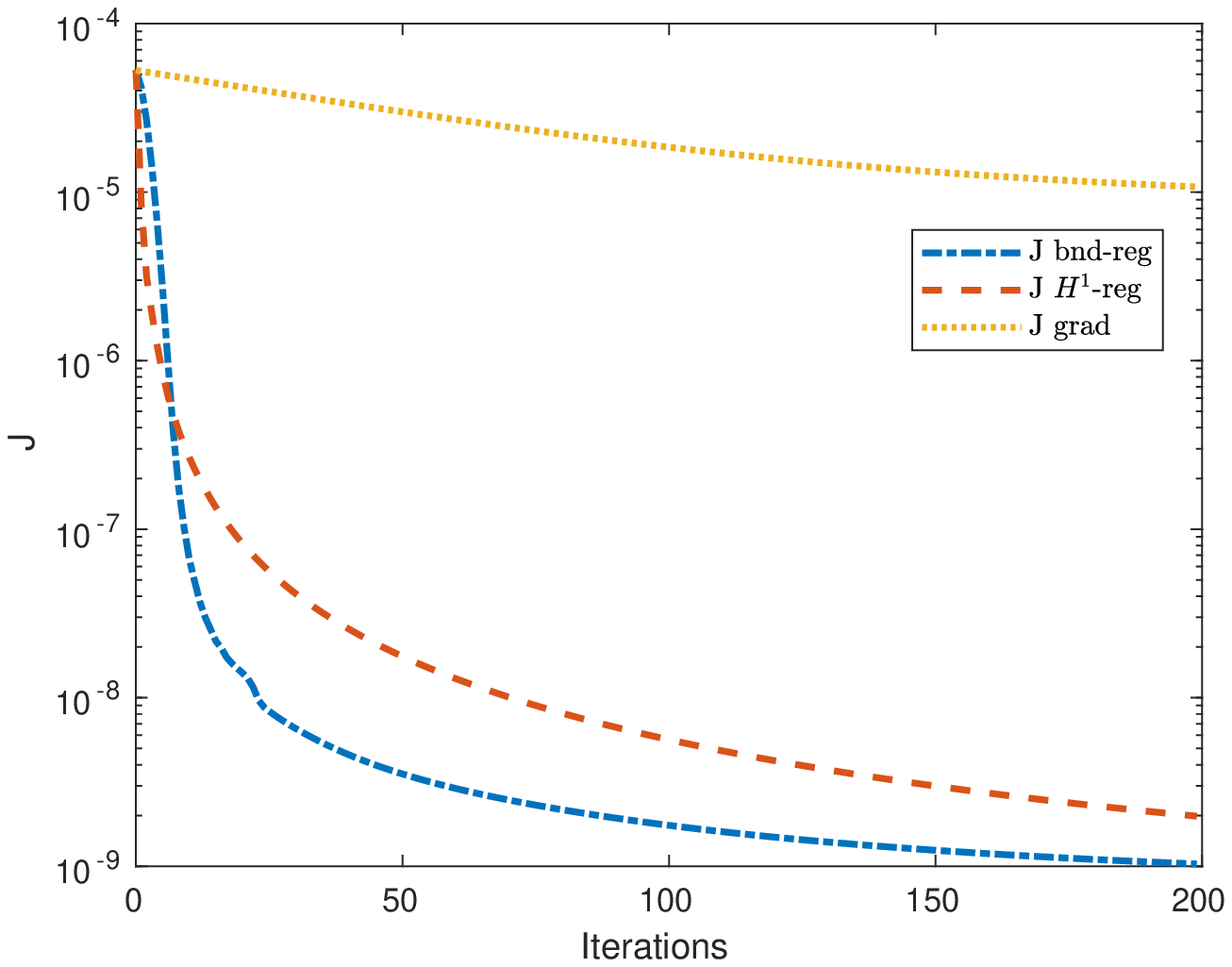}&
        \includegraphics[width=.25\textwidth]{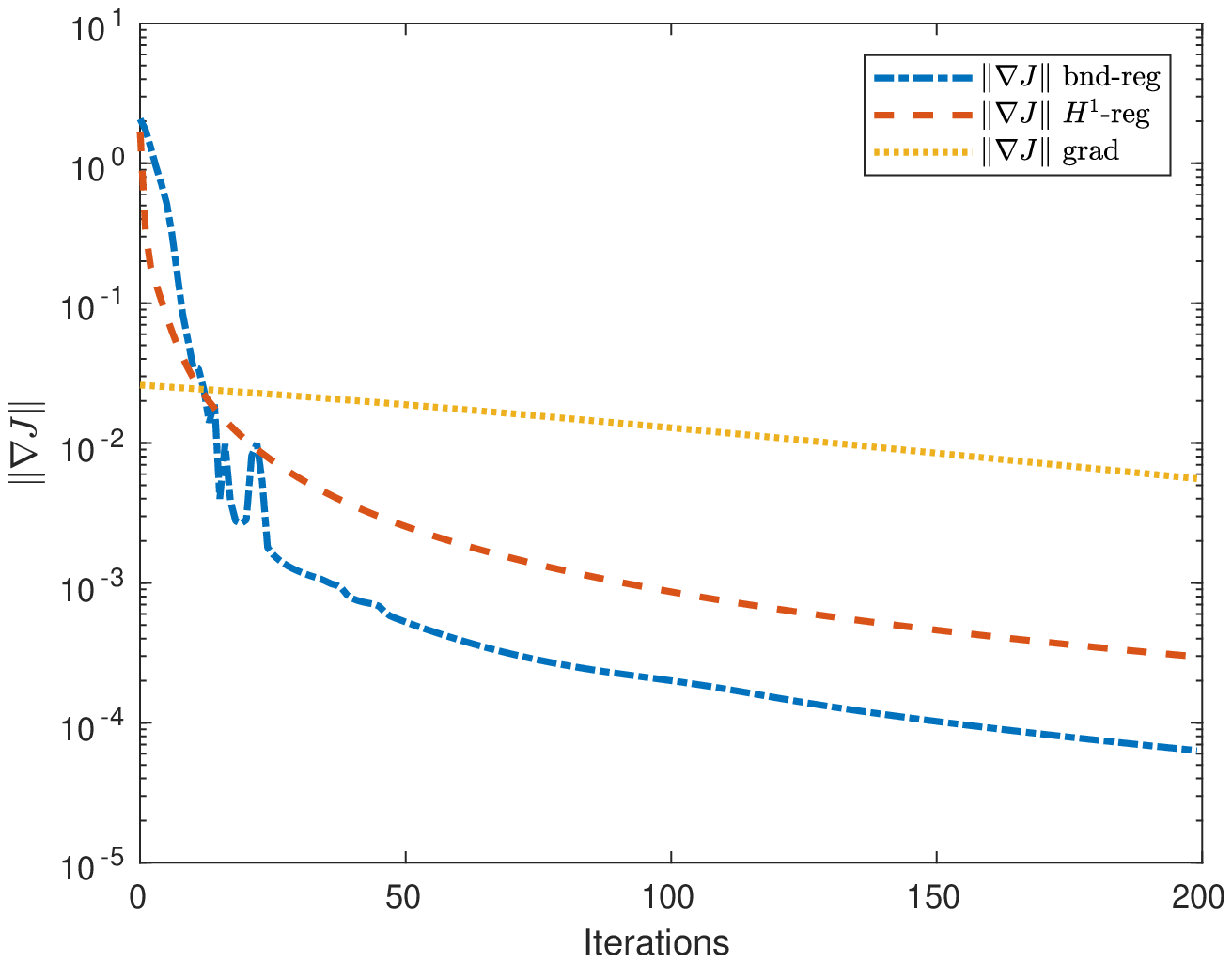} \\
        (a) & (b)
    \end{tabular}
    \caption{Convergence behaviour for shape optimisation problem \eqref{eq_modelPDEopti} with proposed regularisation strategies \eqref{eq_regBnd} and \eqref{eq_regH1} as well as first order method with constant step size $\alpha = 1$. (a) Behaviour of objective function $\Cj$. (b) Behaviour of norm of shape gradient $\| \nabla \Cj(\Omega)\|$.}
    \label{fig_historyEx1}
 \end{figure}

 We ran the optimisation algorithm in three versions. On the one hand, we applied a first order method with constant step size $\alpha = 1$. On the other hand, we applied two second order methods with the two different regularisation strategies for the shape Hessian in \eqref{eq_newtonEq} introduced in \eqref{eq_regH1} and \eqref{eq_regBnd}. We chose the regularisation parameters $\delta$ empirically such that the method performs as well as possible. In the case of \eqref{eq_regH1} we chose $\delta = 0.001$ and in the case of \eqref{eq_regBnd} $\delta = 1$. The experiments were conducted on a finite element mesh consisting of 4886 elements with 2522 vertices and polynomial degree 1.  In Figure \ref{fig_historyEx1}, we can observe the decrease of the objective function as well as of the norm of the shape gradient over 200 iterations for these three algorithmic settings.
 
 Figure \ref{fig_historyEx1} shows the initial design as well as the design after 200 iterations of the second order method with regularisation strategy \eqref{eq_regBnd}. Note that the improved design is very close to $\Omega^* = (0,1)^2$, which is a global solution. The initial design was chosen as the disk of radius $\frac{1}{2}$ centered at the point $\left(\frac{1}{2}, \frac{1}{2}\right)^\top$. The objective value was reduced from $5.297 \cdot 10^{-5}$ to $1.0317 \cdot 10^{-9}$.
 
 \begin{figure}
    \begin{tabular}{cc}
        \includegraphics[width=.23\textwidth, trim = 40 0 50 0, clip]{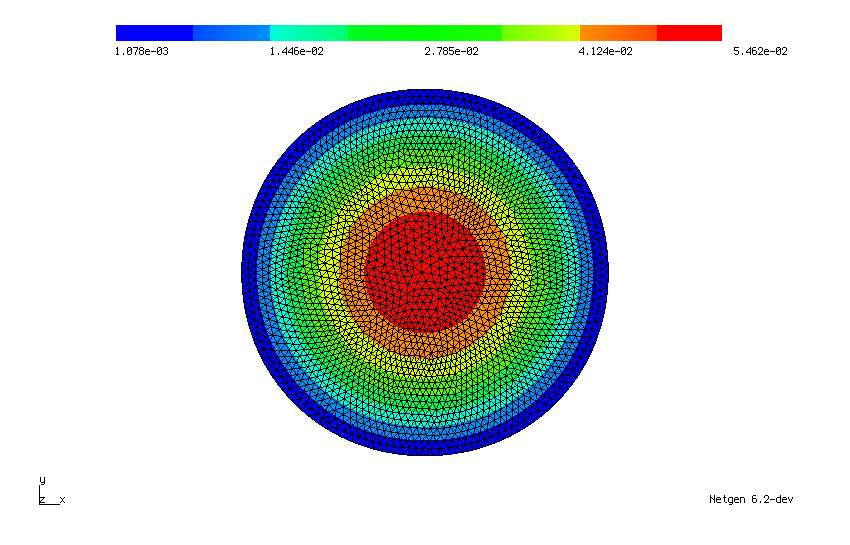}&
        \includegraphics[width=.23\textwidth, trim = 40 0 50 0, clip]{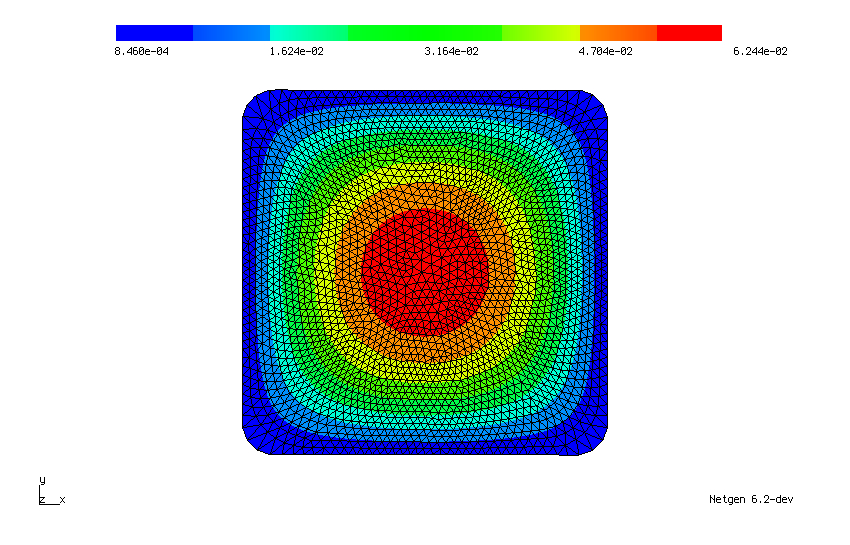}
    \end{tabular}
    \caption{Shape optimisation for problem \eqref{eq_modelPDEopti}. Left: Initial design. Right: Improved design after 200 iterations of second order algorithm with regularisation as proposed in \eqref{eq_regBnd}. Objective value was reduced from $5.297 \cdot 10^{-5}$ to $1.0317 \cdot 10^{-9}$. Color shows solution of constraining PDE \eqref{eq_modelPDEopti_PDE}.}
    \label{fig_historyEx2}
 \end{figure}
 
 \subsection{Nonlinear elasticity}
\begin{figure}
    \begin{tabular}{ccc}
        \includegraphics[width=.15\textwidth, trim = 180 0 180 0, clip]{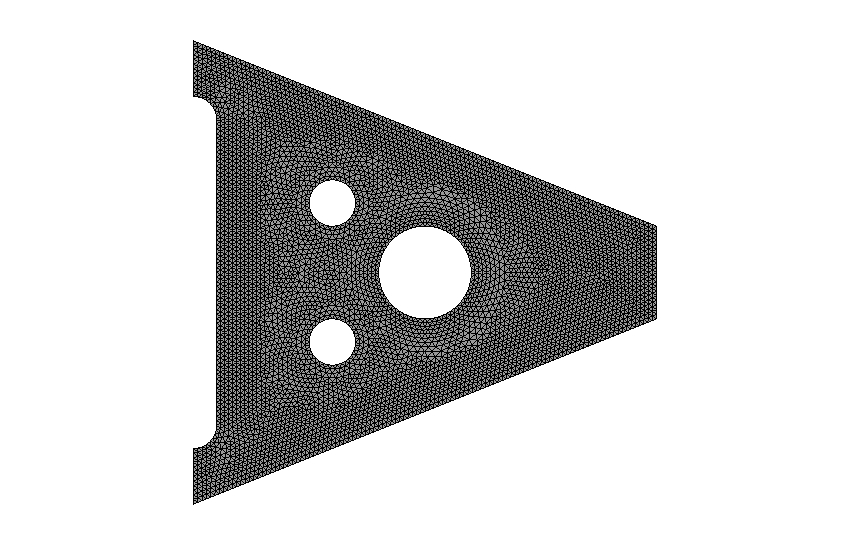}&
        \includegraphics[width=.15\textwidth, trim = 180 0 180 0, clip]{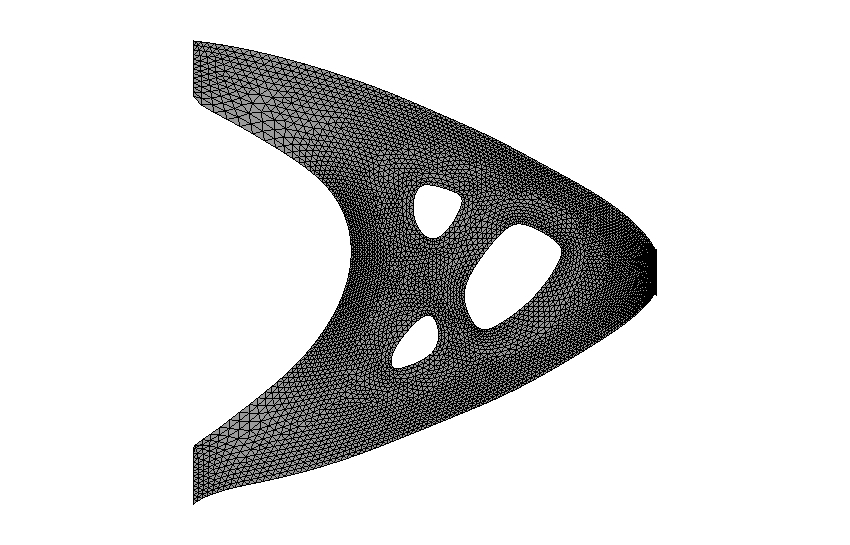}&
        \includegraphics[width=.15\textwidth, trim = 180 0 180 0, clip]{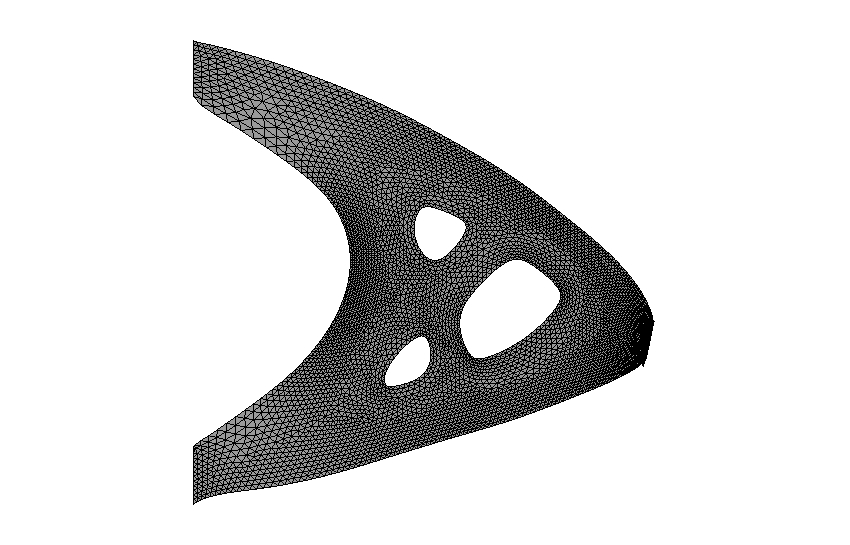}\\
        (a) & (b) & (c)
    \end{tabular}
    \caption{Initial and optimised geometry of cantilever under vertical force on right hand side using St. Venant--Kirchhoff model in nonlinear elasticity. (a) Initial geometry. (b) Optimised geometry (reference configuration). (c) Optimised geometry (deformed configuration).}
    \label{fig_elaNL}
 \end{figure}

 Here, we illustrate the applicability of the automated shape differentiation and optimisation in the more realistic and more complicated setting of nonlinear elasticity in two space dimensions using a Saint Venant–Kirchhoff material with Young's modulus $E=1000$ and Poisson ratio $\nu = 0.3$. We consider a two-dimensional cantilever which is clamped on the upper and lower left parts of the boundary, $\Gamma_l^1 = \{0\} \times (0.88, 1)$ and $\Gamma_l^2 = \{0\} \times (0, 0.12)$, respectively, and is subject to a surface force $g_N = (0,-100)^\top$ on $\Gamma_r= \{1\} \times (0.45,0.55)$. The initial geometry with 3 holes is depicted in Figure \ref{fig_elaNL} (a). Let $\Gamma_l := \Gamma_l^1 \cup \Gamma_l^2$ and $H_{\Gamma_l}^1(\Omega)^2$ the subspace of $H^1(\Omega)^2$ with vanishing trace on $\Gamma_l$. The displacement $u\in H_{\Gamma_l}^1(\Omega)^2$ under the surface force $g_N$ is given as the solution to the boundary value problem
 \ben \label{eq_nlEla}
    \int_\Omega S(u) : \nabla v \; \mbox dx = \int_{\Gamma_r} g_N \cdot v \; \mbox ds 
 \een 
 for all $v \in H_{\Gamma_l}^1(\Omega)^2$. Here, $S(u)$ denotes the Saint Venant--Kirchhoff stress tensor
 \ben
    S(u) = (I_2+\nabla u)  \left[\lambda \mbox{Tr}\left( \frac{1}{2} (C(u)-I_2) \right)I_2 + \mu(C(u)-I_2) \right],
 \een
 where $C(u) = (I_2 + \nabla u)^\top (I_2 + \nabla u)$ and $I_2$ is the identity matrix, see also \cite[Sec. 8]{AllaireJouveToader2004}, and $\lambda$ and $\mu$ denote the Lam\'e constants,
 \ben
    \lambda = \frac{E \nu }{(1+\nu)(1-2\nu)},\qquad \mu = \frac{E}{2 (1+\nu)}.
 \een   
 We minimise the functional 
 \ben 
    J(\Omega, u) = \int_\Omega S(u) : \nabla u \; \mbox dx   +\alpha \int_\Omega 1 \; \mbox dx
 \een
 with $\alpha = 2.5$ subject to \eqref{eq_nlEla} which amounts to maximising the structure's stiffness while bounding the allowed amount of material used.
 
 We remark that the well-posedness of \eqref{eq_nlEla} is not clear, see also the discussion in \cite[Sec. 8]{AllaireJouveToader2004}. Nevertheless, application of the automated shape differentiation and optimisation yields a significant improvement of the initial design.
The highly nonlinear PDE constraint \eqref{eq_nlEla} is solved by Newton's method. In order to have good starting values, a load stepping strategy is employed, i.e., the load on the right hand side is gradually increased, the PDE is solved and the solution is used as an initial guess for the next load step. This is repeated until the full load is applied.
With these ingredients at hand, Algorithm \ref{alg:gradient} (i.e. code lines \ref{lst_algo1_l1}--\ref{lst_algo1_lN}) can be run. We chose the algorithmic parameters \texttt{alpha = 0.1} (as an initial value), \texttt{alpha\_incr\_factor = 1} (i.e. no increase), \texttt{gamma = 1e-4} and \texttt{epsilon = 1e-7}. Moreover, we used \eqref{eq_shapeGrad_H1} with an additional Cauchy-Riemann term as in \eqref{eq_shapeGrad_ela_CR} with weight $\gamma_{CR} = 10$. The objective value was reduced from $3.125$ to $2.635$ (volume term from $1.290$ to $1.096$) in 15 iterations of Algorithm \ref{alg:gradient}. The results were obtained on a mesh consisting of 10614 elements and 5540 vertices using piecewise linear, globally continuous finite elements.

 \subsection{Helmholtz equation}
\begin{figure}
        \begin{tabular}{ccc}
            \includegraphics[width=.18\textwidth, trim = 250 0 250 0, clip]{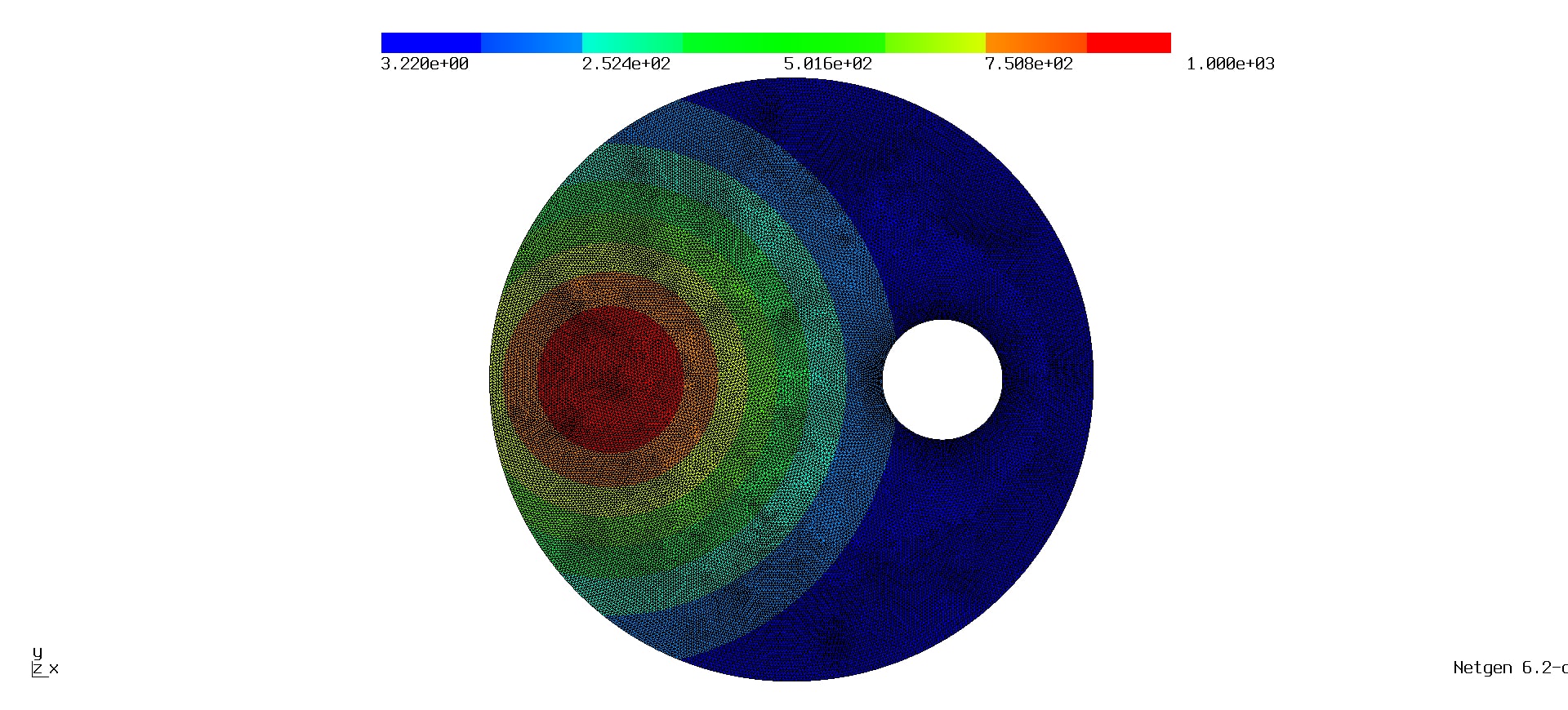} &
            \includegraphics[width=.13\textwidth, trim = 400 0 400 0, clip]{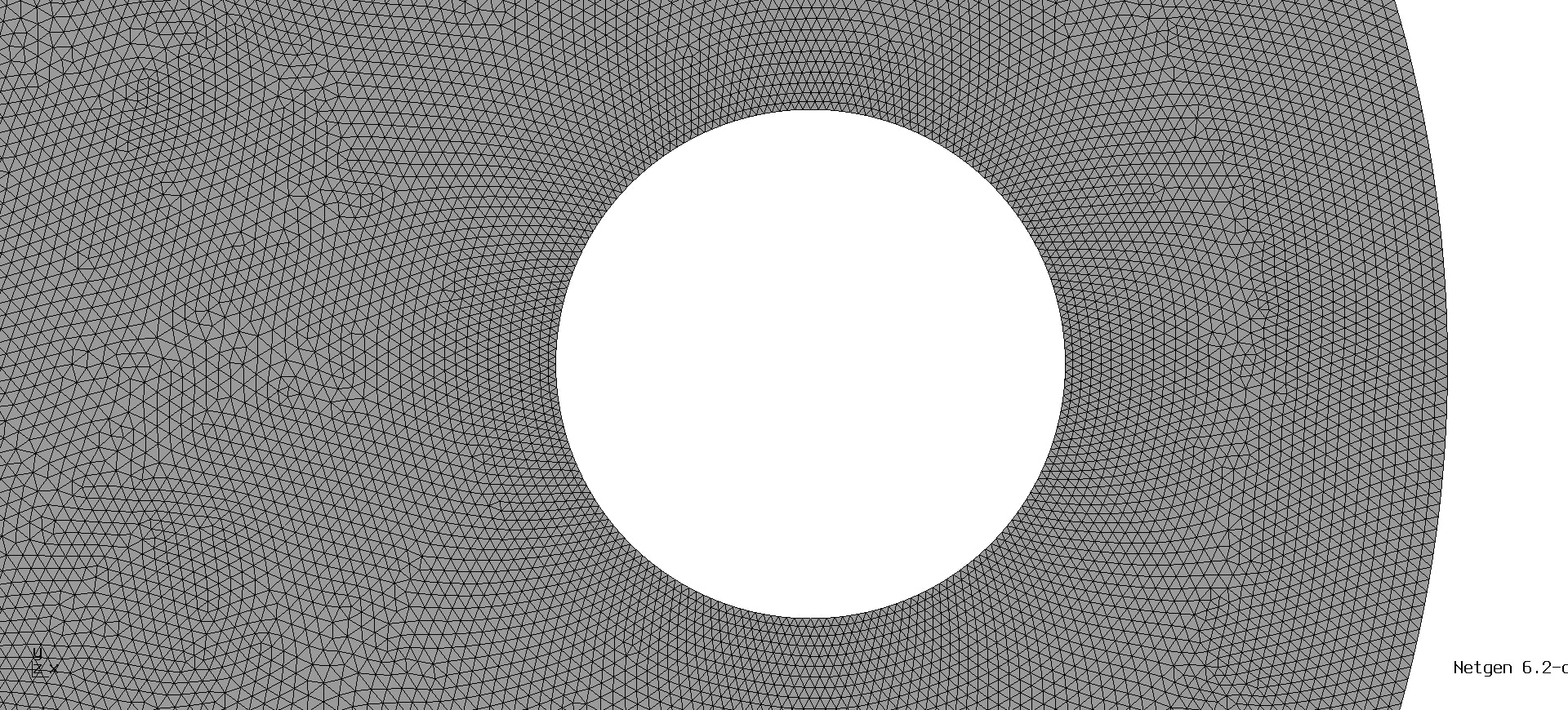} &
            \includegraphics[width=.13\textwidth, trim = 400 0 400 0, clip]{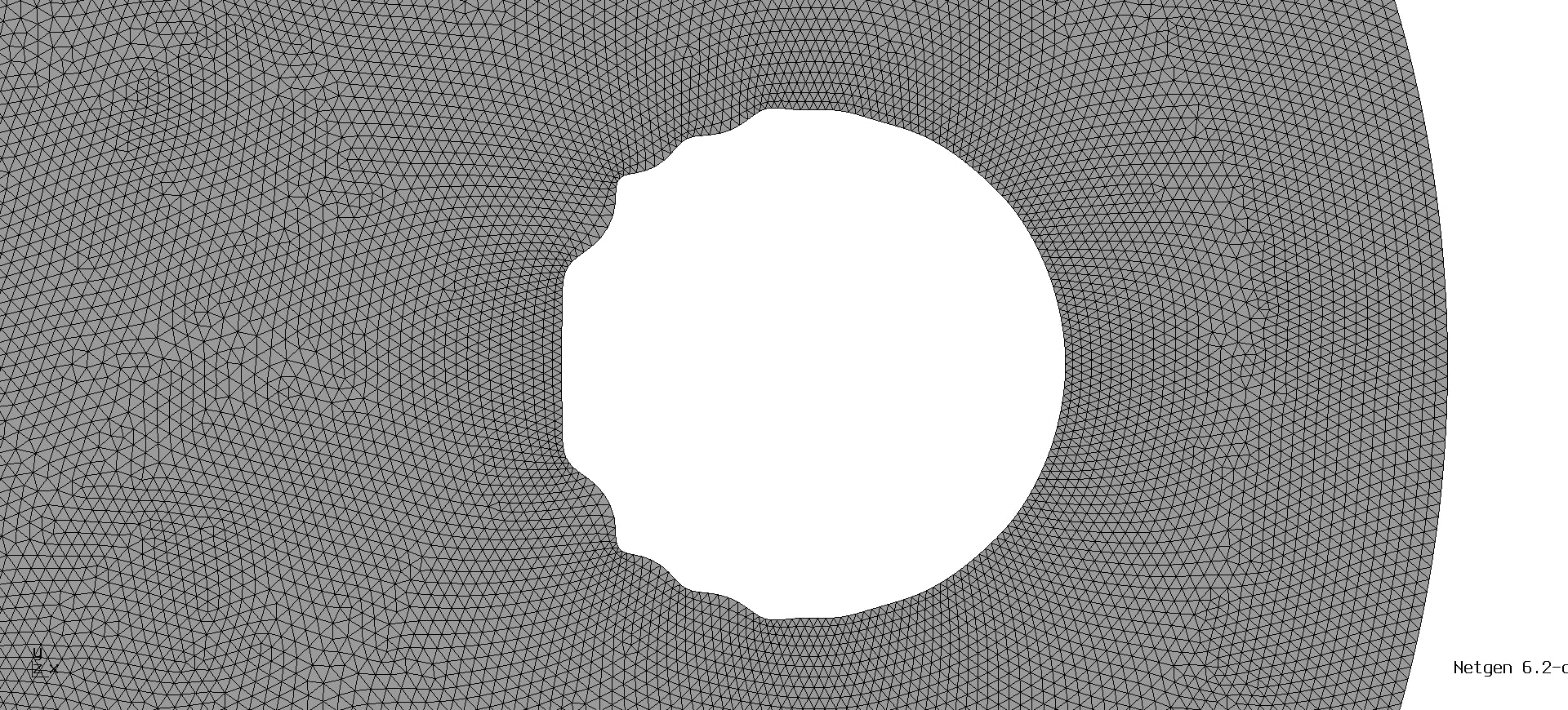} \\
            (a) & (b) & (c)
        \end{tabular}
    \caption{(a) Geometry with right hand side $f$. (b) Initial shape of scatterer (zoom of geometry in (a)). (c) Optimised shape of scatterer (zoom).}
    \label{fig_04hh_shapes}
 \end{figure}

 In this section, we consider the problem of finding the optimal shape of a scattering object. More precisely, we consider the minimisation of the functional
 \ben \label{eq_cost_helm}
    \int_{\Gamma_r} u \overbar u \; \mbox ds
\een
subject to the Helmholtz equation with impedance boundary conditions on the outer boundary: Find $u \in H^1(\Omega, \VC)$ such that 
\ben \label{eq_helmholtz}
    \int_\Omega [ \nabla u \cdot \nabla \bar w - \omega^2 u \bar w \big]
\, \mbox dx -
i \,\omega\, \int_{\Gamma} u \bar w \, \mbox ds = \int_{\Omega} f \bar w
\een
for all $w \in H^1(\Omega, \VC)$. Here, $\overbar w$ denotes the complex conjugate of a complex-valued function $w$, $\omega$ denotes the wave number, $i$ denotes the complex unit and the function $f$ on the right hand side is chosen as
\ben
    f(x_1, x_2) = 10^3 \cdot e^{-9((x_1-0.2)^2 + (x_2-0.5)^2)},
\een
see Figure \ref{fig_04hh_shapes}(a). Furthermore, $$\Omega = B((0.5,0.5)^\top,1) \setminus B((0.75, 0.5)^\top, 0.15)$$ denotes the domain of interest, $\Gamma = \{(x_1, x_2): x_1^2 + x_2^2 =1 \}$ the outer boundary and $\Gamma_r= \{(x_1, x_2): x_1^2 + x_2^2 =1, x_1 \geq 0 \}$ the right half of the outer boundary. Here, only the inner boundary $\partial \Omega \setminus \Gamma$ is subject to the shape optimisation. Thus, the aim of this model problem is to find a shape of the scattering object such that the waves are reflected away from $\Gamma_r$. 

\begin{figure}
        \begin{tabular}{cc}
            \includegraphics[width=.25\textwidth, trim = 160 0 160 0, clip]{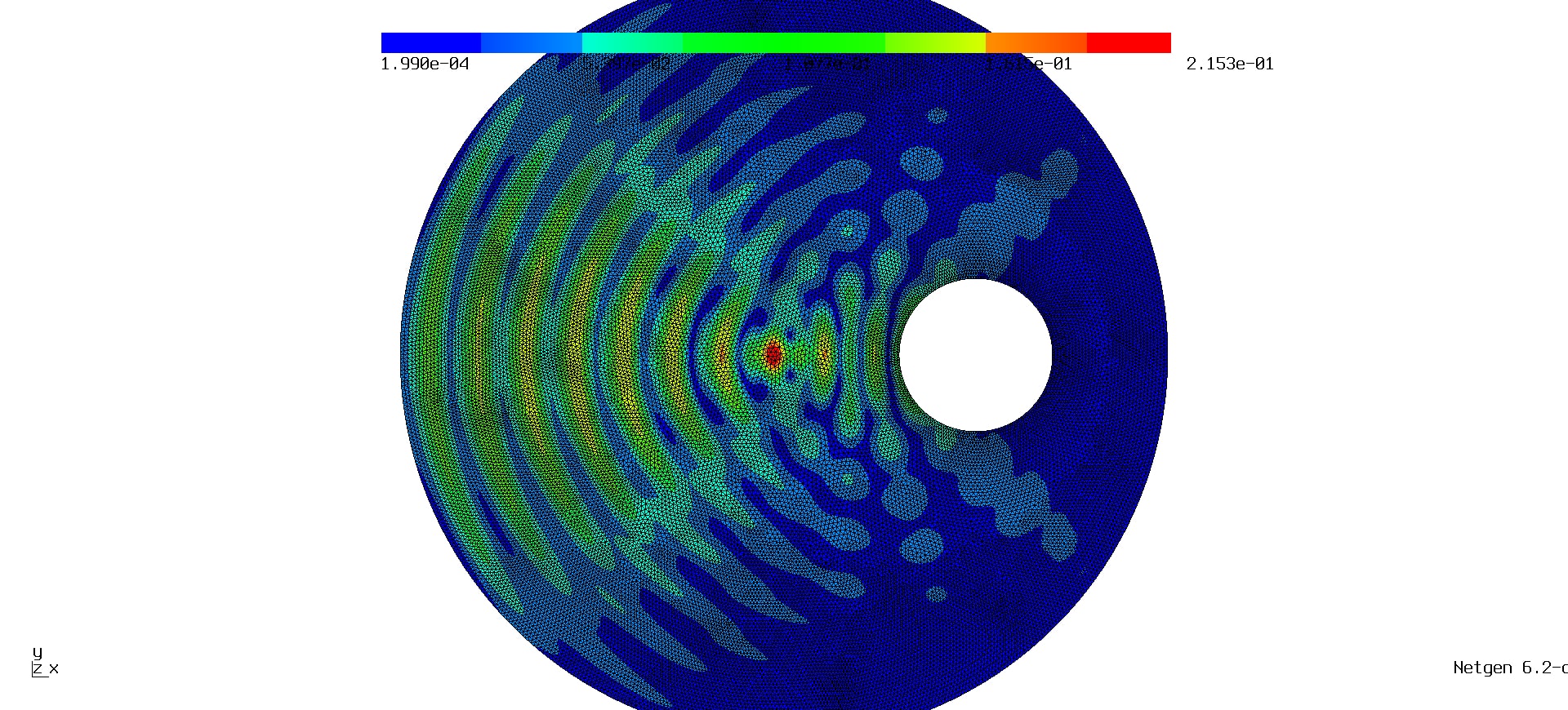} &
            \includegraphics[width=.25\textwidth, trim = 160 0 160 0, clip]{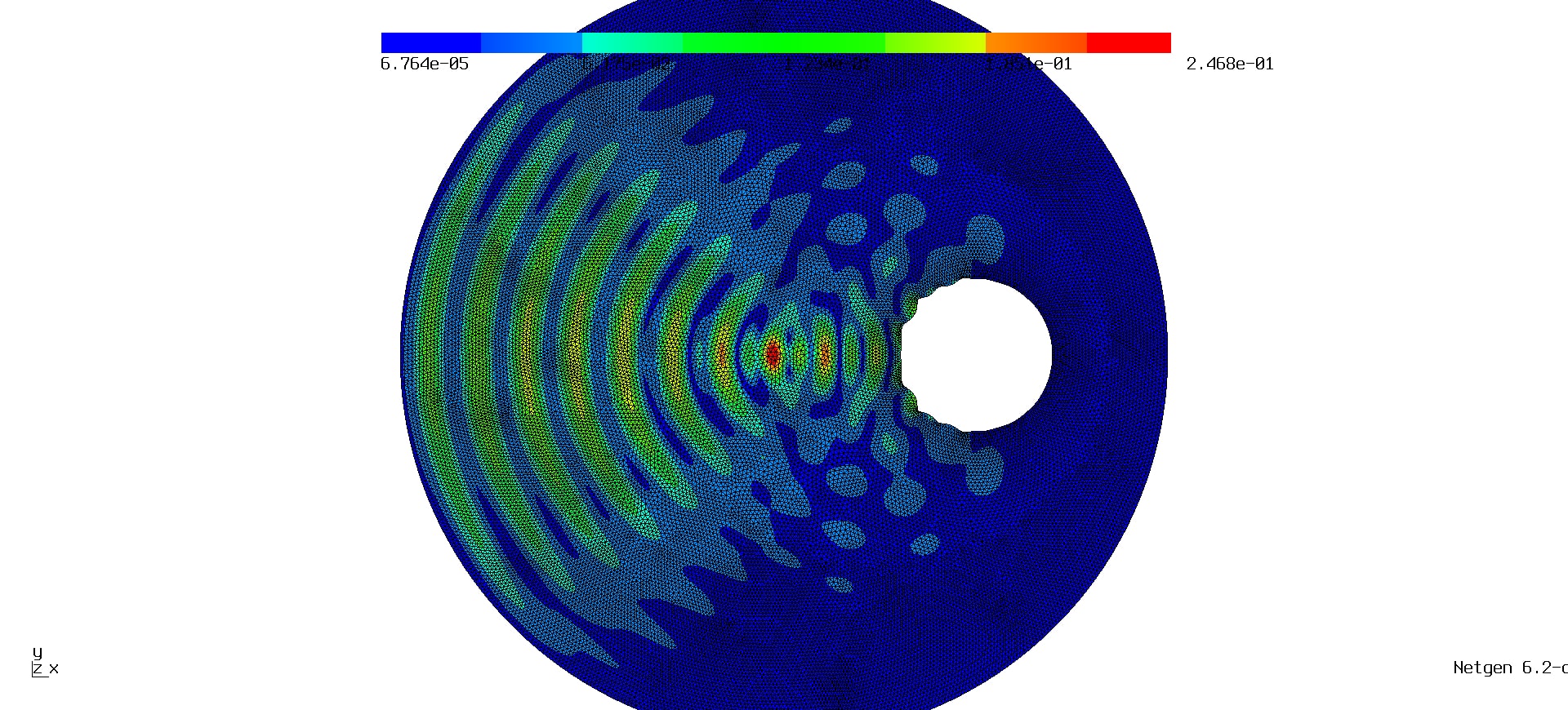} \\
            (a) & (b) 
        \end{tabular}
    \caption{(a) Absolute value of state $u$ for initial configuration. (b) Absolute value of state $u$ for optimised configuration.}
    \label{fig_04hh_states}
 \end{figure}

Figure \ref{fig_04hh_shapes} (b) and (c) show the initial and final shape of the scattering object, respectively. Figure \ref{fig_04hh_states} shows the norm of the state for the initial configuration (circular shape of scattering object) and for the optimised configuration. The objective value was reduced from $3.44 \cdot 10^{-3}$ to $3.31 \cdot 10^{-3}$. The forward simulations were performed using piecewise linear finite elements on a triangular grid consisting of with 34803 degrees of freedom. The optimisation stopped after 12 iterations.
  

 \subsection{Application to an Electric Machine}
\begin{figure}
        \begin{tabular}{c}
            \includegraphics[width=.47\textwidth]{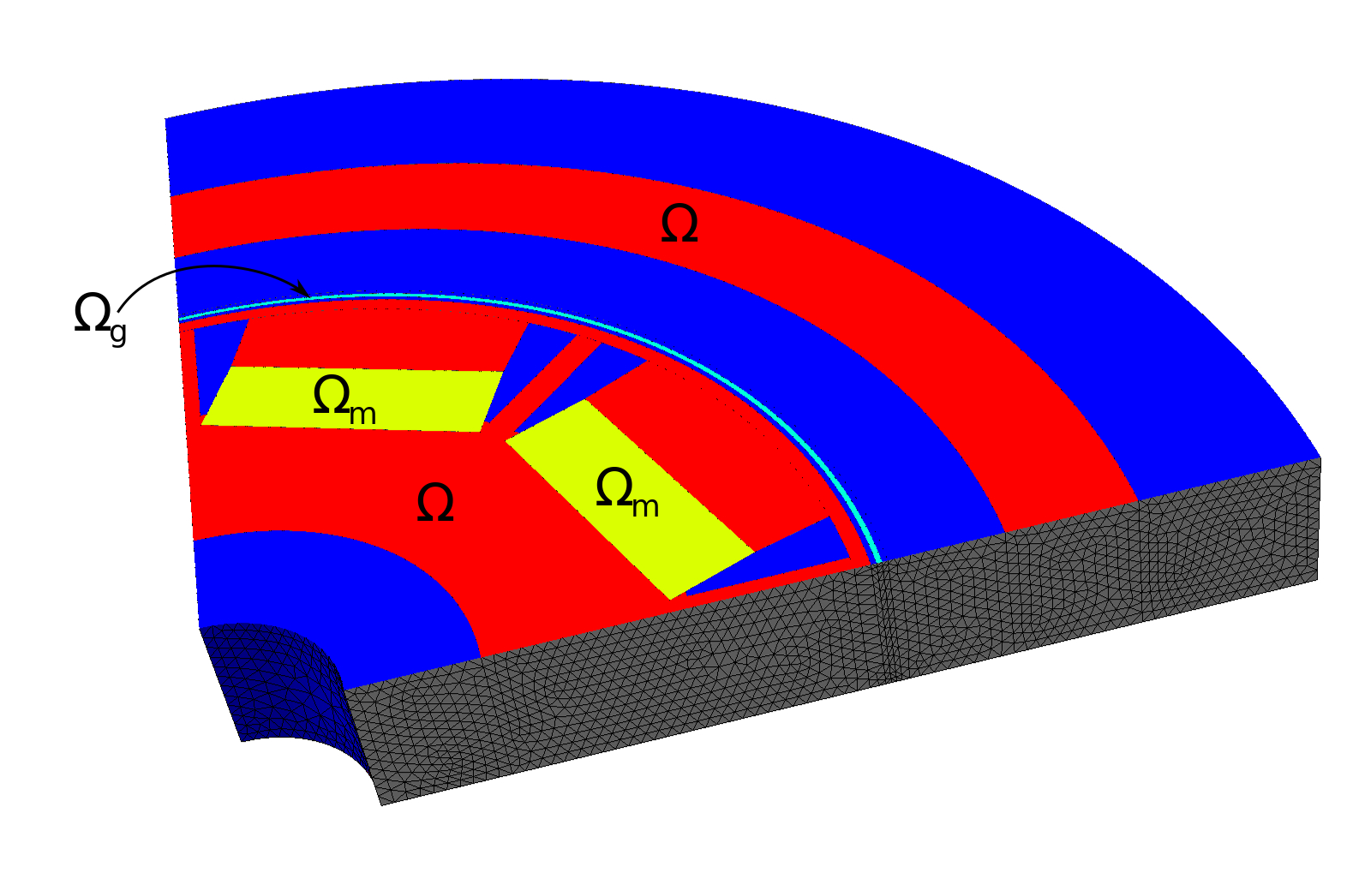}
        \end{tabular}
    \caption{Geometry of electric motor with subdomains in 2D cross section. The ferromagnetic subdomains $\Omega$ are depicted in red, $\Omega_m$ corresponds to the permanent magnets. The rest of the computational domain represents air. Further, $\Omega_g$ represents the air gap region that is relevant for the cost function \eqref{eq_costMotor}. }
    \label{fig_motorSubdom}
\end{figure}

 In this section, we consider the setting of three-dimensional nonlinear magnetostatics in $H(\curl,\Dsf)$ as it appears in the simulation of electric machines. Let $\Dsf \subset  \VR^3$ denote the computational domain, which consists of ferromagnetic material, air regions and permanent magnets, see Figure \ref{fig_motorSubdom}. Our aim is to minimise the functional
 \ben \label{eq_costMotor}
    \int_{\Omega_g} | \curl u \cdot \normal - B_d^n|^2 \; \mbox dx,
 \een
 where $\Omega_g$ denotes the air gap region of the machine, $\normal$ denotes an extension of the normal vector to the interior of $\Omega_g$, $B_d^n: \Omega_g \rightarrow \VR^3$ is a given smooth function and $u \in H_0(\curl,\Dsf)$ is the solution to the boundary value problem
 \ben \label{eq_magnetostatics}
    \int_D \nu_{\Omega}(|\curl u|) \curl u \cdot \curl w + \delta u \cdot w \; \mbox dx = \int_{\Omega_m} M \cdot \curl w \; \mbox dx
\een
for all $w \in H_0(\curl,\Dsf)$. Here, $\Omega \subset \Dsf$ denotes the union of the ferromagnetic parts of the electric machine, $\Omega_m$ denotes the permanent magnets subdomain and
\ben
    \nu_\Omega = \chi_\Omega(x) \hat \nu (|\curl u|) + \chi_{\Dsf \setminus \Omega}(x) \nu_0
\een
denotes the magnetic reluctivity, which is a nonlinear function $\hat \nu$ inside the ferromagnetic regions and equal to a constant $\nu_0$ elsewhere. Further, $\delta >0$ is a small regularisation parameter and $M:\Dsf \rightarrow \VR^3$ denotes the magnetisation in the permanent magnets. The nonlinear function $\hat \nu$ satisfies a Lipschitz condition and a strong monotonicity condition such that problem \eqref{eq_magnetostatics} is well-posed. The goal of minimising the cost function \eqref{eq_costMotor} is to obtain a design which exhibits a smooth rotation pattern. Note that in this particular example we do not consider rotation of the machine, but rather a fixed rotor position, and there are no electric currents present. We refer the reader to \cite[Sec. 6]{GanglSturm2019} for a more detailed description of the problem and to \cite{GLLMS2015} for a 2D version of the same problem.

As mentioned in Section \ref{sec_PDE_constr_shape}, the transformation $\Phi_t$ used in \eqref{eq_Gtup} depends on the differential operator. For the $\curl$-operator, we have 
\begin{align*}
\Phi_t(u) &= \partial T_t^{-\top} (u \circ T_t^{-1}) \qquad \qquad \mbox{ and } \\
(\curl (\Phi_t(u)))\circ T_t &= \frac{1}{\mbox{det}(\partial T_t)} \partial T_t \curl(u),
\end{align*}
see e.g. \cite[Section 3.9]{Monk2003}. Thus, the variational equation \eqref{eq_magnetostatics} can be defined as follows.
\begin{lstlisting}[firstnumber=last]
from math import pi
nu0 = 1e7 / (4*pi)
delta = 0.1

F = Id(3)
c1 = 1/Det(F) * F
c2 = Inv(F).trans

def EquationIron(u,w):
    return ( nuIron(Norm(c1*curl(u))) * (c1*curl(u))*(c1*curl(w)) + delta*(c2*u)*(c2*w) )*Det(F)*dx("iron") 
    
def EquationAir(u,w):
    return (nu0*(c1*curl(u))*(c1*curl(w)) + delta*(c2*u)*(c2*w))*Det(F)*dx("air|airgap") 

def EquationMagnets(u,w):
    return (nu0*(c1*curl(u))*(c1*curl(w)) + delta*(c2*u)*(c2*w)-magn*(c1*curl(w))) *Det(F)*dx("magnets") 
 
def Equation(u,w):
   return EquationIron(u,w) + EquationAir(u,w) + EquationMagnets(u,w)
\end{lstlisting}
Here, the computational domain consists of a subdomain representing the ferromagnetic part of the machine (\texttt{``iron''}) and a subdomain comprising the permanent magnets (\texttt{``magnets''}); the union of all air subdomains, including the air gap between rotating and fixed part of the machine, is given by \texttt{``air|airgap''}, see Figure \ref{fig_motorSubdom}.

\begin{figure}
        \begin{tabular}{cc}
            \includegraphics[width=.23\textwidth]{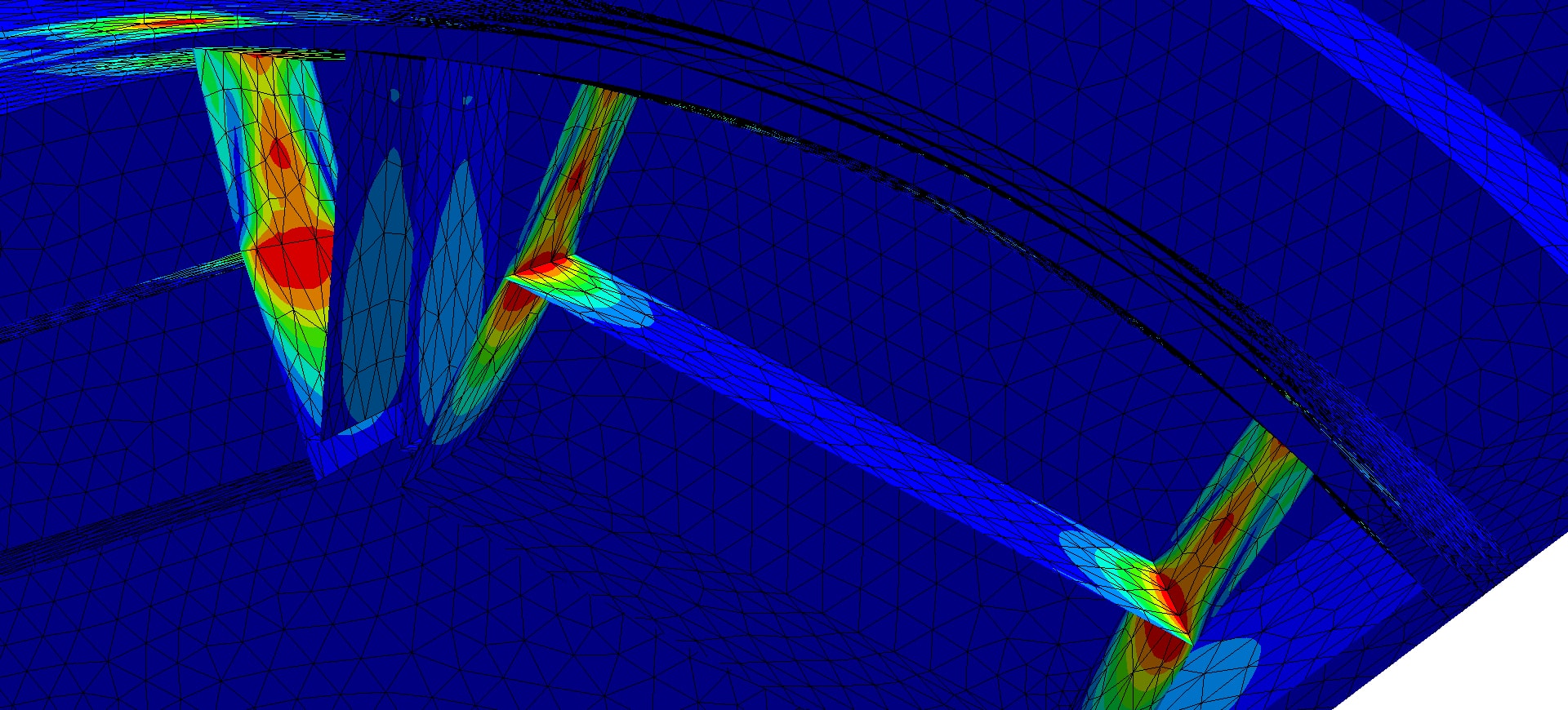} &
            \includegraphics[width=.23\textwidth]{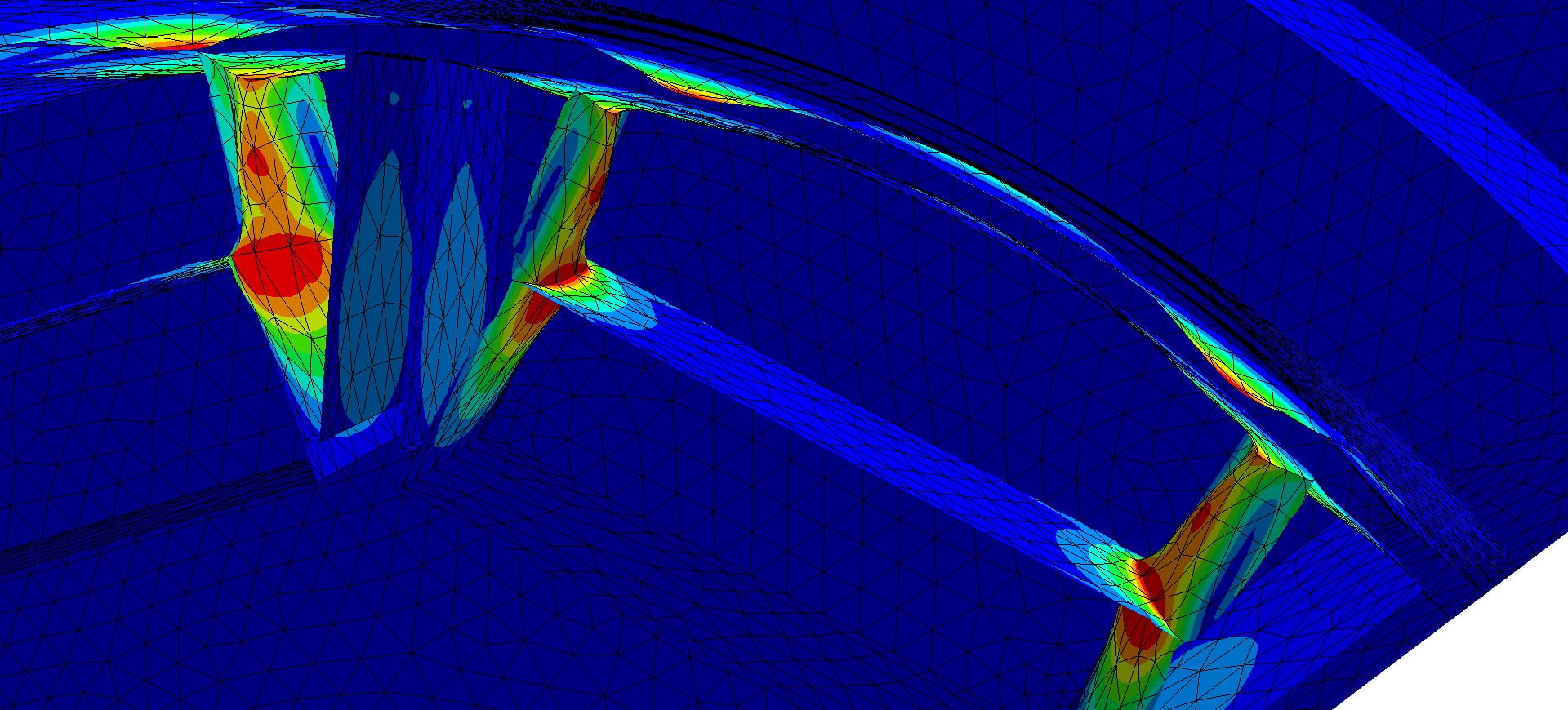} \\
            (a) & (b)
        \end{tabular}
    \caption{(a) Initial design of electric machine. (b) Optimised design.  }
    \label{fig_motor_v5}
 \end{figure}

\begin{figure}
        \begin{tabular}{ccc}
            \includegraphics[width=.14\textwidth, trim = 30 0 30 0, clip]{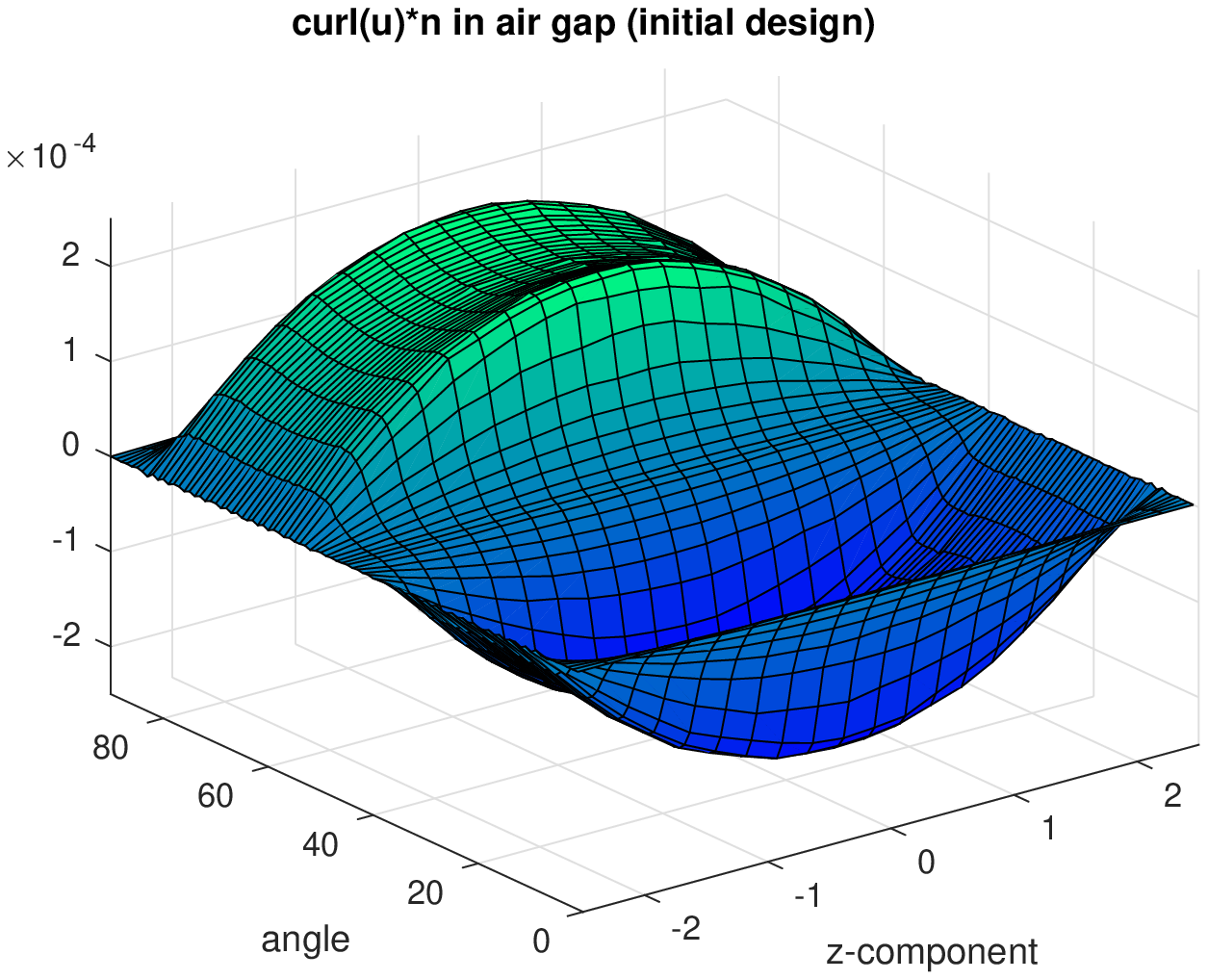} &
            \includegraphics[width=.14\textwidth, trim = 30 0 30 0, clip]{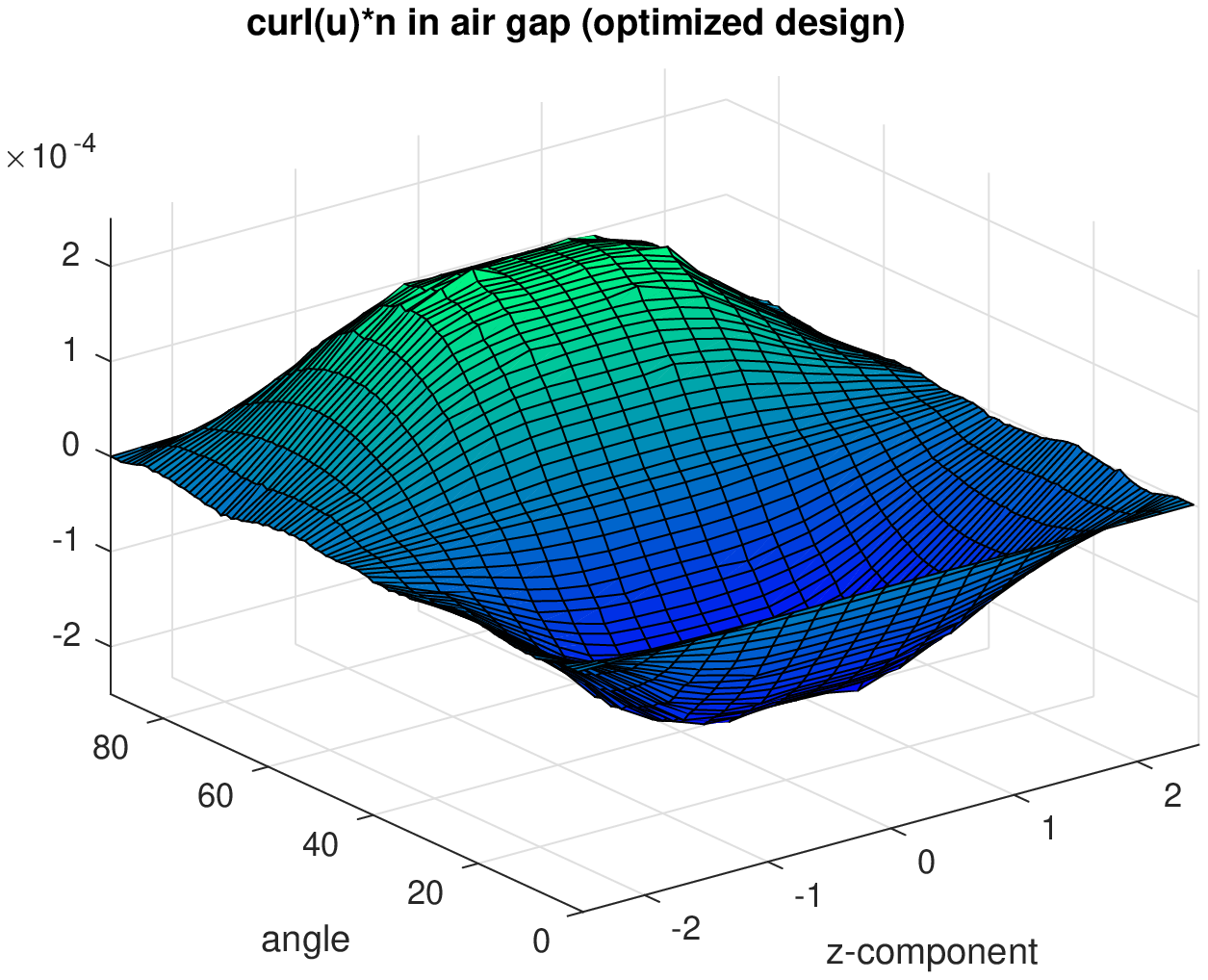} &
            \includegraphics[width=.14\textwidth, trim = 30 0 30 0, clip]{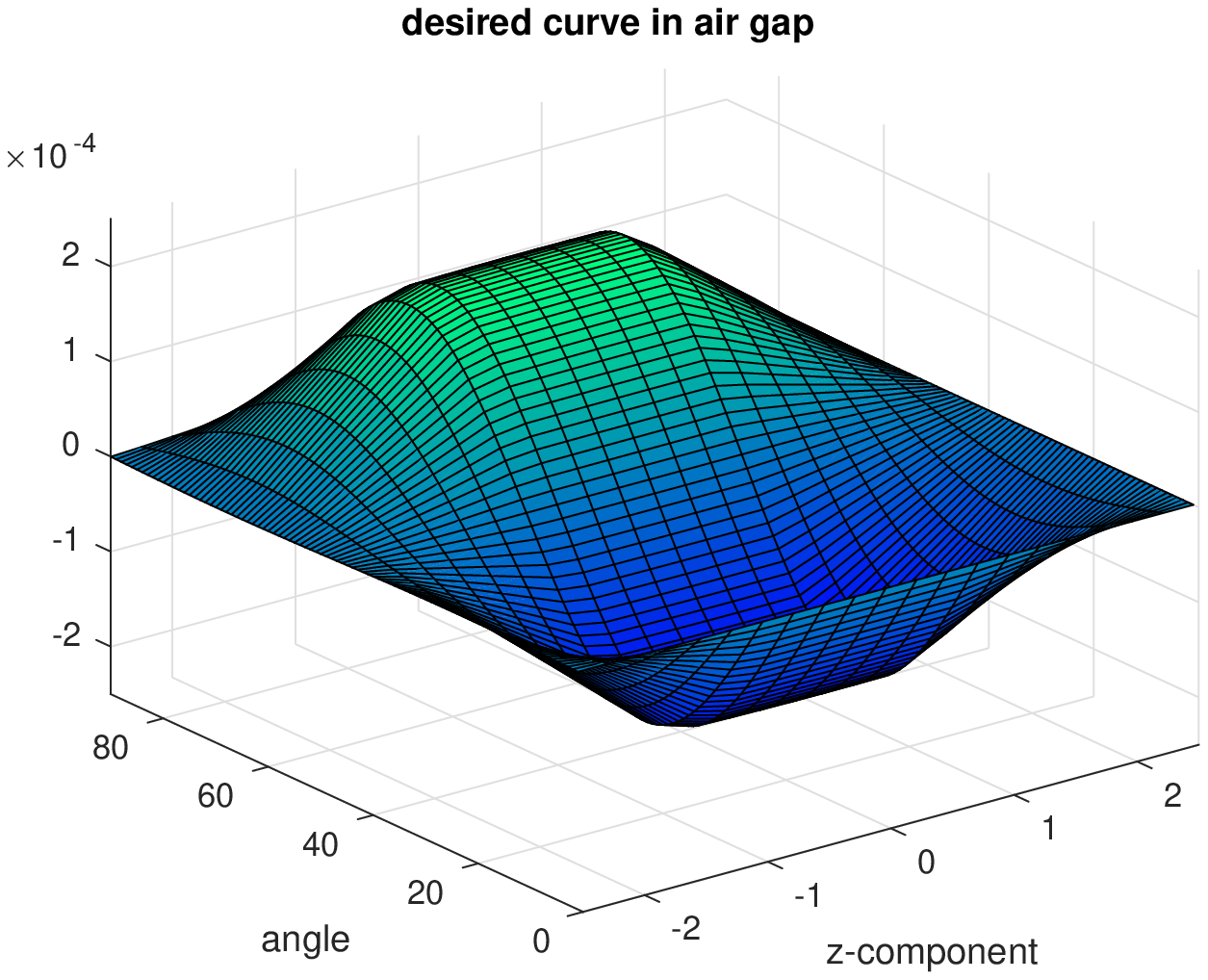}  \\
            (a) & (b) & (c)
        \end{tabular}
    \caption{Improvement of $\curl u \cdot n$ as a function of $z$ and the angle $\varphi$ for a fixed radius $r$ inside $\Omega_g$ compared to desired curve $B_d^n$. (a) $\curl u \cdot n$ for initial configuration. (b) $\curl u \cdot n$ for optimised configuration after 10 iterations. (c) Desired curve $B_d^n$ in polar coordinates as function of $z$ and the angle $\varphi$ for a fixed radius $r$.}
    \label{fig_curlun}
 \end{figure}

Moreover, \texttt{nuIron} denotes the nonlinear reluctivity function $\hat \nu$ and \texttt{magn} contains the magnetisation direction of the permanent magnets. Likewise, the objective function can be defined as follows,
\begin{lstlisting}[firstnumber=last]
def Cost(u):
    return (InnerProduct(c1*curl(u),n2D) - Bd)*Det(F)*dx("airgap")
\end{lstlisting}
where \texttt{n2D} and \texttt{Bd} represent the extension of the normal vector to the interior of the air gap and the desired curve, respectively. For the definition of all quantities, we refer the reader to Online Resource 1.
The shape differentiation as well as the optimisation loop now works in the same way as in the previous examples. Figure \ref{fig_motor_v5} shows the initial design of the motor as well as the optimised design obtained after 11 iterations of Algorithm \ref{alg:gradient} with $\gamma = 0$. The experiment was conducted using a tetrahedral finite element mesh consisting of 13440 vertices, 57806 elements and N\'ed\'elec elements of order 2 (resulting in a total of 323808 degrees of freedom).
The objective value was reduced from $2.5944 \cdot 10^{-8}$ to $4.565 \cdot 10^{-10}$ in the course of the first order optimisation algorithm after 11 iterations. It can be seen from Figure \ref{fig_curlun} that the difference between the quantity $\curl(u)\cdot n$ and the desired curve $B_d^n$ inside the air gap decreases significantly.

 \subsection{Surface PDEs}
\begin{figure}
    \begin{tabular}{cc}
        \includegraphics[width=.22\textwidth]{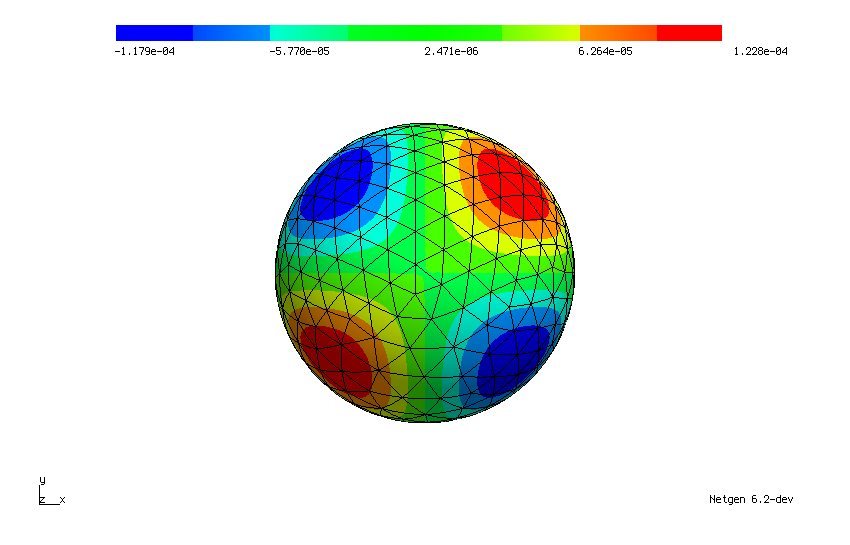}& \includegraphics[width=.22\textwidth]{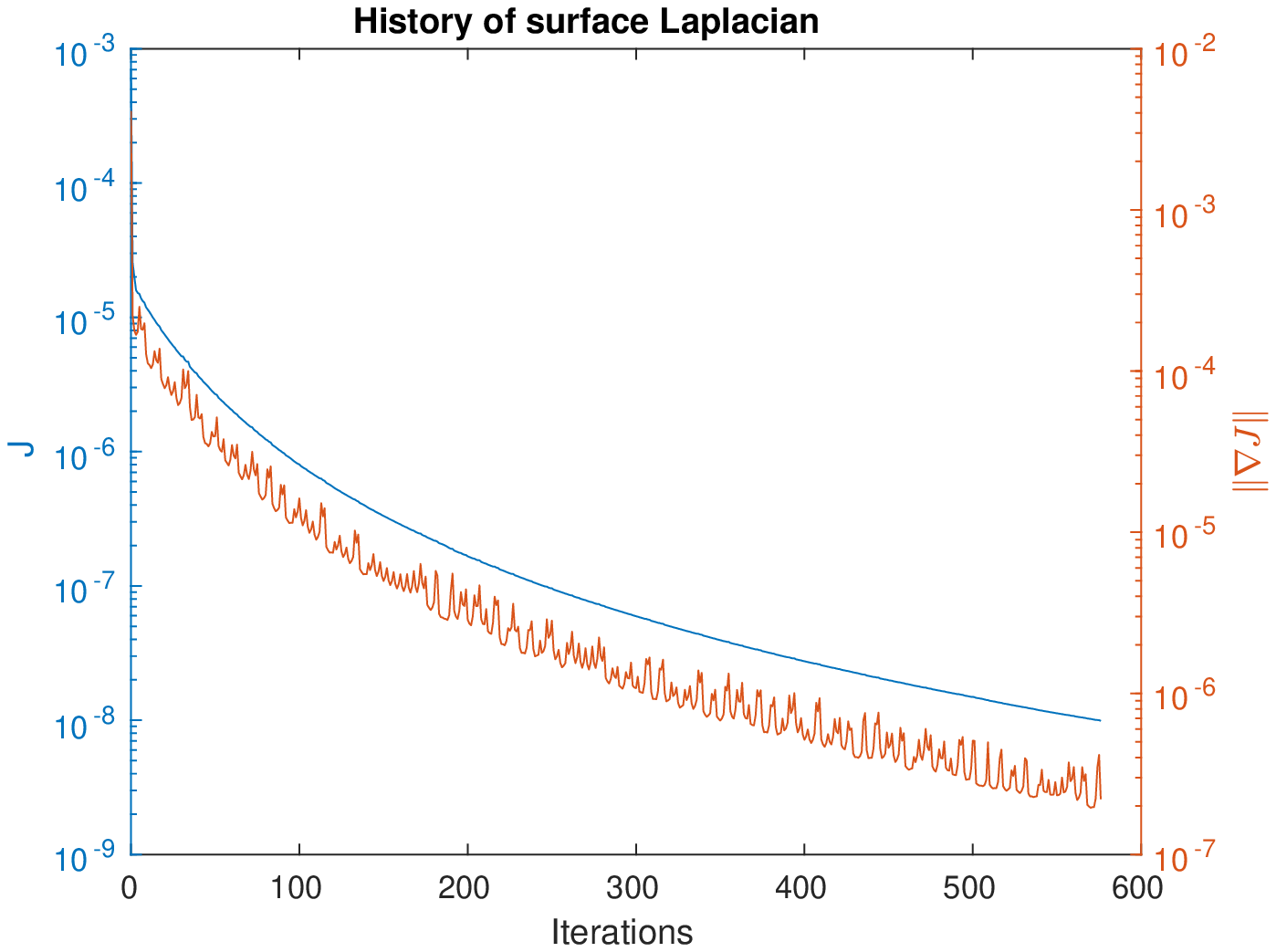}
         \\
        (a)&(b)
    \end{tabular}
    \caption{(a) Initial geometry for shape optimisation with respect to surface PDE \eqref{eq_surface_J}--\eqref{eq_surface_pde}. (b) History of objective value and norm of shape gradient using a first order algorithm with line search.}
    \label{fig_surface1}
 \end{figure}

\begin{figure}
    \begin{tabular}{cc}
        \includegraphics[width=.25\textwidth]{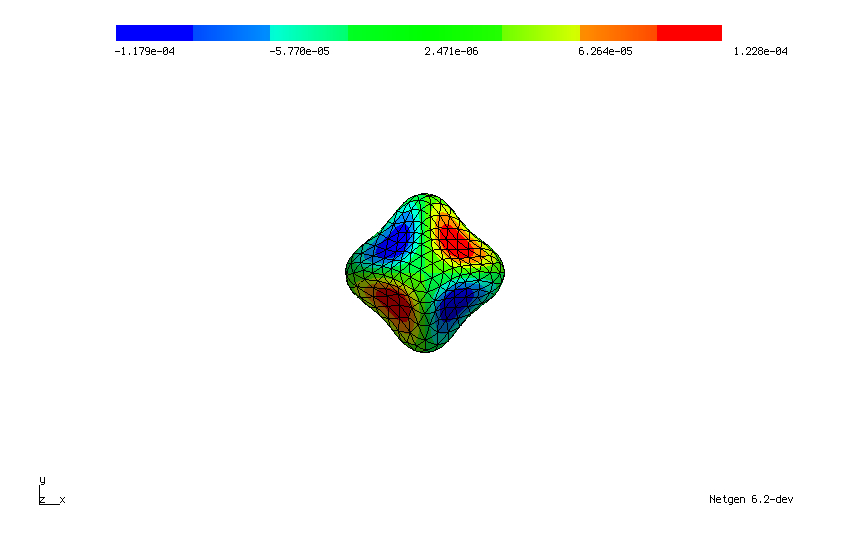}&
        \includegraphics[width=.25\textwidth]{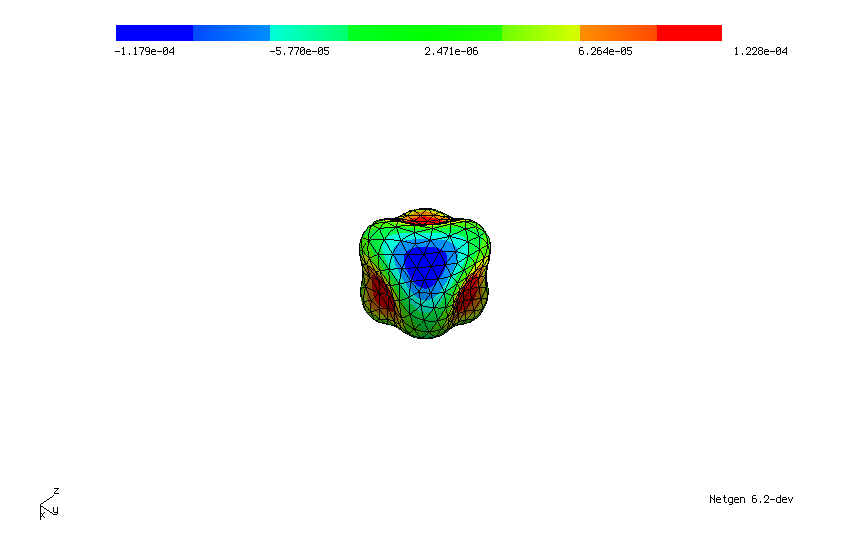}\\
        (a) & (b)        
    \end{tabular}
    \caption{(a) Final design after 575 iterations. (b) Different view of (a).}
    \label{fig_surface2}
 \end{figure}

 Next, we also show the application of a shape optimisation algorithm to a problem constrained by a surface PDE. We solve problem \eqref{eq_surface_J}--\eqref{eq_surface_pde} with $u_d = 0$, $f(x_1, x_2, x_3) = x_1 x_2 x_3$ and initial shape $M = S^2$ the unit sphere in $\VR^3$.  
 We applied a first order algorithm with a line search. Figure \ref{fig_surface1} shows the initial geometry as well as the decrease of the objective function and of the norm of the shape gradient. The objective value was reduced from $7.08 \cdot 10^{-4}$ to $9.88 \cdot 10^{-9}$. Figure \ref{fig_surface2} shows the final design which was obtained after 575 iterations from two different perspectives. The experiment was conducted using a surface mesh with 332 vertices and 660 faces and polynomial degree 3 (resulting in 2972 degrees of freedom).

 \subsection{Time-dependent PDE using space-time method} \label{sec_spacetime}
 In this section we illustrate a non-standard situation where the fully automated shape differentiation using the command $\texttt{.DiffShape()}$ fails, but the semi-automated way can be used to compute the shape derivative.
 
 The situation is that of a parabolic PDE constraint in a space-time setting where the time variable is considered as just another space variable. Let $T>0$ and $\Omega \subset \VR^d$ be given and define the space-time cylinder $Q := \Omega \times (0,T) \subset \VR^{d+1}$. For given smooth functions $u_d$ and $f$ defined on $Q$, we consider the problem
 \begin{subequations} \label{eq_optiHeatEq}
 \begin{align}
    \underset{\Omega}{\mbox{min }} &\int_Q |u-u_d|^2 \, \mbox d (x,\tau) \\
    \mbox{ s.t. } \int_Q \partial_\tau u v &+ \nabla_x u \cdot \nabla_x v \, \mbox d(x,\tau) = \int_Q f v \, \mbox d (x,\tau) \label{eq_optiHeatEq_pde}
 \end{align}
 \end{subequations}
 for all $v$ in the Bochner space $L^2(0,T; H^1_0(\Omega))$ where the state $u$ is to be sought in the Bochner space $L^2(0,T; H^1_0(\Omega)) \cap H^1(0,T; H^{-1}(\Omega))$ with the initial condition $u(x,0) = 0$. Here, $\nabla_x = (\partial_{x_1}, \dots, \partial_{x_d})^\top$ denotes the spatial gradient and $\partial_\tau$ the time derivative. Note that we denote the time variable by $\tau$ in order not to interfere with the shape parameter $t$. We refer the interested reader to \cite{Steinbach2015} for details on this space-time formulation of the PDE constraint. As it is outlined there, the PDE can be solved numerically by choosing the same ansatz and test space consisting of piecewise linear and globally continuous finite element functions on $Q$.
 
 For simplicity, we restrict ourselves to the case where $d=1$, i.e. to the case where the spatial domain $\Omega$ is an interval. We are interested in the shape derivative of problem \eqref{eq_optiHeatEq} with respect to spatial perturbations, i.e., with respect to transformations of the form
 \begin{align*}
    T_t(x,\tau) = \left( \begin{array}{c} x + t V(x, \tau) \\ \tau \end{array} \right)
 \end{align*}
 where $V \in C^{0,1}(Q, \VR^d)$ and $t \geq 0$. We recall the notation $F_t(x, \tau) = \partial T_t(x, \tau)$. By this choice of transformation $T_t$ we exclude an unwanted deformation of the space-time cylinder into the time direction as the time horizon $T>0$ is assumed to be fixed.
 
Following the lines of the previous examples, we can define the cost function, the PDE and the Lagrangian:
\begin{lstlisting}[firstnumber = last]
F11 = Parameter(1)
F12 = Parameter(0)
F21 = Parameter(0)
F22 = Parameter(1)
F = CoefficientFunction((F11, F12, F21, F22), dims=(2,2))

def Cost_vol(u): 
    return (u-ud)*(u-ud)*Det(F)*dx

def Equation_vol(u,v):
   return ( (Inv(F.trans)*grad(u))[1]*v + (Inv(F.trans)*grad(u))[0]*(Inv(F.trans)*grad(v))[0]-f*v)*Det(F)*dx
   
G_vol = Cost_vol(gfu) + Equation_vol(gfu,gfp)
\end{lstlisting}
Here, \texttt{gfu} denotes the solution to the state equation \eqref{eq_optiHeatEq_pde} and \texttt{gfp} the solution to the adjoint equation, which is posed backward in time and reads in its strong form,
\begin{align*}
    -\partial_\tau p - \Delta p &= - 2(u-u_d) && (x,\tau) \in Q, \\
    p(x, \tau) &= 0 && (x,\tau) \in \partial \Omega \times (0,T),\\
    p(x,T) &= 0 && x \in \Omega.
\end{align*}
We can compute the shape derivative similarly to the previous examples by means of formula \eqref{eq_dGdTdGdF}, i.e.,
\begin{align}
     \left. \frac{d}{dt} \mathcal J(\Omega_t)  \right|_{t=0}
   &= \left. \left( \frac{dG}{dT_t} \frac{d T_t}{dt} + \frac{dG}{dF_t} \frac{d F_t}{dt}\right)\right|_{t=0}. \label{eq_dGdTdGdF2}
\end{align}
However, it must be noted that in this special situation we have
\begin{align}
    \frac{d T_t}{dt} = \left( \begin{array}{c} V \\ 0 \end{array} \right) \; \mbox{ and } \; \frac{d F_t}{dt} = \left( \begin{array}{cc} \partial_x V & \partial_\tau V\\ 0&0 \end{array} \right). \label{eq_dTdtdFdt}
\end{align}
The shape derivative can now be obtained as follows: Given a mesh of the space-time cylinder $Q$, we define an $\VR^d$-valued $H^1$-space to represent the vector field $V$ (here we assumed $d=1$, thus we are facing a scalar-valued space). The shape derivative is a linear functional on this space and is obtained by plugging in \eqref{eq_dTdtdFdt} into \eqref{eq_dGdTdGdF2}:
\begin{lstlisting}[firstnumber = last]   
VEC1 = H1(mesh, order=1)    
V = VEC1.TestFunction()

dJOmega_para = LinearForm(VEC1);(*@\label{lst_dJOmega_para1}@*)
dJOmega_para += G_vol.Diff(x, V)
dJOmega_para += G_vol.Diff(F11, grad(V)[0])
dJOmega_para += G_vol.Diff(F12, grad(V)[1])(*@\label{lst_dJOmega_para4}@*)
\end{lstlisting}
 
 \begin{remark}
    The fully automated shape differentiation command \texttt{.DiffShape(V)} cannot be used here because the vector field $V$ has less components than the dimension of the mesh. On the other hand, if $V$ was chosen as a vector field with $d+1$ components, the command \texttt{.DiffShape(V)} would evaluate formula \eqref{eq_dGdTdGdF2}, but would assume $\frac{d T_t}{dt} = V = (V_x, V_\tau)^\top$ and $$\frac{d F_t}{dt} = \partial V = \left( \begin{array}{cc} \partial_x V_x & \partial_\tau V_x\\ \partial_x V_\tau & \partial_\tau V_\tau \end{array} \right)$$ and could not take into account the special situation at hand as shown in \eqref{eq_dTdtdFdt}. This example is meant to illustrate the greater flexibility of the semi-automated compared to the fully automated shape differentiation.
 \end{remark}

 Code lines \ref{lst_dJOmega_para1}-\ref{lst_dJOmega_para4} show the computation of the shape derivative in the direction of an $\VR^d$-valued function $V = V(x,\tau)$ (recall $d=1$ here). However, using $\tau$-dependent vector fields would result in time-dependent optimal shapes, which is often not desired. Rather, one is interested in vector fields of the form $V = V(x) \neq V(x,\tau)$ which still yield a descent, i.e. $D \mathcal J(\Omega)( V) < 0$. This can be achieved as follows:
 \begin{enumerate}
  \item Compute a time-dependent shape gradient $\tilde W$ by solving
  \begin{align}\label{eq_avg_step1}
    \int_Q \partial \tilde W : \partial \tilde V + \tilde W \cdot \tilde V = D\mathcal J(\Omega)( \tilde V) \quad \mbox{for all } \tilde V,
  \end{align}
  \item Set $W(x, \tau) = \frac{1}{T} \int_0^T \tilde W(x, s) \mbox ds$.
 \end{enumerate}
 Note that $W(x,\tau)$ is constant in $\tau$. Then we see by plugging in $\tilde V = -W$ in \eqref{eq_avg_step1} that
 \begin{align*}
    D \mathcal J(\Omega)(-W)  <0,
 \end{align*}
 thus $-W$ is a descent direction.
 
 We used this strategy to solve problem \eqref{eq_optiHeatEq} for $d=1$ with the data $u_d(x, \tau) = x(1-x)\tau$, $f(x, \tau ) = x(1-x)+2\tau$ numerically starting out from the initial domain $\Omega_{init} = (0.2, 0.8)$ and the fixed time interval $(0,T) = (0,1)$. Note that the data is chosen such that the domain $\Omega^\star = (0,1)$ is a global solution to problem \eqref{eq_optiHeatEq}.
 
\begin{figure}
    \begin{center}
        \begin{tabular}{cc}
            \includegraphics[width=.22\textwidth, trim = 100 0 60 0, clip]{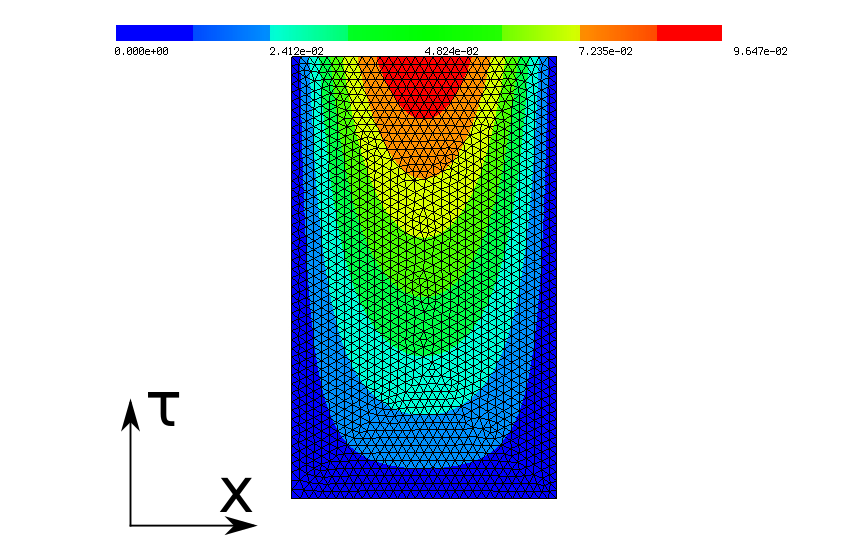} & \includegraphics[width=.22\textwidth, trim = 100 0 60 0, clip]{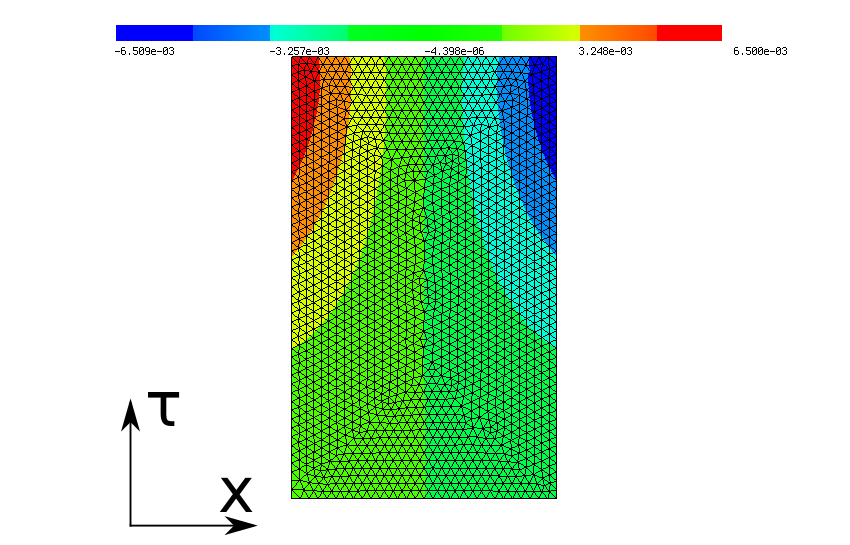} \\
            (a) & (b)
        \end{tabular}
    \end{center}
    \caption{Space-time cylinder in initial configuration. (a) Solution to state equation \eqref{eq_optiHeatEq_pde}. (b) Time-dependent descent vector field $\tilde W$ obtained by solving \eqref{eq_avg_step1}.}
    \label{fig_spacetime1}
 \end{figure}

\begin{figure}
    \begin{center}
        \begin{tabular}{cc}
            \includegraphics[width=.22\textwidth, trim = 100 0 60 0, clip]{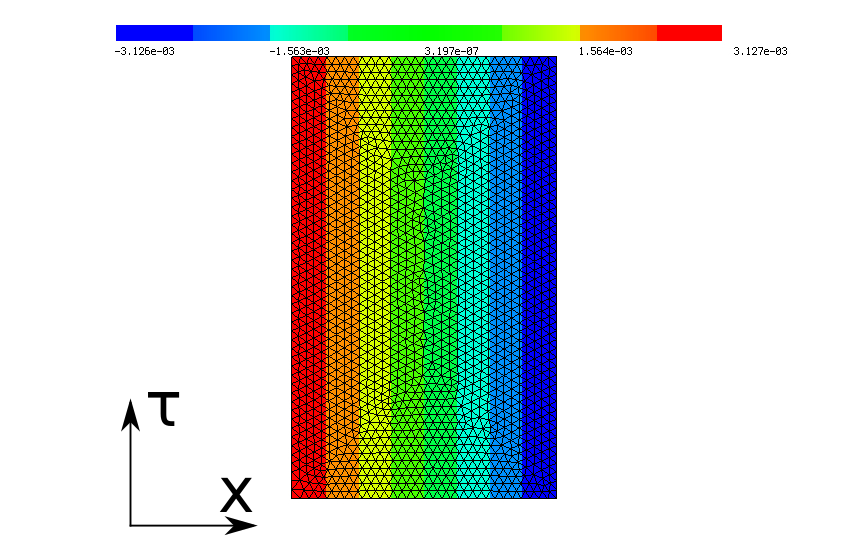} & \includegraphics[width=.22\textwidth, trim = 100 0 60 0, clip]{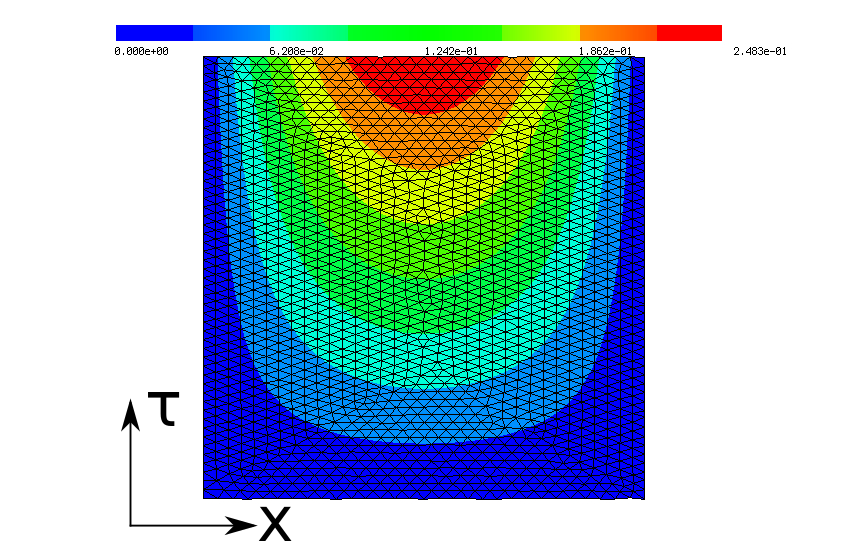} \\
            (a) & (b)
        \end{tabular}
    \end{center}
    \caption{Space-time cylinder in initial and final configuration. (a) Time-independent descent vector field $W$ obtained by averaging $\tilde W$. (b) Solution to state equation \eqref{eq_optiHeatEq_pde} on final design.}
    \label{fig_spacetime2}
 \end{figure}

 Figure \ref{fig_spacetime1} shows the initial design together with the solution to the state equation and the time-dependent descent vector field $\tilde W$ obtained as solution of \eqref{eq_avg_step1}. Figure \ref{fig_spacetime2}(a) shows the averaged vector field $W$ which is independent of $\tau$ and yields a uniform deformation of the space-time cylinder. The final design after 293 iterations can be seen in Figure \ref{fig_spacetime2}(b). The cost function was reduced from $4.65 \cdot 10^{-3}$ to $9.95 \cdot 10^{-9}$.
 

  For more details on the implementation of this example, we refer to Online Resource 1 and for more details on shape optimisation in a space-time setting, we refer the interested reader to \cite{master_CK}.

 \section*{Conclusion}
We showed how to obtain first and second order shape derivatives for unconstrained and PDE-constrained shape optimisation problems in a semi-automatic and fully automatic way in the finite element software package \ngsolve{}. We verified the proposed method numerically by Taylor tests and by showing its successful application to several shape optimisation problems. We believe that this intuitive approach can help research scientists working in the field of shape optimisation to further improve numerical methods on the one hand, and product engineers working with \ngsolve{} to design devices in an optimal fashion on the other hand.

\section*{Acknowledgements}
Michael Neunteufel has been funded by the Austrian Science Fund (FWF) project W1245. Moreover, we would like to thank Christian K\"othe for his contribution to Section \ref{sec_spacetime}.

\section*{Replication of results}
The python scripts which were used for the results presented in this paper are available in Online Resource 1. All computations were performed using NGSolve version V6.2.2004.

\section*{Conflicts of interest/Competing interests}
On behalf of all authors, the corresponding author states that there is no conflict of interest.

\bibliography{template_rev2}
\bibliographystyle{plain}

\end{document}